\documentclass[11pt]{amsart}

\newif\ifcomments
\commentstrue 

\makeatletter
\AtBeginDocument{%

  \let\original@@tocwrite=\@tocwrite
  \newif\ifHKVflag
  \def\HKV@uniqtoken{\HKV@uniqtoken}
  \def\HKV@endmark{\HKV@endmark}
  \long\def\HKV@getfirsttoken#1#{\HKV@getfirsttoken@#1\bgroup\HKV@endmark}
  \long\def\HKV@getfirsttoken@#1#2\HKV@endmark#3\HKV@endmark{#1}
  \renewcommand{\@tocwrite}[2]{%
    \begingroup
      \HKVflagfalse 
      \@ifempty{#2}{}{%
        \expandafter\expandafter\expandafter\ifx\HKV@getfirsttoken#2\HKV@uniqtoken{}\HKV@endmark\Sectionformat\expandafter\@firstoftwo\else\expandafter\@secondoftwo\fi
        {%
          \def\Sectionformat##1##2{\@ifempty{##1}{}{\HKVflagtrue}}%
          #2
        }{\HKVflagtrue}%
      }%
      \def\@tempa{}%
      \ifHKVflag\def\@tempa{\original@@tocwrite{#1}{#2}}\fi
    \expandafter\endgroup
    \@tempa
  }%
}

\let\original@maketitle=\maketitle
\let\origina@@maketitle=\@maketitle
\let\original@email=\email
\let\original@address=\address
\expandafter\def\expandafter\appendix\expandafter{\expandafter\enddoc@text\appendix
\let\addresses\@empty
\let\email=\original@email
\let\address=\original@address
\bigskip
}
\makeatother

\usepackage[left=1.4in,top=1.3in,right=1.4in,nofoot]{geometry}
\usepackage{amssymb}
\usepackage[all,cmtip,poly]{xy}
\usepackage[usenames,dvipsnames]{xcolor}
\usepackage{hyperref}
\usepackage{enumitem}
\usepackage{mathrsfs}
\usepackage{tikz-cd}
\usepackage{dsfont}
\usepackage{mathtools}

\newcommand{\WD}{{\operatorname{WD}}}
\newcommand{\rec}{{\operatorname{rec}}}

\newcommand{\m}{\mathfrak{m}}
\newcommand{\n}{\mathfrak{n}}
\newcommand{\p}{\mathfrak{p}}
\newcommand{\B}{\mathfrak{B}}
\newcommand{\fU}{\mathfrak{U}}

\newcommand{\ang}[1]{\langle #1 \rangle}

\newcommand{\T}{\mathbb{T}}

\newcommand{\cA}{\mathcal{A}}

\newcommand{\cC}{\mathcal{C}}
\newcommand{\cD}{\mathcal{D}}

\newcommand{\cN}{\mathcal{N}}
\newcommand{\cO}{\mathcal{O}}
\newcommand{\cP}{\mathcal{P}}

\newcommand{\cW}{\mathcal{W}}

\newcommand{\cZ}{\mathcal{Z}}

\newcommand{\bG}{\mathbf{G}}
\newcommand{\bH}{\mathbf{H}}

\newcommand{\bZ}{\mathbf{Z}}

\newcommand{\bPGL}{\mathbf{PGL}}
\newcommand{\bPU}{\mathbf{PU}}
\newcommand{\bPSO}{\mathbf{PSO}}

\newcommand{\A}{\mathbb A}
\newcommand{\C}{\mathbb C}
\newcommand{\FF}{\mathbb F}
\newcommand{\F}{\wt F}
\newcommand{\fp}{\FF_p}

\newcommand{\fpb}{\overline \FF_p} 

\newcommand{\HH}{\mathcal H}

\renewcommand{\O}{{\mathcal O}}

\newcommand{\Q}{\mathbb Q}

\newcommand{\qp}{{\Q_p}}

\newcommand{\qpb}{\overline\Q_p}

\newcommand{\Z}{\mathbb Z}

\newcommand{\zpb}{\overline\Z_p}

\newcommand{\eq}{\Leftrightarrow}
\newcommand{\into}{\hookrightarrow}

\newcommand{\onto}{\twoheadrightarrow}
\newcommand{\congto}{\xrightarrow{\,\sim\,}}

\newcommand{\tocong}{\xleftarrow{\,\sim\,}}
\newcommand{\vp}{\varphi}
\newcommand{\ve}{\varepsilon}

\newcommand{\s}{^\times}

\newcommand{\Fss}{^{\mathrm{F-ss}}}
\newcommand{\op}{^{\mathrm{op}}}
\newcommand{\ab}{^{\mathrm{ab}}}
\newcommand{\cts}{_{\mathrm{cts}}}

\newcommand{\dual}{^\vee}

\newcommand{\subalg}{_{\mathrm{alg}}}

\newcommand{\PU}{\operatorname{PU}}
\newcommand{\GL}{\operatorname{GL}}

\newcommand{\PGL}{\operatorname{PGL}} 
\newcommand{\HT}{\operatorname{HT}}

\newcommand{\Res}{\operatorname{Res}}

\DeclareMathOperator{\Mod}{Mod}
\DeclareMathOperator{\Rep}{Rep}

\DeclareMathOperator{\val}{val}

\DeclareMathOperator{\Aut}{Aut}
\DeclareMathOperator{\End}{End}
\DeclareMathOperator{\Hom}{Hom}

\DeclareMathOperator{\diag}{diag}

\DeclareMathOperator{\Gal}{Gal}

\DeclareMathOperator{\tr}{tr}

\DeclareMathOperator{\BC}{BC}

\DeclareMathOperator{\Ind}{Ind}

\DeclareMathOperator{\Frob}{Frob}

\DeclareMathOperator{\LJ}{LJ}
\DeclareMathOperator{\JL}{JL}
\DeclareMathOperator{\Nrd}{Nrd}
\DeclareMathOperator{\Sp}{Sp}

\DeclareMathOperator{\chr}{char}
\DeclareMathOperator{\lcm}{lcm}
\DeclareMathOperator{\vol}{vol}
\DeclareMathOperator{\Tor}{Tor}
  \newcommand{\ilim}{\varinjlim}   

\renewcommand{\o}[1]{\overline{#1}}
\renewcommand{\u}[1]{\underline{#1}}
\newcommand{\wt}[1]{\widetilde{#1}}
\newcommand{\wh}[1]{\widehat{#1}}

\newcommand{\ad}{^\mathrm{ad}}
\renewcommand{\sc}{^\mathrm{sc}}

\newcommand{\crs}{^{\mathrm{cr}}}
\newcommand{\st}{\mathrm{st}}

\newcommand{\textA}{\mathrm{A}}
\newcommand{\textB}{\mathrm{B}}
\newcommand{\textC}{\mathrm{C}}
\newcommand{\textD}{\mathrm{D}}
\newcommand{\textE}{\mathrm{E}}
\newcommand{\textF}{\mathrm{F}}
\newcommand{\textG}{\mathrm{G}}

\newcommand{\Dst}{\textD_{\st}}

\ifcomments
  \newcommand{\mar}[1]{\marginpar{\raggedright\tiny #1}}
\else
  \newcommand{\mar}[1]{}
\fi

\renewcommand{\(}{\textup{(}}
\renewcommand{\)}{\textup{)}}

\makeatletter
\newcommand{\mylabel}[2]{#2\def\@currentlabel{#2}\label{#1}}
\makeatother

\newcounter{step} 
\renewcommand{\thestep}{\arabic{step}}
\newenvironment{step}[1]%
{
\refstepcounter{step}
{\textit{Step \thestep:}#1}
}%

\usepackage{hyperref}
\usepackage{enumitem}

\theoremstyle{plain} 
 
\newtheorem{lm}[equation]{Lemma}
\newtheorem{lem}[equation]{Lemma}
\newtheorem{lemma}[equation]{Lemma}

\newtheorem{prop}[equation]{Proposition}
\newtheorem{proposition}[equation]{Proposition}
\newtheorem{thm}[equation]{Theorem}
\newtheorem{theorem}[equation]{Theorem}
\newtheorem{thmx}{Theorem}

\newtheorem{coroll}[equation]{Corollary}
\newtheorem{corollary}[equation]{Corollary}

\theoremstyle{definition}

\newtheorem{definition}[equation]{Definition}
\theoremstyle{remark}
\newtheorem{rk}[equation]{Remark}

\newtheorem{remark}[equation]{Remark}

\numberwithin{equation}{subsection}
\numberwithin{figure}{subsection}

\begin{document}

\title[Existence of supersingular representations]{On the existence of admissible supersingular representations of $p$-adic reductive groups}
\author{Florian Herzig} 
\address[F.\ Herzig]{Department of Mathematics, University of Toronto,
  40 St.\ George Street, Room 6290, Toronto, ON M5S 2E4, Canada}
\thanks{The first-named author was partially supported by a Simons Fellowship and an NSERC grant.  The second-named author was supported by the National Science Foundation under Award No.\ DMS-1400779, along with Manish Patnaik's Subbarao Professorship and NSERC grant RGPIN 4622.}
\email{herzig@math.toronto.edu}
\author{Karol Kozio\l}
\address[K.\ Kozio\l]{Department of Mathematics, University of Michigan, Ann Arbor, MI 48109-1043} 
\email{kkoziol@umich.edu}
\author{Marie-France Vign\'eras}
\address[M.-F.\ Vign\'eras]{Institut de Math\'ematiques de Jussieu, 4 place Jussieu, 75005 Paris, France}
\email{vigneras@math.jussieu.fr}

\begin{abstract}
  Suppose that $\bG$ is a connected reductive group over a finite extension $F/\qp$, and that $C$ is a field
  of characteristic $p$.
  We prove that the group $\bG(F)$ admits an irreducible admissible supercuspidal, or equivalently
  supersingular, representation over $C$.
\end{abstract}

\maketitle
\tableofcontents

\section{Introduction}
\label{sec:introduction}

Suppose that $F$ is a non-archimedean field of residue characteristic $p$ and that $\bG$ is a connected reductive algebraic
group over $F$.  
There has been a growing interest in understanding the smooth representation theory of the $p$-adic group $G := \bG(F)$ over a field
$C$ of characteristic $p$, going back to the work of Barthel--Livn\'e \cite{bib:BL-general} and Breuil \cite{Breuil} in
the case of $\bG = \mathbf{GL}_2$.  The latter work in particular demonstrated the relevance of the mod $p$ representation theory of
$p$-adic reductive groups to the $p$-adic Langlands program.

The results of \cite{bib:ahhv} (when $C$ is algebraically closed) and \cite{HenV} (for a general $C$ of characteristic $p$)
give a classification of irreducible admissible representations in terms of supercuspidal
$C$-representations of Levi subgroups of $G$. Here, an irreducible admissible smooth representation $\pi$ is said to be
\emph{supercuspidal} if it does not occur as subquotient of any parabolic induction $\Ind^G_P \sigma$, where $P$ is a proper
parabolic subgroup of $G$ and $\sigma$ an irreducible admissible representation of the Levi quotient of $P$.  
Unfortunately, so far, the supercuspidal representations themselves remain mostly mysterious, outside anisotropic groups,
$\GL_2(\qp)$ (\cite{bib:BL-general}, \cite{Breuil}), and some related cases (\cite{Abdellatif}, \cite{Cheng}, \cite{Koziol}, \cite{KoziolXu}). If $F/\qp$ is a non-trivial unramified extension, then
irreducible supercuspidal representations of $\GL_2(F)$ were first constructed by Pa\v{s}k\={u}nas \cite{paskunas:diags}; 
however, it seems hopelessly complicated to classify them \cite{BreuilPaskunas}, \cite{bib:Hu}. One additional challenge in
constructing supercuspidal representations is that irreducible smooth representations need not be admissible in general
(unlike what happens over $\C$), as was shown recently by Daniel Le \cite{le-nonadm}.

There are two ways to characterize supercuspidality in terms of Hecke actions. The first description assumes $C$ is algebraically closed and
uses weights and Hecke eigenvalues for any fixed choice $K$ of special parahoric subgroup (a weight is then an irreducible
representation of $K$). It was shown to be equivalent to supercuspidality in \cite{bib:ahhv}. The second description uses the
center of the pro-$p$ Iwahori--Hecke algebra. The equivalence between the second Hecke description and supercuspidality was shown in
\cite{OV} (when $C$ is algebraically closed) and~\cite{HenV} (for a general $C$ of characteristic $p$). In either description, supercuspidality is
characterized by the vanishing of certain Hecke operators, which is why supercuspidal representations are also called
\emph{supersingular}.

Our main theorem is the following:

\begin{thmx}\label{thm:main}
Suppose $F$ is of characteristic $0$, $\bG$ is any connected reductive algebraic group over $F$, and $C$ any field of characteristic
  $p$.  Then $G$ admits an irreducible admissible supersingular, or equivalently supercuspidal, representation over $C$.
\end{thmx}

This theorem is new outside the low rank cases mentioned above.  It provides an affirmative answer to Question 3 of \cite{ahhv:questions} when $\chr F = 0$, and carries out the announcement contained in \cite[\S\ III.26]{bib:ahhv}.  Note also that the analogous theorem for supercuspidal representations over $\mathbb{C}$ was proved by Beuzart-Plessis \cite{BP}. 

We now briefly explain our argument, which uses several completely different ideas. First, in
Section~\ref{sec:superc-repr} we reduce to the cases where $C$ is finite and $\bG$ is absolutely simple adjoint. If $\bG$ is
moreover anisotropic, then $G$ is compact and any irreducible smooth representation of $G$ is finite-dimensional (hence
admissible) and supercuspidal. If $\bG$ is isotropic, we distinguish three cases.

For most groups $\bG$ we show in Section~\ref{sec:proof-main-theorem} that there exists a discrete series representation $\pi$ of
$G$ over $\C$ that admits invariants under an Iwahori subgroup $\B$, and that has moreover the following property: the module
$\pi^{\B}$ of the Iwahori--Hecke algebra $H(G,\B)$ admits a $\Z[q^{1/2}]$-integral structure whose reduction modulo the
maximal ideal of $\Z[q^{1/2}]$ with residue field $\fp$ is supersingular. The Hecke modules $\pi^{\B}$ are constructed either
from characters (using \cite{Borel}) or reflection modules (using \cite{Lusztig} and \cite{GS}; the latter is needed
to handle unramified non-split forms of $\bPSO_8$).

Suppose from now on that $F$ is of characteristic zero, i.e.\ that $F/\qp$ is a finite extension.  
The $p$-adic version of the de George--Wallach limit multiplicity formula (\cite[App.\ 3]{DKV} plus \cite[Thm.\
K]{Kazhdan}) implies that the representation $\pi$ above embeds in $C^\infty(\Gamma\backslash G, \C)$ for some discrete cocompact subgroup $\Gamma$ of
$G$ (as $\chr F = 0$).
By construction we deduce that the Hecke module $C^\infty(\Gamma\backslash G/\B, \fp) = C^\infty(\Gamma\backslash G, \fp)^{\mathfrak{B}}$ of $\B$-invariants admits a
supersingular submodule. Crucially, by cocompactness of $\Gamma$ we know that $C^\infty(\Gamma\backslash G, \fp)$ is an
admissible representation of $G$. Picking any non-zero supersingular vector $v \in C^\infty(\Gamma\backslash G/\B, \fp)$, the
$G$-subrepresentation of $C^\infty(\Gamma\backslash G, \fp)$ generated by $v$ admits an irreducible quotient, which is
admissible (as $\chr F = 0$) and supersingular.

Unfortunately, this argument does not work for all groups $\bG$. We have the following exceptional cases:
\begin{enumerate}
\item $\PGL_n(D)$, where $n \ge 2$ and $D$ a central division algebra over $F$;
\item $\PU(h)$, where $h$ is a split hermitian form in $3$ variables over a ramified quadratic extension of $F$ or
a non-split hermitian form in $4$ variables over the unramified quadratic extension of $F$.
\end{enumerate}
Note that for the group $\PGL_n(D)$ with $n \ge 2$ the only discrete series representations $\pi$ having $\B$-invariant vectors
are the unramified twists of the Steinberg representation (by Proposition~\ref{prop fact2}(i) and the classification of Bernstein--Zelevinsky and Tadi\'c),
but then $\pi^{\B}$ is one-dimensional with non-supersingular reduction.

In the second exceptional case, where $G \cong \PU(h)$ for certain hermitian forms $h$, we use the theory of coefficient
systems and diagrams, building on ideas of Pa\v{s}k\={u}nas \cite{paskunas:diags}. 
See Section \ref{sec:supers-repr-rank}. Note that $\bG$ is of relative rank 1, so the
adjoint Bruhat--Tits building of $G$ is a tree, and our method works for all such groups.  In order to carry it out, we may apply the reductions in Section \ref{sec:superc-repr} and assume that $\bG$ is absolutely simple and simply connected.  Given a supersingular module $\Xi$ for the \emph{pro-$p$} Iwahori--Hecke algebra of $G$, we naturally construct a $G$-equivariant coefficient system (or cosheaf) $\cD_{\Xi}$ on the Bruhat--Tits tree of $G$.  The homology of $\cD_{\Xi}$ admits a smooth $G$-action, and any irreducible admissible quotient will be supersingular (by Proposition \ref{prop ss=sc}). To construct such a quotient, we define an auxiliary coefficient system $\cD'$, which is built out of injective envelopes of representations of certain parahoric subgroups, along with a morphism $\cD_{\Xi} \rightarrow \cD'$.  The image of the induced map on homology is admissible, and admits an irreducible quotient $\pi'$ which is itself admissible (as $\chr F = 0$) and supersingular.

In the first exceptional case, where $G \cong \PGL_n(D)$, we use a global method (see Section~\ref{sec:supers-repr-pgl_n}).
We find a totally real number field $\F^+$ and a compact unitary group  
$\u G$ over $\F^+$ such that $\u G(\F_v^+)$ is isomorphic to
$\GL_n(D)$ for a suitable place $v|p$ of $\F^+$. Then, fixing a level away from $v$ and taking the limit over all levels at
$v$, the space $S$ of algebraic automorphic forms of $\u G(\A_{\F^+}^\infty)$ over $\fpb$ affords an admissible smooth action of
$\u G(\F_v^+)$. Using automorphic induction and descent we construct an automorphic representation $\pi$ of $\u G(\A_{\F^+})$ whose
associated Galois representation $r_\pi$ has the property that its reduction modulo $p$ is irreducible locally at $v$. From
$\pi$ we get a maximal ideal $\m$ in the Hecke algebra (at good places outside $p$), and we claim that any irreducible
subrepresentation of the localization $S_{\m}$ is supercuspidal.

To prove the claim, we use the pro-$p$ Iwahori--Hecke criterion for supercuspidality and argue by contradiction.  If one of
the relevant Hecke operators has a non-zero eigenvalue, we lift to characteristic zero by a Deligne--Serre argument and
construct an automorphic representation $\pi'$ with Galois representation $r_{\pi'}$ having the same reduction as $r_\pi$
modulo $p$. Using local-global compatibility at $p$ for $r_{\pi'}$ and some basic $p$-adic Hodge theory we show that the
non-zero Hecke eigenvalue in characteristic $p$ implies that $r_{\pi'}$ is reducible locally at $v$, obtaining the desired
contradiction.

For our automorphic base change and descent argument we require results going slightly beyond \cite{labesse},
since our group $\u G$ is typically not quasi-split at all finite places. In the appendix, Sug Woo Shin explains the necessary
modifications.

Finally, we remark that we would expect Theorem~\ref{thm:main} to be true even when $\chr F = p$. So far this only
seems to be known for the groups $\GL_2(F)$ \cite{paskunas:diags}, outside trivial cases.
We crucially use that $\chr F = 0$ in (at least) the following ways: 
\begin{enumerate}
\item the existence of discrete cocompact subgroups, which fails for most groups if $\chr F = p$ \cite[\S3.4]{BorelHarder}, \cite[Cor.\ IX.4.8(iv)]{Margulis},
\item admissibility is preserved under passing to a quotient representation, 
\item the automorphic method in case of the group $\PGL_n(D)$.
\end{enumerate}

\subsection{Acknowledgements}\label{sec:acknowledgements}
The third-named author thanks Boris Pioline, Gordan Savin and the organizers of the conference on Automorphic Forms, Mock Modular Forms and String Theory in Banff (10/2017) for emails, discussions and their invitation, which led to closer examination of unramified minimal representations corresponding to the reflection modules of the affine Hecke algebras.  
She also thanks Volker Heiermann for a discussion on discrete modules, G\"unter Harder for emails on discrete cocompact subgroups, Jean-Loup Waldspurger for recollections on the antipode, and Guy Henniart for providing a proof below.

The first-named author thanks the University of Paris-Sud and the Mathematics Institute of Jussieu--Paris Rive Gauche, 
where some of this work was carried out.

We thank Noriyuki Abe for some helpful comments.

\subsection{Notation}\label{sec:notation}

Fix a prime number $p$, and let $F$ be a non-archimedean local field of residue characteristic $p$
(we will later assume that $\chr F = 0$, i.e.,\ 
that $F$ is a finite extension of $\Q_p$).   The field $F$ comes equipped with ring of integers $\O_F$ and residue field $k_F$ of cardinality $q$, a power of $p$. We fix a uniformizer $\varpi$, and let $\textnormal{val}_F$ and $|\cdot|_F$ denote the normalized valuation and normalized absolute value of $F$, respectively.

If $\bf{H}$ is an algebraic $F$-group, we denote by $H$ its group of $F$-points $\mathbf{H}(F)$.

Let $\bf G$ be a connected reductive $F$-group, $\bf T$ a maximal $F$-split subtorus of $\bf G$, $\bf B$ a minimal
$F$-parabolic subgroup of $\bf G$ containing $\bf T$, and $x_0$ a special point of the apartment of the adjoint Bruhat--Tits
building defined by $\bf T$.  We associate to $x_0$ and the triple ($\bf G$, $\bf T$, $\bf B $) the following data: 
\begin{itemize}
\item the center $\bf Z(G)$ of $\bf G$, 
\item the root system $\Phi \subset X^*(T)$, 
\item the set of simple roots $\Delta \subset \Phi$, 
\item the centralizer $\bf Z$ of $\bf T$, 
\item the normalizer $\boldsymbol{\mathcal{N}}$ of $\bf T$, 
\item the unipotent radical $\bf U$ of $\bf B$ (hence ${\bf B}=\bf Z \bf U$), and the opposite unipotent radical $\mathbf{U}\op$, 
\item the triples $({\bf{G}^{\rm sc}}, {\bf{T}^{\rm sc}}, {\bf{B}^{\rm sc}})$ and $({\bf{G}^{\rm ad}}, {\bf{T}^{\rm ad}}, {\bf{B}^{\rm ad}}),$ corresponding to the simply-connected covering of the derived subgroup and the adjoint group of $\bf G$,
\item the apartment $\mathscr{A} := X_*(T)/X_*(Z(G)^\circ)\otimes_{\mathbb{Z}}\mathbb{R}$ associated to $\mathbf{T}$ in the adjoint Bruhat--Tits building,
\item the alcove $\mathcal C$ of $\mathscr{A}$ with vertex $x_0$ lying in the dominant Weyl chamber with vertex $x_0$,
\item the Iwahori subgroups $\mathfrak B$ and $\mathfrak{B}^{\rm sc}$ of $G$ and $G^{\rm sc}$, respectively, fixing $\mathcal C$ pointwise,
\item the pro-$p$-Sylow subgroup $\mathfrak{U}$ of $\mathfrak{B}$.
\end{itemize}

Given a field $L$, we denote by $\overline{L}$ a fixed choice of algebraic closure.  We fix a field $C$ of characteristic $c \in \{0,2,3,5,7,\ldots\}$, which will serve as the field of coefficients for the modules and representations appearing below.  In our main result we will assume $c = p$.

Suppose $K$ is a compact open subgroup of $G$, and $R$ is a commutative ring.  We define the Hecke algebra associated to this data to be the $R$-algebra 
$$H_R(G,K) := \End_G R[K\backslash G].$$  
If $R = \mathbb{Z}$, we simply write $H(G,K)$.  In our applications below, we will often assume that $K = \mathfrak{B}$ or $K = \mathfrak{U}$.

Given a module or algebra $X$ over some ring $R$ and a ring map $R \rightarrow R'$, we let $X_{R'} := X\otimes_R R'$ denote the base change.

Other notation will be introduced as necessary in subsequent sections.

\section{Iwahori--Hecke algebras}
\label{sec:IHalgs}

In this section we review some basic facts concerning Iwahori--Hecke algebras and their (supersingular) modules.  We will use these algebras extensively in our construction of supercuspidal $G$-representations.  See \cite{Vig1}, \cite{Vigcenter}, and \cite{Vigss} for references.

\subsection{Definitions}
\label{sec:iwahori-hecke-alg}

Recall that we have defined the Iwahori--Hecke ring as
$$H (G, \mathfrak B) = \End_{G}\mathbb{Z}[\mathfrak B \backslash G].$$ 
We have an analogous ring $H(G\sc, \mathfrak{B}\sc)$ for the simply-connected group.  The natural ring homomorphism $H(G\sc, \mathfrak B\sc) \to H(G, \mathfrak B )$ (induced by the covering $G\sc \to G$ of the derived subgroup) is injective, so we identify $H(G\sc, \mathfrak B\sc)$ with a subring of $H(G, \mathfrak{B})$.  We first discuss presentations for these rings.

There is a canonical isomorphism
$$j\sc : H(G\sc, \mathfrak B\sc) \congto H (W,S, q_s),$$ 
where $H := H (W,S, q_s)$ is the Hecke ring of an affine Coxeter system $(W,S)$ with parameters $\{q_s := q^{d_s}\}_{s\in S}$.  The $d_s$ are positive integers, which we will abusively also refer to as the parameters of $G$.  Thus, $H (W,S, q_s)$ is a free $\mathbb{Z}$-module with basis $\{T_w\}_{w\in W}$, satisfying the braid and quadratic relations:
\begin{alignat*}{2}
 T_w T_{w'} & = T_{ww' } & & \quad  \textrm{for}~ w,w'\in W,~ \ell(w)+\ell(w')=\ell(ww'), \\  
  (T_s - q_s)(T_s + 1) & = 0 & & \quad \text{for}~ s\in S.
\end{alignat*}
Here $\ell: W \rightarrow \mathbb{Z}_{\geq 0}$ denotes the length function with respect to $S$.  We identify $H(G\sc, \mathfrak B\sc)$ with $H$ via $j\sc$.

In order to describe $H(G,\mathfrak{B})$, we require a larger affine Weyl group.  We define the \emph{extended affine Weyl group} to be
$$\widetilde{W} := \mathcal{N}/Z_0,$$
where $Z_0$ is the unique parahoric subgroup of $Z$.  The group $\widetilde{W}$ acts on the apartment $\mathscr{A}$, and permutes the alcoves of $\mathscr{A}$ transitively.  We let $\Omega$ denote the subgroup of $\widetilde{W}$ stabilizing $\mathcal{C}$.  The affine Weyl group $W$ is isomorphic to a normal subgroup of $\widetilde{W}$, and permutes the alcoves \emph{simply} transitively.  We therefore have a semidirect product decomposition
$$\widetilde{W} = W \rtimes \Omega.$$
The function $\ell$ extends to $\widetilde{W}$ by setting $\ell(uw) = \ell(wu) = \ell(w)$ for $u\in \Omega, w\in W$.  In particular, we see that $\Omega$ is the group of length-zero elements of $\widetilde{W}$.

Let $\Sigma$ denote the reduced root system whose extended Dynkin diagram $\mathsf{Dyn}$ is equal to the Dynkin diagram of $(W,S)$, and let $\mathsf{Dyn}'$ denote the Dynkin diagram $\mathsf{Dyn}$ decorated with the parameters $\{d_s\}_{s\in S}$.  The quotient of $\Omega$ by the pointwise stabilizer of $\mathcal C$ in $\Omega$ is isomorphic to a subgroup $\Psi$ of $\Aut(W,S,d_s)$, the group of automorphisms of $\mathsf{Dyn}'$.  Thus, $\Omega$ acts on $\mathsf{Dyn}'$ and consequently on $H(W,S,q_s)$, and the isomorphism $j\sc$ extends to an isomorphism
\begin{equation}\label{eq Had}
j : H (G, \mathfrak B) \congto \mathbb Z[\Omega] \mathbin{\widetilde{\otimes}} H (W,S, q_s),
\end{equation} 
where $\widetilde{\otimes}$ denotes the twisted tensor product.  The generalized affine Hecke ring $\widetilde{H} := \mathbb Z[\Omega] \mathbin{\widetilde{\otimes}} H (W,S, q_s)$ as above is the free $\mathbb Z$-module with basis $\{T_w\}_{w\in \widetilde{W}}$, satisfying the braid and quadratic relations:
\begin{alignat}{2}
 T_w T_{w'} & = T_{ww' } & &  \quad \textrm{for}~ w,w'\in \widetilde W,~ \ell(w)+\ell(w')=\ell(ww'), \label{eq braid}\\  
  (T_s - q_s)(T_s + 1) & = 0 & &  \quad \text{for}~  s\in S.\label{eq Tsq}
\end{alignat} 
Thus, we see that the Iwahori--Hecke ring $H(G,\mathfrak{B})$ is determined by the type of $\Sigma$, the parameters $\{d_s\}_{s\in S}$, and the action of $\Omega$ on $\mathsf{Dyn}'$.

The group $\widetilde{W}$ forms a system of representatives for the space of double cosets $\mathfrak{B}\backslash G/\mathfrak{B}$.  Under the isomorphism $j$, the element $T_w\in \widetilde{H}$ for $w\in \widetilde{W}$ corresponds to the endomorphism sending the characteristic function of $\mathfrak{B}$ to the characteristic function of $\mathfrak{B}n\mathfrak{B}$, where $n\in \mathcal{N}$ lifts $w$.

Finally, let $\mathbf{P} = \mathbf{MN}$ denote a standard parabolic $F$-subgroup of $\mathbf{G}$ (meaning $\mathbf{B} \subset \mathbf{P}$), and suppose that $\mathbf{M}$ contains $\mathbf{T}$.  Then the group $M\cap \mathfrak{B}$ is an Iwahori subgroup of $M$, and we may form the algebra
$$H(M, M\cap \mathfrak{B}) = \textnormal{End}_M\mathbb{Z}[(M\cap \mathfrak{B})\backslash M].$$
It is \emph{not} a subalgebra of $H(G,\mathfrak{B})$ in general.  The basis of $H(M, M\cap \mathfrak{B})$ will be denoted $T_w^M$, where $w$ is an element of the extended affine Weyl group associated to $M$.

\subsection{Dominant monoids}
\label{sec:dominant-monoids}

The subgroup 
$$\Lambda := Z/Z_0$$ 
of $\widetilde{W} = \mathcal{N}/Z_0$ is commutative and finitely generated, and its torsion subgroup is equal to $\widetilde{Z}_0/Z_0$, where $\widetilde{Z}_0$ denotes the maximal compact subgroup of $Z$.  (When the group $\bf G$ is $F$-split or semisimple and simply connected, we have $Z_0 = \widetilde{Z}_0$.)  The short exact sequence 
$$1 \rightarrow \Lambda \rightarrow \mathcal{N}/Z_0 \rightarrow \mathcal{N}/Z \rightarrow 1$$ 
splits, identifying the (finite) Weyl group $W_0 := \mathcal{N}/Z$ of $\Sigma$ with $\textnormal{Stab}_{W}(x_0)$.  We therefore obtain semidirect product decompositions 
$$\Lambda \rtimes W_0  = \widetilde{W}$$
and
$$\Lambda\sc\rtimes W_0 = W,$$ 
where $\Lambda\sc := \Lambda \cap W$.

Given a subgroup $J$ of $Z$, we define 
$$\Lambda_J := JZ_0/Z_0 \subset \Lambda.$$   
We analyze $\Lambda_J$ for various groups $J$ presently.

Let $T_0$ denote the maximal compact subgroup of $T$, and note that $T_0 = Z_0 \cap T$.  This implies that the inclusion $T \hookrightarrow TZ_0$ induces an isomorphism $T/T_0 \congto TZ_0/Z_0 = \Lambda_T$, and therefore the map 
\begin{align}
 X_*(T) & \stackrel{\sim}{\longrightarrow} \Lambda_T  \label{cochar isom} \\
\mu & \longmapsto \lambda_\mu  := \mu(\varpi)Z_0/Z_0 \notag
\end{align}
is a $W_0$-equivariant isomorphism.

Recall that we have a unique homomorphism 
$$\nu: Z \rightarrow \mathscr{A},$$
determined by the condition
$$\langle \alpha, \nu(t)\rangle = -\textnormal{val}_F(\alpha(t))$$
for $t\in T$ and $\alpha\in \Phi$.  We claim that the kernel of $\nu$ is the saturation of $Z(G)\widetilde{Z}_0$ in $Z$, i.e., the set of all elements $z\in Z$ such that $z^n\in Z(G)\widetilde{Z}_0$ for some $n\geq 1$.  Indeed, the kernel of $\nu$ contains $Z(G)$ and $\widetilde{Z}_0$, and the group $Z/Z(G)\widetilde{Z}_0$ is commutative and finitely generated.  This gives an induced map
\begin{equation}\label{nu}
\nu: Z/Z(G)\widetilde{Z}_0 \rightarrow \mathscr{A}.
\end{equation}
We note the following three facts: (1) the image of $T\ad$ in $Z/Z(G)\widetilde{Z}_0$ is of finite index (cf.\ comments following (16) in \cite{Vig1}); (2) the $\mathbb{Z}$-span of the coroots $\Phi^\vee$ is of finite index in $X_*(T\ad)$; (3) $\nu(\alpha^\vee(\varpi^{-1})) = \alpha^\vee$ for any coroot $\alpha^\vee \in X_*(T)$.  Combining these, we see that the image of \eqref{nu} has the same rank as $Z/Z(G)\widetilde{Z}_0$, which is equal to the rank of $X_*(T\ad)$.  Therefore, the kernel of \eqref{nu} is exactly the torsion subgroup of $Z/Z(G)\widetilde{Z}_0$.  This gives the claim.

Since $Z_0$ is contained in the kernel of $\nu$, the group $\Lambda$ acts by translation on $\mathscr{A}$ via $\nu$. 
Therefore, $\Lambda_{\ker\nu}$ is the pointwise stabilizer of $\mathcal{C}$ in $\Lambda$.  Similarly, one easily checks that $\Lambda_{\ker\nu}$ is the pointwise stabilizer of $\mathcal{C}$ in $\Omega$.  (In fact, we have $\Lambda \cap \Omega = \Lambda_{\ker\nu}$, cf.\ \cite[Cor.\ 5.11]{Vig1}.)
Hence, we obtain
\begin{equation}\label{eq fix 1}
  \Omega /\Lambda_{\ker\nu} \congto \Psi,
\end{equation}
and the embeddings of $\Lambda$ and $\Omega$ into $\widetilde W$ induce
\begin{equation}\label{eq fix 2} 
\Lambda /(\Lambda_{\ker\nu} \times \Lambda\sc) \congto \widetilde{W} /(\Lambda_{\ker\nu} \times W) \tocong \Omega/\Lambda_{\ker\nu}.
\end{equation}

An element $\lambda\in \Lambda$ is called \emph{dominant} (and $\lambda^{-1}$ is called \emph{anti-dominant}), if 
$$z(U\cap {\mathfrak B}) z^{-1} \subset U\cap {\mathfrak B}$$ 
for any $z\in Z$ which lifts $\lambda$. We let $\Lambda^+$ denote the monoid consisting of dominant elements of $\Lambda$, and
similarly for any subgroup $\Lambda' \le \Lambda$ we define $\Lambda'^+ := \Lambda' \cap \Lambda^+$.
Using the isomorphism \eqref{cochar isom}, we say a cocharacter $\mu\in X_*(T)$ is \emph{dominant} if $\lambda_\mu$ is, and let $X_*(T)^+$ denote the monoid consisting of dominant elements of $X_*(T)$.  
The group of invertible elements in the dominant monoid $\Lambda^+$ is exactly the subgroup $\Lambda_{\ker\nu}$, and the invertible elements of $\Lambda^{\textnormal{sc},+}$ are trivial.

\begin{lemma}\label{lem fin} 
  The subgroup $\Lambda_{Z(G) } \times \Lambda\sc$ \(resp.\ $\Lambda_T$\) of $\Lambda$ is finitely generated of finite index.  
  The submonoid $\Lambda_{Z(G) } \times \Lambda^{\textnormal{sc},+}$ \(resp.\ $\Lambda_T^+$\) of the dominant monoid $\Lambda^+$ is finitely generated of finite index.
\end{lemma}

Here, we say that a submonoid $N$ of a commutative monoid $M$ has finite index if $M = \bigcup_{i=1}^n (N+x_i)$ for some $x_i \in M$.
If $M$ is finitely generated, then $dM$ is of finite index in $M$ for all $d \ge 1$.

\begin{proof} 
The groups $\ker\nu/Z(G)\widetilde{Z}_0$ and $\widetilde{Z}_0/Z_0$ are finite, and equations \eqref{eq fix 1} and \eqref{eq fix 2} imply that $\Lambda /( \Lambda_{\ker\nu} \times
  \Lambda\sc)$ is isomorphic to the finite group $\Psi$.  Thus, we see that the commutative group $\Lambda_{Z(G) } \times \Lambda\sc$ is a finitely generated, finite index subgroup of $\Lambda$.  Similarly, $\Lambda_T$ is finite free and it is well known that it is of finite index in $\Lambda$. 
Gordan's lemma implies the second assertion (as in the proof of \cite[7.2 Lem.]{HV1}).
\end{proof}

\subsection{Bernstein elements}
\label{sec:bernstein-elements}

Let $w\in \widetilde{W}$, and let $w = us_1\cdots s_n$ be a reduced expression, with $u\in \Omega, s_i\in S$.  We set $q_w := q_{s_1}\cdots q_{s_n},$ and define 
$$T_s^* := T_s - q_s + 1 \quad\textnormal{and}\quad T_w^* := T_uT_{s_1}^*\cdots T_{s_n}^*.$$  
Then $T_w T^*_{w^{-1}} = q_w,$ and the linear map defined by $T_w \mapsto (-1)^{\ell(w)}T_w^*$ is an automorphism of $\widetilde{H}$.

For $\mu\in X_*(T)$, we let $\mathcal O_\mu\subset \Lambda$ denote the $W_0$-orbit of $\lambda_\mu$.  We then define
$$z_\mu:= \sum_{\lambda\in \mathcal O_\mu} E_{\lambda},$$  
where $E_\lambda$ are the integral Bernstein elements of $\widetilde{H}$ corresponding to the spherical orientation induced by $\Delta$ (\cite[Cor.\ 5.28, Ex.\ 5.30]{Vig1}).  
Precisely, they are characterized by the relations
\begin{align}\label{eq E}
  E_{\lambda} & = \begin{cases} T_\lambda & \text{if $\lambda$ is anti-dominant,}\\ 
  T_\lambda^* & \text{if $\lambda$ is dominant,} \end{cases} \\
\label{bern product}  E_{\lambda_1} E_{\lambda_2} & = (q_{\lambda_1}q_{\lambda_2}q_{\lambda_1 \lambda_2}^{-1})^{1/2} E_{\lambda_1 \lambda_2} \quad \text{if}~\lambda_1, \lambda_2\in \Lambda,
\end{align} 
where we take the positive square root.  (If $\lambda_1, \lambda_2$ are both dominant (or anti-dominant), then $E_{\lambda_1} E_{\lambda_2} =E_{\lambda_1 \lambda_2}$.)  The elements $z_\mu$ are central in $\widetilde{H}$, and when $\mu\in X_*(T\sc)$, $z_\mu$ lies in $H$.

We let $\mathcal A$ denote the commutative subring of the generalized affine Hecke ring $\widetilde{H}$ with $\mathbb Z$-basis $\{E_{\lambda}\}_{\lambda\in \Lambda}$.  When $G = Z$, we have $\widetilde{H} = \mathcal{A} \cong \mathbb Z[\Lambda]$, but $\mathcal A$ is not isomorphic to $\mathbb Z[\Lambda]$ in general.   
The rings $\mathcal{A}, \widetilde{H}$, and the center of $\widetilde{H}$ are finitely generated modules over the central subring with basis $\{z_\mu\}_{\mu\in X_*(T)}$, which is itself a finitely-generated ring.

\subsection{Supersingular modules}
\label{sec:supers-discr-modul}

We now discuss supersingular Hecke modules.

Recall that $C$ is our coefficient field of characteristic $c$.  We define $H_C := H \otimes C$ and $\widetilde{H}_C := \widetilde{H} \otimes C$, which are isomorphic to the Iwahori--Hecke algebras $H_C(G\sc, \mathfrak{B}\sc)$ and $H_C(G, \mathfrak{B})$, respectively.

\begin{definition}[{cf.\ \cite[\S 5.1(3)]{OV}}]
\label{def ssing}
Let $M$ be a non-zero right $\widetilde{H}_C $-module.  A non-zero element $v\in M$ is called \emph{supersingular} if $v \cdot z_\mu^n=0$ for all $\mu\in X_*(T)^+$ such that $-\mu \not\in X_*(T)^+$, and all sufficiently large $n$.  The $\widetilde{H}_C $-module $M$ is called \emph{supersingular} if all its non-zero elements are supersingular. We make a similar definition for modules over $H_C$, using the monoid $X_*(T\sc)^+$.  
\end{definition}

We remark that the definition of a supersingular module differs slightly from that of \cite[Def.\ 6.10]{Vigss}. There it was
required that $c = p$, and that $M \cdot z_\mu^n =0$ for all $\mu \in X_*(T)^+$ such that $-\mu \not\in X_*(T)^+$ and $n$ sufficiently large.
 
\begin{lm}\label{lm:fact}\hfill
\begin{enumerate}
  \item Any simple $\widetilde{H}_C $-module is finite dimensional, and is semisimple as an $H_C $-module. \label{lm:fact 1}
  \item If $c \nmid p|W_0|$, then $\widetilde{H}_C$ does not admit any simple supersingular modules.  \label{lm:fact 2}
  \item If $c = p$, a simple $\widetilde{H}_C $-module is supersingular if and only if its restriction to $H_C$ is supersingular. \label{lm:fact 3}
  \end{enumerate}
\end{lm}

\begin{proof}
(i) The first statement follows from \cite[\S 5.3]{VigDurham}.  For the second part, note that there exists a finite index subgroup ${\Omega}'$ of ${\Omega}$ which acts trivially on $H$ (for example, we may take $\Omega' = \Lambda_{\ker\nu}$).  Set $H'_C := C[\Omega'] \otimes_C H_C$. Any simple $H_C$-module $N$ extends trivially to an $H'_C$-module $N'$, and the restriction of $N'\otimes_{H'_C} \widetilde{H}_C$ to $H_C$ is a finite direct sum $\bigoplus_{u\in \Omega/\Omega'} N^u$ of (simple) conjugates $N^u$ of $N$ by elements $u\in \Omega$. If $M$ is a simple $\widetilde{H}_C$-module and $N$ is contained in $M|_{H_C}$, then $M$ is a quotient of $N'\otimes_{H'_C} \widetilde{H}_C$ (and thus the restriction of $M$ is semisimple).
        
(ii) It suffices to assume $C$ is algebraically closed.  Let $M$ denote a simple supersingular module.  Since the center of $\widetilde{H}$ is commutative and $M$ is finite dimensional, there exists an eigenvector $v \in M$ with eigenvalues $\chi$ for the action of the center. 
Supersingularity then implies
\begin{equation}\label{ssing char neq p}
0 = v\cdot z_{\mu'} = \chi(z_{\mu'})v
\end{equation}
for any $\mu' \in X_*(T)^+$ such that $-\mu' \not\in X_*(T)^+$.

Choose $\mu \in X_*(T)^+$ lying in the interior of the dominant Weyl chamber, so in particular $-\mu \not\in X_*(T)^+$, and let $w_\circ\in W_0$ denote the longest element.
We claim that
\begin{equation}\label{ssing char neq p 2}
z_{\mu}z_{-w_\circ(\mu)} = q_{\lambda_\mu}|W_0| z_0 + \sum_{\substack{\mu'\in X_*(T)^+ \\ \ell(\lambda_{\mu'}) > 0}} a_{\mu'} z_{\mu'}
\end{equation}
for some $a_{\mu'} \in \mathbb{Z}$.  To see this, note that the product of the orbits $\mathcal{O}_\mu \cdot \mathcal{O}_{-w_\circ(\mu)}$ consists of elements of the form $\lambda_{w(\mu)}\lambda_{-w'w_\circ(\mu)}$, where $w,w'\in W_0$.  If the length of $\lambda_{w(\mu)}\lambda_{-w'w_\circ(\mu)}$ is 0, then \cite[Cor.\ 5.11]{Vig1} implies $w(\mu) - w'w_\circ(\mu)$ is orthogonal to every simple root.  Since this element is also a sum of coroots, we conclude that $w(\mu) - w'w_\circ(\mu) = 0$, which implies $w = w'w_\circ$, as the $W_0$-stabilizer of $\mu$ is trivial.  The product formula \eqref{bern product} then gives equation \eqref{ssing char neq p 2}.

Now, for $\mu' \in X_*(T)^+$, the condition $-\mu' \not\in X_*(T)^+$ is equivalent to $\ell(\lambda_{\mu'}) > 0$.  
Applying $\chi$ to both sides of \eqref{ssing char neq p 2} and using \eqref{ssing char neq p} (for varying $\mu'$) gives $q_{\lambda_\mu}|W_0| = 0$, a contradiction.  
        
(iii)  This follows from \cite[Cor.\ 6.13]{Vigss} and part (i).  
\end{proof}

\begin{rk}
The conditions in part \ref{lm:fact 2} of the above lemma are necessary: when $c \neq p$ divides $|W_0|$, there exist non-zero supersingular modules.  For an example, suppose $\mathbf{G} = \mathbf{SL}_2$, $q$ is odd, and $c = 2$.  Then $H_C = \widetilde{H}_C$ admits a unique character $\chi$, which sends $T_s$ to 1 for each $s\in S$.  If we let $\mu := (1,-1) \in X_*(T)^+$, then
$$z_\mu = T_{s_1}T_{s_2} + T_{s_2}T_{s_1},$$
where $S = \{s_1, s_2\}$.  Thus, we have $\chi(z_\mu) = 0$.  By induction, and using the assumption $c = 2$, we see that the element $z_{k\mu}$ lies in the ideal of the center generated by $z_\mu$, for every $k \geq 1$.  From this, we conclude that $\chi$ is supersingular.
\end{rk}

\section{On supercuspidal representations}
\label{sec:superc-repr}

The aim of this section is to collect various results concerning supercuspidal representations.  We first state Proposition \ref{prop ss=sc}, which gives a convenient criterion for checking that a given irreducible admissible representation is supercuspidal when $\chr C = p$.  Propositions \ref{prop Cfinite} and \ref{prop redadjoint} allow us to make further reductions: in order to prove that $\bG(F)$ admits an irreducible admissible supercuspidal $C$-representation when $\chr F = 0$ and $\chr C = p$, it suffices to assume that $C$ is finite and $\bG$ is absolutely simple, adjoint, and isotropic.

\subsection{Supercuspidality criterion}
\label{sec:superc-crit}

We begin with a definition.

\begin{definition} 
Let $R$ be a subfield of $C$.  We say that a $C$-representation $\pi$ of $G$ \emph{descends to $R$} if there exists an $R$-representation $\tau$ of $G$
  and a $G$-equivariant $C$-linear isomorphism 
  $$\varphi: C\otimes_{R} \tau \congto \pi.$$ 
  We call $\varphi$ (and more often $\tau$) an \emph{$R$-structure} of $\pi$, or a \emph{descent of $\pi$ to $R$}.
\end{definition}

We now describe the scalar extension of an irreducible admissible $C$-representation $\pi$ of $G$ \cite{HenV}.  Given such a $\pi$, the commutant $D := \End_C(\pi)$ is a division algebra of finite dimension over $C$. Let $E$ denote the center of $D$, $E_s/C$ the maximal separable extension contained in $E/C$ and $\delta$ the reduced degree of $D/E$. Let $\overline{L}$ be an algebraically closed field containing $E$ and $\pi_{\overline{L}}$ the scalar extension of $\pi$ from $C$ to $\overline{L}$.

\begin{prop}[{\cite[Thms.\ I.1, III.4]{HenV}}] \label{prop fact HenV} 
  The length of $\pi_{\overline{L}}$ is $\delta [E:C]$ and
  $$\pi_{\overline{L}} \cong  \bigoplus_{i\in \Hom_C(E_s,\overline{L} )}\pi_i^{\oplus \delta}$$ 
  where each $\pi_i$ is indecomposable with commutant
  $\overline{L}\otimes_{i,E_s}E$, descends to a finite extension of $C$, has length $[E:E_s]$, and its irreducible
  subquotients are pairwise isomorphic, say to $\rho_i$. The $\rho_i$ are admissible, with commutant $\overline{L}$,
  $\Aut_C(\overline{L})$-conjugate, pairwise non-isomorphic, and descend to a finite extension of $C$. Any descent of $\rho_i$
  to a finite extension $C'/C$, viewed as $C$-representation of $G$, is $\pi$-isotypic of finite length.
\end{prop}

\begin{proof}
  By \cite[Thms.\ I.1, III.4]{HenV}, it suffices to prove that if $\rho_i$ descends to a $C'$-representation
  $\rho'_i$ with $C'/C$ finite, then $\rho'_i$ is $\pi$-isotypic of finite length.
  Then $(\rho'_i)_{\overline{L}}$ injects into $\pi_{\overline{L}}$, and so $\rho'_i$ injects into $\pi_{C'}$ by \cite[Rk.\ II.2]{HenV}, which implies
  the claim.
\end{proof}

In particular, any irreducible admissible $C$-representation $\pi$ with commutant $C$ is absolutely irreducible in the sense
that its base change $\pi_L$ is irreducible for any field extension $L/C$. For example, this holds when $C$ is algebraically closed.

Given an irreducible admissible $C$-representation $\pi$, the space $\pi^{\mathfrak{U}}$ of $\mathfrak{U}$-invariants comes equipped with a right action of the pro-$p$ Iwahori--Hecke algebra $H_C(G,\mathfrak{U})$.  This algebra has a similar structure to that of $H_C(G,\mathfrak{B})$.  In particular, we have analogous definitions of the Bernstein elements $E_\lambda$ ($\lambda \in \Lambda$) and the central elements $z_\mu$ ($\mu \in X_*(T)$), as well as an analogous notion of supersingularity for right $H_C(G,\mathfrak{U})$-modules (cf.\ Definition \ref{def ssing}).  We say an irreducible admissible $C$-representation $\pi$ is \emph{supersingular} if the right $H_C(G,\mathfrak{U})$-module $\pi^{\mathfrak{U}}$ is supersingular.

Finally, recall that an irreducible admissible $C$-representation $\pi$ of $G$ is said to be \emph{supercuspidal} if it is not a subquotient of $\Ind_P^G\tau$ for any parabolic subgroup $P = MN \subsetneq G$ and any irreducible admissible representation $\tau$ of the Levi subgroup $M$.

\begin{proposition}[Supercuspidality criterion]\label{prop ss=sc}
  Assume $c = p$. Suppose that $\pi$ is an irreducible admissible $C$-representation of $G$. The following are equivalent:
  \begin{enumerate}
  \item $\pi$ is supercuspidal;
  \item $\pi$ is supersingular; 
  \item $\pi^{\fU}$ contains a non-zero supersingular element;
  \item every subquotient of $\pi^{\fU}$ is supersingular;
  \item some subquotient of $\pi^{\fU}$ is supersingular. 
  \end{enumerate}
\end{proposition}

\begin{proof}
  We have (i)$\eq$(ii)$\eq$(iii) by \cite[Thms.\ I.13, III.17]{HenV}.  Since (ii)$\Rightarrow$(iv)$\Rightarrow$(v), it
  suffices to show that (v)$\Rightarrow$(ii).  Let $\overline{C}$ denote an algebraic closure of $C$. Say $\pi^{\fU}$ has
  supersingular subquotient $M$. Then $(\pi^{\fU})_{\overline{C}} \cong (\pi_{\overline{C}})^{\fU}$ has subquotient $M_{\overline{C}}$, and
  $M_{\overline{C}}$ is clearly supersingular.  By Proposition~\ref{prop fact HenV} there exists an irreducible admissible constituent $\rho$ of $\pi_{\overline{C}}$
  such that the $H_C(G,\mathfrak{U})$-module $\rho^{\mathfrak{U}}$ shares an irreducible constituent with $M_{\overline{C}}$.  
  In particular, $\rho^{\mathfrak{U}}$ has a supersingular subquotient, and \cite[Thm.\ 3]{OV} implies $\rho$ is supersingular.  
  Then \cite[Lem.\ III.16 2)]{HenV} implies that $\pi$ is supersingular.
\end{proof}

\begin{remark}
When $\pi^{\mathfrak{B}} \neq 0$, the above criterion holds with $\pi^{\fU}$ replaced by $\pi^{\mathfrak{B}}$ in items (iii), (iv), and (v).  This follows from the fact that $\pi^{\mathfrak{B}}$ is a direct summand of $\pi^{\mathfrak{U}}$ as an $H_C(G,\mathfrak{U})$-module, and the action of $H_C(G,\mathfrak{U})$ on $\pi^{\mathfrak{B}}$ factors through $H_C(G,\mathfrak{B})$.  
\end{remark}

We now discuss how supercuspidality behaves under extension of scalars.  We require a preliminary lemma.

 \begin{lm}\label{lm:scalar-res}
   Suppose that $C'/C$ is a finite extension and that $\pi'$ is an irreducible admissible $C'$-representation of $G$. Then
   $\pi'|_{C[G]} \cong \pi^{\oplus n}$ for some irreducible admissible $C$-representation $\pi$ of $G$ and some $n \ge 1$.
 \end{lm}

 \begin{proof}
   Let $\overline{C}$ be an algebraic closure of $C$. Then the finite-dimensional $\overline{C}$-algebra $A := C' \otimes_C \overline{C}$ is of finite length over itself.
   The simple $A$-modules are given by $\overline{C}$ with $C'$ acting via the various $C$-embeddings $C' \to \overline{C}$. 
   It follows that $\pi'|_{C[G]} \otimes_C \overline{C} \cong \pi' \otimes_{C'} A$ is of finite length as $\overline{C}$-representation
   by Proposition~\ref{prop fact HenV}. So $\pi'|_{C[G]}$ is of finite length. If $\pi$ denotes an irreducible submodule,
   then $\sum_i \lambda_i \pi = \pi'|_{C[G]}$, where $\{\lambda_i\}_{i=1}^m$ is a basis of $C'/C$. It follows that $\pi'|_{C[G]} \cong \pi^{\oplus n}$
   for some $n \leq m$. Moreover $\pi$ is admissible, as $\pi'|_{C[G]}$ is.
 \end{proof}

\begin{proposition}\label{prop HenV} Let $\overline{L}$ denote an algebraically closed field containing $C$.
If $c \ne p$, we assume that $\overline{L} = \overline{C}$ is an algebraic closure of $C$.

A $C$-representation $\pi$ is supercuspidal if and only if some irreducible subquotient $\rho$ of $\pi_{\overline{L}}$
is supercuspidal, if and only if every irreducible subquotient $\rho$ of $\pi_{\overline{L}}$ is supercuspidal.
\end{proposition}

\begin{proof}
  If $c = p$, we note that $\pi$ is supercuspidal if and only if $\pi$ is supersingular by Proposition \ref{prop ss=sc}.
  This is equivalent to some/every subquotient of $\pi_{\overline{L}}$ being supersingular \cite[Lem.\ III.16 2)]{HenV},
  or equivalently supercuspidal (again by \cite[Thm.\ I.13]{HenV}).

  Now suppose that $c \ne p$ and $\overline{L} = \overline{C}$.
  Recall that parabolic induction $\Ind_P^G$ is exact, and commutes with scalar extensions and restrictions \cite[Prop.\ III.12(i)]{HenV}. 
   If $\pi$ is not supercuspidal, then $\pi$ is a subquotient of $\Ind_P^G \tau$ for some proper parabolic $P = MN$ and
   irreducible admissible $C$-representation $\tau$ of $M$. Then $\pi_{\overline{C}}$ is a subquotient of
   $(\Ind_P^G\tau)_{\overline{C}} \cong \Ind_P^G(\tau_{\overline{C}})$. In
   particular, each irreducible (admissible) subquotient $\pi'$ of $\pi_{\overline{C}}$ is a subquotient of $\Ind_P^G \tau'$ for
   some irreducible (admissible) subquotient $\tau'$ of $\tau_{\overline{C}}$. Hence none of the $\pi'$ are supercuspidal.

For the converse, suppose by contradiction that $\pi_{\overline{C}}$ has an
   irreducible subquotient $\rho$ that is not supercuspidal, i.e.\ $\rho$ is a subquotient of $\Ind_P^G \tau$ for some proper
   parabolic $P = MN$ and irreducible admissible $\overline{C}$-representation $\tau$ of $M$. By \cite[II.4.7]{Viglivre} (as $c \ne p$),
   respectively by Proposition~\ref{prop fact HenV}, we can choose a
   finite extension $C'/C$ with $C' \subset \overline{C}$ such that $\tau$, respectively all irreducible constituents of $\Ind_P^G \tau$ and $\pi_{\overline{C}}$, can be
   defined over $C'$.  Write $\tau \cong (\tau')_{\overline{C}}$ for some $C'$-representation $\tau'$. Say the irreducible
   subquotients of $\Ind_P^G \tau'$ are $\sigma_1$, \dots, $\sigma_n$. So by our choice of $C'$, we know that
   $\rho \cong (\sigma_i)_{\overline{C}}$ for some $i$.  As $\sigma_i$ is a subquotient of $\Ind_P^G \tau'$, 
   we see that $\sigma_i|_{C[G]}$ is a subquotient of $\Ind_P^G(\tau'|_{C[M]})$.  But
   $\sigma_i|_{C[G]}$ is $\pi$-isotypic by Proposition~\ref{prop fact HenV} and $\tau'|_{C[M]}$ has finite length by Lemma~\ref{lm:scalar-res}, so $\pi$ is a
   subquotient of $\Ind_P^G \tau''$ for some irreducible (admissible) subquotient $\tau''$ of $\tau'|_{C[M]}$.
\end{proof}

\subsection{Change of coefficient field}
\label{S Cfinite}

This section contains the proof of the following result.

\begin{proposition}[Change of coefficient field]\label{prop Cfinite}\

  \begin{enumerate}
  \item If $G$ admits an irreducible admissible supercuspidal representation over some \emph{finite} field of characteristic $p$,
    then $G$ admits an irreducible admissible supercuspidal representation over any field of characteristic $p$.
  \item 
    If $G$ admits an irreducible admissible supercuspidal representation over some field of characteristic $c \ne p$, then
    $G$ admits an irreducible admissible supercuspidal representation over any algebraic extension of the prime field of characteristic $c$.
  \end{enumerate}
\end{proposition}

\begin{proof}
  Let $F_c$ be the prime field of characteristic $c$ (so that $F_0 = \mathbb{Q}$ and $F_c = \mathbb{F}_c$ if $c\neq 0$).  

  \step{}\label{step:sc-charl}
  We show that, if $c\neq p$ and $G$ admits an irreducible admissible supercuspidal $C$-representation $\pi$, then $G$
  admits one over a finite extension of $F_c$.
 
  Indeed, by Proposition \ref{prop HenV} we can suppose $C$ is algebraically closed.  We claim that we may twist $\pi$ by a $C$-character $\chi$ of $G$, so that the central character of $\pi\otimes\chi$ takes values in $\overline{F_c}$.  
To see this, we first note that there exists a subgroup ${}^\circ G$ of $G$ such that (1) $G/{}^\circ G \cong \mathbb{Z}^r$ for some $r\geq 0$; (2) the restriction to $Z(G)$ of the map $u: G \twoheadrightarrow \mathbb{Z}^r$ has image of finite index; (3) $\ker(u|_{Z(G)}) = Z(G) \cap {}^\circ G$ is compact. (For all of this, see \cite[\S 1.12, 2.3]{bernstein-centre}.)  Let $\mathfrak{L} := \textnormal{im}(u|_{Z(G)}) \subset \mathbb{Z}^r$ denote the image of $u|_{Z(G)}$.  Since $C$ is algebraically closed, the restriction map 
$$\Hom(\mathbb{Z}^r, C^\times) \xrightarrow{\textnormal{res}} \Hom(\mathfrak{L}, C^\times)$$ 
is surjective.  Let $\omega_\pi$ denote the central character of the irreducible admissible $C$-representation $\pi$, and note that $\omega_\pi|_{Z(G)\cap {}^\circ G}$ takes values in $\overline{F_c}$ (since $\pi$ is smooth and $Z(G)\cap {}^\circ G$ is compact).  Choose a splitting $v$ of the surjection $u:Z(G) \twoheadrightarrow \mathfrak{L}$, and let $\chi''\in \Hom(\mathfrak{L},C^\times)$ denote the character $\omega_\pi^{-1}\circ v$.  We then let $\chi'\in \Hom(\mathbb{Z}^r, C^\times)$ denote any preimage of $\chi''$ under $\textnormal{res}$, and let $\chi:G \rightarrow C^\times$ be the inflation of $\chi'$ to $G$ via $u$.  Using that $\omega_{\pi\otimes\chi} = \omega_\pi \chi$ and $\omega_{\pi\otimes\chi}|_{Z(G)\cap {}^\circ G} = \omega_\pi|_{Z(G) \cap {}^\circ G}$, the construction of $\chi$ implies $\omega_{\pi\otimes\chi}(z) \in \overline{F_c}$ for all $z\in Z(G)$.  
  
We may therefore assume that the central character of $\pi$ takes values in $\overline{F_c}$.  As $c \ne p$, by \cite[II.4.9]{Viglivre} the representation $\pi$ descends to a finite extension $F_c'/F_c$. Since descent preserves irreducibility, admissibility and supercuspidality, we obtain an irreducible admissible supercuspidal $F_c'$-representation of $G$.

  \step{}\label{step:sc-descend}
  We show that if $G$ admits an irreducible admissible supercuspidal representation over a finite extension
  of $F_c$ then $G$ admits such a representation over $F_c$.
 
Suppose $C/F_c$ is a finite field extension and $\pi$ an irreducible admissible $C$-representation of $G$. By Lemma \ref{lm:scalar-res}, $\pi|_{F_c[G]}$ contains an irreducible admissible $F_c$-representation $\pi'$.  By adjunction,  $\pi$ is a quotient of the scalar extension $\pi'_C$ of $\pi'$ from $F_c$ to $C$. 
  
  We now show that if $\pi$ is supercuspidal, then $\pi'$ is also supercuspidal.  Assume that $\pi'$ is not supercuspidal, so that it is a subquotient of $\Ind_P^G \tau'$, where $P$ is a proper parabolic subgroup of $G$ and $\tau'$ is an irreducible admissible $F_c$-representation of the Levi subgroup $M$ of $P$. Since parabolic induction is compatible with scalar extension from $F_c$ to $C$, the representation $\pi'_C$ is a subquotient of $\Ind_P^G \tau'_C$, and therefore the same is true of $\pi$. The $C$-representation $\tau'_C$ of $M$ has finite length and its irreducible subquotients are admissible by \cite[Thm.\ III.4]{HenV}. Hence, $\pi$ is a subquotient of $\Ind_P^G \rho$ for some irreducible admissible subquotient $\rho$ of $\tau'_C$, and we conclude that $\pi$ is not supercuspidal.
     
  \step{}\label{step:sc-basechange} 
  We show that if $G$ admits an irreducible admissible supercuspidal $\fp$-representation (resp., $F_c$-representation, where $c \neq p$), then $G$ does so over
  any field of characteristic $p$ (resp., any algebraic extension of $F_c$). More generally we show that if $L/C$ is any field extension, assumed to be algebraic if $c\neq p$, and $G$ admits an 
  irreducible admissible supercuspidal $C$-representation then the same is true over $L$.

  Let $L/C$ be a field extension as above, and choose compatible algebraic closures $\overline{L}/\overline{C}$.  Suppose $\pi$ is an irreducible admissible supercuspidal $C$-representation of $G$, and let $\tau$ be an irreducible subquotient of the scalar extension $\pi_L$ of $\pi$ from $C$ to $L$. By \cite[Lem.\ III.1(ii)]{HenV}, $\tau$ is admissible.  The scalar extension $\tau_{\overline{L}}$ of $\tau$ from $L$ to $\overline{L}$ is a subquotient of the scalar extension $\pi_{\overline{L}}$ of $\pi_L$ from $L$ to $\overline{L}$ (the latter being equal to the scalar extension of $\pi$ from $C$ to $\overline{L}$).  By Propositions \ref{prop fact HenV} and \ref{prop HenV}, $\pi_{\overline{L}}$ has finite length and its irreducible subquotients are admissible and supercuspidal.
  Therefore, the same is true of $\tau_{\overline{L}}$.  
  By Proposition \ref{prop HenV} this implies that $\tau$ is supercuspidal.
\end{proof}

We now use extension of scalars to prove the following lemma, which will be used in the proof of Prop.\ \ref{prop HG}.
 
\begin{lemma} \label{lem te} Let $\pi$ be an irreducible admissible $C$-representation of $G$ and $H$ a finite commutative
  quotient of $G$.  Then the representation $\pi \otimes_C C[H]$ of $G$, with the natural action of $G$ on $C[H]$, has finite
  length and its irreducible subquotients are admissible.
\end{lemma}

\begin{proof}  
  The scalar extension of the $C$-representation $\pi$ (resp.\ $C[H]$) to $\overline{C}$ has finite length with irreducible admissible quotients
  $\pi_i$ (resp.\ $\chi_j$, of dimension $1$). Therefore $(\pi \otimes_C C[H])_{\overline{C}}\cong \pi_{\overline{C}}
  \otimes_{\overline{C}}{\overline{C}}[H]$ has finite length with irreducible admissible subquotients (namely, the $\pi_i\otimes_{\overline{C}} \chi_j$),
  implying the same for $\pi \otimes_C C[H]$.
\end{proof}

\subsection{Reduction to an absolutely simple adjoint group}
\label{sec:reduct-an-absol}

As is well known, the adjoint group $\bG\ad$ of $\bG$ is $F$-isomorphic to a finite direct product of connected reductive $F$-groups 
\begin{equation}\label{eq Gad} 
\bG\ad\cong \bH \times \prod_i \Res_{F'_i/F}(\bG'_i),
\end{equation} 
where $\bH$ is anisotropic, the $F'_i/F$ are finite separable extensions, and $\Res_{F'_i/F}(\bG'_i)$ are scalar restrictions from $F'_i$ to $F$ of isotropic, absolutely simple, connected adjoint $F'_i$-groups $\bG'_i$.

\begin{proposition} \label{prop redadjoint} 
Assume that the field $C$ is algebraically closed or finite, and that $\chr F = 0$. 
If, for each $i$, the group $\bG'_i(F'_i)$ admits an irreducible admissible supercuspidal $C$-representation, then $G$ admits an irreducible admissible supercuspidal $C$-representation.
\end{proposition}

The proposition is the combination of Propositions \ref{prop GtimesH}, \ref{prop GHfinite}, \ref{prop HtoG}, and \ref{prop RF'toF} below, corresponding to the operations of finite product, central extension, and scalar restriction (all when $C$ algebraically closed or finite). We also note that if $\bG$ is anisotropic, then $G$ is compact and any irreducible smooth representation of $G$ is finite-dimensional (hence admissible) and supercuspidal.

\subsubsection{Finite product} \label{sssect finprod}

Let ${\mathbf G_1}$ and ${\mathbf G_2}$ be two connected reductive $F$-groups, and $\sigma$ and $\tau$ irreducible admissible $C$-representations of $G_1$ and $G_2$, respectively.  
 
\begin{proposition} \label{prop GtimesH} Assume that $C$ is algebraically closed.
  \begin{enumerate}
  \item The tensor product $\sigma \otimes_C\tau$ is an irreducible admissible $C$-representation of $G_1\times G_2$.
  \item Every irreducible admissible $C$-representation of $G_1\times G_2$ is of this form.
  \item The $C$-representation $\sigma \otimes_C\tau$ determines $\sigma$ and $\tau$ \(up to isomorphism\).
  \item The $C$-representation $\sigma \otimes_C\tau$ is supercuspidal if and only if $\sigma$ and $\tau$ are supercuspidal.
  \end{enumerate}
\end{proposition}
 
\begin{proof} 
Note first that $\sigma\otimes_C \tau$ is admissible: for compact open subgroups $K_1$ of $G_1$ and $K_2$ of $G_2$, we have a natural isomorphism (\cite[\S
  12.2 Lem.\ 1]{BkiA8})
  $$\Hom_{K_1}(\mathbf{1}_{K_1}, \sigma )\otimes_C  \Hom_{K_2}(\mathbf{1}_{K_2},  \tau) \congto  \Hom_{K_1\times K_2}(\mathbf{1}_{K_1}  \otimes_C \mathbf{1}_{K_2}, \sigma \otimes_C \tau),$$
  where $\mathbf{1}_{K_i}$ denotes the trivial representation of $K_i$.  Thus, the admissibility of $\sigma$ and $\tau$ implies the admissibility of $\sigma \otimes_C \tau$.

  Suppose now $C$ algebraically closed.
 
(i)  Proposition \ref{prop fact HenV} implies that the commutant of $\sigma$ is $C$.  Irreducibility of $\sigma\otimes_C \tau$ then follows from \cite[\S 12.2 Cor.\ 1]{BkiA8}.  

(ii) Let $\pi$ be an irreducible admissible $C$-representation of $G_1\times G_2$, and let $K_1,K_2$ be any compact open
  subgroups of $G_1, G_2$, respectively, such that $\pi^{K_1\times K_2}\neq 0$.

  If $c=p$, the $C$-representation of $G_1$ given by $\pi^{1\times K_2}$ is admissible (since $\pi^{K_1'\times K_2}$ is finite
  dimensional for any $K_1'$). By \cite[Lemma 7.10]{HV2}, it contains an irreducible admissible $C$-subrepresentation $\sigma$. Set
  $\tau:= \Hom_{G_1}(\sigma, \pi) \ne 0$, with the natural action of $G_2$. The representation $\sigma \otimes_C\tau$ embeds
  naturally in $\pi$. As $\pi$ is irreducible, it is isomorphic to $\sigma \otimes_C\tau$, and $\tau$ is irreducible. As
  $\pi$ is admissible, $\tau$ is admissible as well. (This proof is due to Henniart.)
 
  If $c\neq p$, the space $\pi^{K_1\times K_2}$ is a simple right $H_C(G_1\times G_2, K_1\times K_2)$-module (\cite[I.4.4, I.6.3]{Viglivre}), and we have 
 $$H_C(G_1\times G_2, K_1\times K_2)\cong H_C(G_1, K_1)\otimes_C H_C(G_2,K_2).$$
 By \cite[\S 12.1 Thm.\ 1]{BkiA8}, the finite-dimensional simple $H_C(G_1, K_1)\otimes_C H_C(G_2,K_2)$-modules factor, meaning $\pi^{K_1\times K_2}\cong \sigma^{K_1} \otimes_C \tau^{K_2}$ for irreducible admissible $C$-representations $\sigma, \tau$ of $G_1,G_2$, respectively (this uses \cite[I.4.4, I.6.3]{Viglivre} again).  Thus, we obtain $\pi \cong \sigma \otimes_C \tau$.

(iii) As a $C$-representation of $G_1$, $\sigma\otimes_C\tau$ is $\sigma$-isotypic.  Similarly, as a $C$-representation of $G_2$,  $\sigma\otimes_C\tau$ is $\tau$-isotypic. The result follows.

(iv) The parabolic subgroups of $G_1\times G_2$ are products of parabolic subgroups of $G_1$ and of $G_2$. Let $P, Q$ be
  parabolic subgroups of $G_1,G_2$, respectively, with Levi subgroups $M, L$, respectively, and let $\pi'$ be an irreducible admissible $C$-representation
  of the product $M\times L$.  By part (ii), the $C$-representation $\pi'$ factors, say $\pi' = \sigma' \otimes_C \tau'$ for irreducible admissible
  $C$-representations $\sigma'$ of $M$ and $\tau'$ of $L$.   We then obtain a natural isomorphism
$$\Ind_P^{G_1}\sigma' \otimes_C \Ind_Q^{G_2}\tau' \congto \Ind_{P\times Q}^{G_1\times G_2}\pi'.$$ 
  Since the inductions on the left-hand side have finite length, we conclude that the irreducible subquotients of $\Ind_{P\times Q}^{G_1\times G_2}\pi'$ are tensor products of the irreducible subquotients of $\Ind_P ^{G_1} \sigma'$ and of $\Ind_Q^{G_2} \tau'$, which gives the result.  
\end{proof}

We assume from now until the end of \S\ref{sssect finprod} that $C$ is a finite field.  

\begin{proposition}\label{prop GHfin}
  Assume that $C$ is finite.
  Let $\pi$ be an irreducible admissible $C$-representation of $G$. The commutant of $\pi$ is a finite field extension $D$ of $C$
  and the scalar extension $\pi_D$ of $\pi$ from $C$ to $D$ is isomorphic to
  $$\pi_D\cong \bigoplus _{i\in \Gal(D/C)}\pi_i,$$
  where the $\pi_i$ are irreducible admissible $D$-representations of $G$.  Moreover, the $\pi_i$ each have commutant $D$, are pairwise non-isomorphic, form a
  single $\Gal(D/C)$-orbit, and, viewed as $C$-representations, are isomorphic to $\pi$.
\end{proposition}

\begin{proof}
  The commutant $D$ of $\pi$ is a division algebra of finite dimension over $C$. Since the Brauer group of a finite field is trivial, $D$ is a finite Galois extension of $C$.  The result now follows from \cite[Thms.\ I.1, III.4]{HenV} by taking $R'=D$.  (Note also that as a $C$-representation, $\pi_D$ is $\pi$-isotypic of length $[D:C]$.)
\end{proof}

Recall that we have fixed irreducible admissible $C$-representations $\sigma$ and $\tau$ of $G_1$ and $G_2$, respectively.  Their respective commutants $D_\sigma$ and $D_\tau$ are finite extensions of $C$ of dimensions $d_\sigma$ and $d_\tau$, respectively. We embed them into $\overline{C}$, and consider:
\begin{itemize}
\item the field $D$ generated by $D_\sigma$ and $D_\tau$, which has $C$-dimension $\lcm(d_\sigma, d_\tau)$, 
\item the field $D' := D_\sigma\cap D_\tau$, which has $C$-dimension $\gcd(d_\sigma, d_\tau)$. 
\end{itemize}
The fields $D_\sigma, D_\tau$ are linearly disjoint over $D'$, we have $D_\sigma\otimes_{D'}D_\tau\cong D$ and 
\begin{equation}\label{eq tensor} 
D_\sigma \otimes_C D_\tau \cong \prod_{k = 1}^{[D':C]} D.
\end{equation}

\begin{proposition}\label{prop GHfinite}   
  Assume that $C$ is finite.
  The $C$-representation $\sigma\otimes_C \tau$ of $G_1\times G_2$ is isomorphic to
  $$\sigma\otimes_C \tau \cong \bigoplus_{k =1}^{\gcd(d_\sigma, d_\tau)} \ \pi_k,$$ 
where the $\pi_k$ are irreducible admissible $C$-representations with commutant $D$, which are pairwise non-isomorphic. The $C$-representations $\sigma$
  and $\tau$ are supercuspidal if and only if all the $\pi_k$ are supercuspidal, if and only if some $\pi_k$ is supercuspidal.
\end{proposition}

\begin{proof}
  By Proposition \ref{prop GHfin}, we have
  $$\sigma_D\cong \bigoplus _{i\in \Gal(D_\sigma /C)}\sigma_i,\qquad \tau_D\cong \bigoplus _{j\in \Gal(D_\tau /C)}\tau_j,$$
  where the $\sigma_i$ (resp.\ $\tau_j$) are irreducible admissible $D$-representations of $G_1$ (resp.\ $G_2$) with commutant $D$, which are pairwise non-isomorphic, form a single $\Gal(D/C)$-orbit, descend to $D_\sigma$ (resp.\ $D_\tau$) and their descents, viewed as $C$-representations, are isomorphic to $\sigma$ (resp.\ $\tau$). The $C$-representation $\sigma\otimes_C\tau$ of $G_1\times G_2$ is admissible, and its scalar extension from $C$ to $D$ is equal to
  \begin{equation}\label{eq tensor1}( \sigma\otimes_C \tau)_D\ \cong \ \sigma_D\otimes_D
    \tau_D \cong \bigoplus_{(i,j)\in \Gal (D_\sigma/C)\times \Gal (D_\tau/C)}\ \sigma_i \otimes_D \tau_j.
  \end{equation}
  The $D$-representation $\sigma_i \otimes_D \tau_j$ of $G_1\times G_2$ is admissible and has commutant $D\otimes_D D=D$ (\cite[\S12.2 Lem.\ 1]{BkiA8}). Hence, $\sigma_i \otimes_D \tau_j$ is absolutely irreducible and equation \eqref{eq tensor1} implies $(\sigma\otimes_C \tau)_D$ is semisimple. By \cite[\S12.7 Prop.\ 8]{BkiA8}, this implies that $\sigma\otimes_C \tau$ is semisimple; its commutant is isomorphic to $D_\sigma \otimes_C D_\tau$ by \cite[\S12.2 Lem.\ 1]{BkiA8}.
  From equation \eqref{eq tensor} we see that $\sigma\otimes_C \tau$ has length $[D':C] = \gcd (d_\sigma, d_\tau)$, its irreducible constituents $\pi_k$ are admissible and pairwise non-isomorphic with commutant $D$.

  Applying Proposition \ref{prop GtimesH} over $\overline{C}$ and Proposition \ref{prop HenV} (several times), we see that $\sigma$ and $\tau$ are supercuspidal if and only if some/every $\sigma_i$ and some/every
  $\tau_j$ are supercuspidal, if and only if some/every $\sigma_i\otimes_D\tau_j$ is supercuspidal. From Proposition \ref{prop HenV} again, this is also equivalent to $\pi_k$ being supercuspidal for some/every $k$.
\end{proof}

\subsubsection{Central extension} \label{sssect centralext}

Recall that we have a short exact sequence of $F$-groups
$$1 \rightarrow \mathbf{Z}(\mathbf{G}) \rightarrow \mathbf{G}  \xrightarrow{i} \mathbf{G}^{\rm ad} \rightarrow 1,$$
which induces an exact sequence between $F$-points
$$1\rightarrow Z(G) \rightarrow  G \xrightarrow{i} G^{\rm ad} \rightarrow H^1(F, \mathbf{Z(G)} ).$$
The image $i(G)$ of $G$ is a closed cocompact normal subgroup of $G^{\rm ad}$ and $H^1(F, \mathbf{Z(G)} )$ is commutative.  

Until the end of \S\ref{sssect centralext}, we assume that $\chr F = 0$. The group $H^1(F, \mathbf {Z(G)} )$ is then finite (\cite[Thm.\ 6.14]{PlaRap}), implying that
$i(G)$ is an \emph{open} normal subgroup of $G^{\rm ad}$ and the quotient $G^{\rm ad}/i(G)$ is finite and commutative.  Our next task will be to prove the following:

\begin{proposition} \label{prop HtoG} 
Suppose that $\chr F = 0$. Then $G^{\rm ad}$ admits an irreducible admissible supercuspidal $C$-representation if and only if $G$ admits such a representation such that moreover $Z(G)$ acts trivially.
\end{proposition}
 
Inflation from $i(G)$ to $G$ identifies representations of $i(G)$ with representations of $G$ having trivial $Z(G)$-action; this inflation functor respects irreducibility and admissibility.  The composite functor 
\begin{center}
(inflation from $i(G)$ to $G$) $\circ$ (restriction from $G^{\rm ad}$ to $i(G)$)
\end{center}
from $C$-representations of $G^{\rm ad}$ to representations of $G$ trivial on $Z(G)$ will be denoted by $- \circ i$.

Suppose $\widetilde{\rho}$ is an irreducible admissible $C$-representation of $G$ with trivial action of $Z(G)$.  Then $\widetilde{\rho}$ is the inflation of a representation $\rho$ of the open, normal, finite-index subgroup $i(G)$ of $G^{\rm ad}$. The $C$-representation $\rho$ of $i(G)$ is irreducible and admissible, and therefore the induced representation $\Ind_{i(G)}^{G^{\rm ad}}\rho$ 
of $G^{\rm ad}$ is admissible of finite length.  
Any irreducible quotient $\pi$ of $\Ind_{i(G)}^{G^{\rm ad}}\rho$ is admissible (if $c = p$, this uses the assumption $\chr F = 0$; see \cite[\S 4, Thm.\ 1]{henniart:adm}).  
By adjunction, $\pi|_{i(G)}$ contains a subrepresentation isomorphic to $\rho$ and, by inflation from $i(G)$ to $G$, $\widetilde{\rho}$ is isomorphic to a subquotient of $\pi \circ i$.
  
Conversely, suppose $\pi$ is an irreducible admissible $C$-representation of $G^{\rm ad}$.  The restriction $\pi|_{i(G)}$ of $\pi$ to $i(G)$ is semisimple of finite length, and its irreducible constituents $\rho$ are $G^{\rm ad}$-conjugate and admissible (see \cite[I.6.12]{Viglivre}; note that the condition that the index is invertible in $C$ is not
  necessary and not used in the proof).  Hence, the $C$-representation $\pi \circ i$ of $G$ is semisimple of finite length, and its irreducible constituents are the inflations $\widetilde{\rho}$ of the irreducible constituents $\rho$ of $\pi|_{i(G)}$.

Proposition \ref{prop HtoG} now follows from:

\begin{proposition}\label{prop HG} 
Suppose that $\chr F = 0$ and let $\pi, \rho$ and $\widetilde{\rho}$ be as above.  Then $\pi$ is supercuspidal if and only if some $\widetilde{\rho}$ is supercuspidal, if and only if
  all $\widetilde{\rho}$ are supercuspidal.
\end{proposition}

\begin{proof}
We first check first the compatibility of parabolic induction with $- \circ i$.  The parabolic $F$-subgroups of $\bf G$ and of $\bf G^{\rm ad}$ are in bijection via the map $i$ (\cite[22.6 Thm.]{BorelLAG}). If the parabolic $F$-subgroup $\bf P$ of $\bf G$ corresponds to the parabolic $F$-subgroup $\bf Q$ of $\bf G^{\rm ad}$, then $i$ restricts to an isomorphism between their unipotent radicals, and sends a Levi subgroup $\bf M$ of $\bf P$ onto a Levi subgroup $\bf L$
  of $\bf Q$.  Further, we have an exact sequence between $F$-points:
  $$1\rightarrow Z(G) \rightarrow M \xrightarrow{i} L \rightarrow H^1(F, \bZ(\bG) ).$$
  We have $G^{\rm ad}=Q i(G)$ and $Q\cap i(G)=i(P)=i(M)U,$ where $i(M)$ is an open normal subgroup of $L$ having finite commutative
  quotient, and $U$ is the unipotent radical of $Q$. Thus, if $\sigma$ is a smooth $C$-representation of $L$, the Mackey decomposition implies
  $(\Ind_Q^{G^{\rm ad}}\sigma)|_{i(G)} \cong \Ind_{i(P)}^{i(G)}(\sigma|_{i(M)})$ and, by inflation from $i(G)$ to $G$, we obtain:
  \begin{equation}\label{eq ind}
  (\Ind_Q^{G^{\rm ad}}\sigma) \circ i\cong \Ind_P^G (\sigma\circ i).
  \end{equation}

We may now proceed with the proof.  It suffices to prove: 
\begin{enumerate}
\item if $\pi$ is non-supercuspidal, then all $\widetilde{\rho}$ are non-supercuspidal, 
\item if some $\widetilde{\rho}$ is non-supercuspidal, then $\pi$ is non-supercuspidal. 
\end{enumerate}

To prove (i), let $\pi$ be an irreducible admissible non-supercuspidal $C$-representation of $G^{\rm ad}$, which is isomorphic to a subquotient of $\Ind_Q^{G^{\rm ad}} \sigma$ for $Q \subsetneq G^{\rm ad}$ and $\sigma$ an irreducible admissible $C$-representation of $L$.  Therefore, $\pi\circ i$ is isomorphic to a subquotient of $(\Ind_Q^{G^{\rm ad}} \sigma) \circ i$, and by equation \eqref{eq ind}, each $\widetilde{\rho}$ is isomorphic to a subquotient of $\Ind_P^G \widetilde{\tau}$ for some irreducible subquotient $\widetilde{\tau}$ of $\sigma \circ i$ (depending on $\rho$).  Since $\widetilde{\tau}$ is admissible and $P\subsetneq G$, all the $\widetilde{\rho}$ are non-supercuspidal.

To prove (ii), let $\pi$ be an irreducible admissible $C$-representation of $G^{\rm ad}$ such that some irreducible constituent $\widetilde{\rho}$ of $\pi\circ i$ is non-supercuspidal.  Suppose $\widetilde{\rho}$ is isomorphic to a subquotient of $\Ind_P^G\tau'$ for $P\subsetneq G$ and $\tau'$ an irreducible admissible $C$-representation of $M$.  The central subgroup $Z(G)$ acts trivially on $\widetilde{\rho}$, and hence also on $\tau'$.  Therefore $\tau'= \widetilde{\tau}$ for some irreducible subquotient $\tau$ of $\sigma|_{i(G)}$, where
  $\sigma$ is an irreducible admissible $C$-representation of $L$. The representation $\widetilde{\rho}$ is isomorphic to a subquotient of $\Ind_P^{G}(\sigma \circ i)$. 
By equation \eqref{eq ind} and exactness of parabolic induction, $\Ind_{i(G)}^{G\ad}(\rho)$, and hence its quotient $\pi$, is isomorphic to a subquotient of
  $\Ind_{i(G)}^{G^{\rm ad}}((\Ind_Q^{G^{\rm ad}}\sigma)|_{i(G)})$. This representation is isomorphic to
$$\Ind_{ i(M)U}^{G^{\rm ad}}(\sigma|_{i(M)}) \cong \Ind_Q^{G^{\rm ad}}\big(\Ind_{ i(M)}^{L}(\sigma|_{i(M)})\big) \cong \Ind_Q^{G^{\rm ad}} (\sigma \otimes_C C[i(M)\backslash L]).$$
By Lemma \ref{lem te}, the $C$-representation $\sigma\otimes_C C[i(M)\backslash L]$ of $L$ has finite length and its irreducible subquotients
  $\nu$ are admissible. Therefore $\pi$ is isomorphic to a subquotient of $\Ind_Q^{G^{\rm ad}} \nu$ for some
  $\nu$ and some $Q\subsetneq G^{\rm ad}$, and therefore $\pi$ is non-supercuspidal.
\end{proof}
     
\subsubsection{Scalar restriction} 
\label{sec:scalar-restriction}

Now let $F'/F$ be a finite separable extension, $\bf G'$ a connected reductive $F'$-group and ${\bf G} := \textnormal{Res}_{F'/F}({\bf G'})$ the scalar restriction of $\bf G'$ from $F'$ to $F$.  As topological groups, $G' := {\bf G'}(F')$ is equal to $G := {\bf G}(F)$.  By \cite[6.19. Cor.]{BorelTits}, $G'$ and $G$ have the same parabolic
subgroups.  Hence:

\begin{proposition}\label{prop RF'toF} 
  $G'$ admits an irreducible admissible supercuspidal $C$-representation if and only if $G$ does.
\end{proposition}

\section{Proof of the main theorem for most simple groups}
\label{sec:proof-main-theorem}

\subsection{Discrete Iwahori--Hecke modules}
\label{sec:discr-iwah-hecke}

Let $\Rep_C(G,\mathfrak B)$ denote the category of $C$-representations of $G$ generated by their $\mathfrak B$-invariant vectors, and let $\Mod(H_C(G,\mathfrak B))$ denote the category of right $H_{C}(G,\mathfrak B)$-modules.  The functor of $\mathfrak{B}$-invariants 
\begin{align*}
\Rep_C(G,\mathfrak B) & \rightarrow \Mod(H_C(G,\mathfrak B))\\
\pi & \mapsto \pi^{\mathfrak B}
\end{align*}
admits a left adjoint
\begin{align*}
\mathfrak T:  \Mod(H_C(G,\mathfrak B)) & \rightarrow \Rep_C(G,\mathfrak B)\\
M & \mapsto M \otimes_{ H_{ C}(G,\mathfrak B)} C[\mathfrak B\backslash G].
\end{align*}

\begin{prop}\label{fact1}
  When $c\neq p$, the functor $\pi \mapsto \pi^{\mathfrak B}$ induces a bijection between the isomorphism classes of irreducible $C$-representations $\pi$ of $G$ with $\pi^{\mathfrak B} \neq 0$ and isomorphism classes of simple right $H_C(G,\mathfrak{B})$-modules \(\cite[I.4.4, I.6.3]{Viglivre}\).  When $C=\mathbb{C}$, the functors are inverse equivalences of categories \(cf.\ \cite[Cor.\ 3.9(ii)]{bernstein-centre}; see also \cite[Thms.\ 4.8, 4.4(iii)]{morris}\). 
\end{prop}

\begin{remark}
The above functors are not as well-behaved when $c = p$. In this case, the functor of $\mathfrak{B}$-invariants may not preserve irreducibility.  Similarly, the left adjoint $\mathfrak{T}$ may not preserve irreducibility.
\end{remark}

When $C = \mathbb C$, the Bernstein ring embedding $H_{\mathbb{C}}(Z,Z_0) \xrightarrow{\theta} \widetilde{H}_\mathbb{C}$ is the linear map defined by 
sending $T^Z_\lambda$ to $\theta_\lambda := q_\lambda^{-1/2}E_\lambda$ 
for $\lambda \in \Lambda$.  Its image is equal to $\mathcal{A}_{\mathbb{C}}$.  Note that if $\lambda\in \Lambda$ is anti-dominant and $z\in Z$ lifts $\lambda$, we have $q_\lambda = \delta_B(z)$, where $\delta_B$ denotes the modulus character of $B$.

We now recall some properties of the category $\Rep_{\mathbb{C}}(G,\mathfrak B)$, including Casselman's criterion of square integrability modulo center, before giving the definition of a discrete simple right $H_{\mathbb{C}}(G,\mathfrak B)$-module.  Recall that $\pi_U$ denotes the space of $U$-coinvariants (i.e., the unnormalized Jacquet module) of a representation $\pi$.

\begin{lm}\label{lm:pro-p-invts-jacquet}
  Suppose that $\pi$ is an admissible $\mathbb{C}$-representation of $G$.  Then the natural map $\pi \onto \pi_U$ induces an isomorphism $\vp: \pi^{\mathfrak{B}} \congto \pi_U^{Z_0}$.
  Moreover, we have $\vp(v\cdot \theta_{\lambda^{-1}}) = \delta_B^{-1/2}(t) (t\cdot\vp(v))$ for $\lambda \in \Lambda_T$, $t \in T$ lifting $\lambda$, and $v \in \pi^{\mathfrak{B}}$.
\end{lm}

\begin{proof}
  Recall that $\mathfrak{B}$ has an Iwahori decomposition with respect to $Z$, $U$, $U\op$. Then \cite[Prop.\ 4.1.4]{Cas} implies that the map $\pi \onto \pi_U$ induces an isomorphism 
  $\pi^{\mathfrak{B}} \cdot T_{\lambda^{-1}} \congto \pi_U^{Z_0}$ for $\lambda \in \Lambda_T$
  with $\max_{\alpha\in \Delta} |\alpha(\lambda)|_F$ sufficiently small. By \cite[Prop.\ 4.13(1)]{Vig1} the operator
  $T_{\lambda^{-1}}$ is invertible in $H_{\C}(G,\mathfrak{B})$, so $\pi^{\mathfrak{B}} \cdot T_{\lambda^{-1}} = \pi^{\mathfrak{B}}$. 

  To show the last statement, we may assume that $\lambda \in \Lambda_T^+$.
  Then, in our terminology, \cite[Lemma 4.1.1]{Cas} says that $\vp(|\B t\B/\B|^{-1} [\B t\B]\cdot v) = t\cdot \vp(v)$, where $[\B t\B]$ denotes the usual double
  coset operator on $\pi^\B$. Now $[\B t\B]\cdot v = v\cdot T_{t^{-1}}$ and $T_{t^{-1}} = E_{t^{-1}} = q_{t^{-1}}^{1/2} \theta_{t^{-1}}$.
  Moreover, $|\B t\B/\B| = q_t = q_{t^{-1}} = \delta_B(t^{-1})$. Putting this all together, we obtain the claim.
\end{proof}

\begin{remark}\label{rk:pro-p-invts-jacquet}
The lemma and its proof hold when $\mathfrak{B}$ is replaced by $\mathfrak{U}$ and $Z_0$ is replaced by $Z_0 \cap \mathfrak{U}$.  
\end{remark}

\begin{prop}\label{prop fact2} 
Let $\pi$ be an irreducible $\mathbb{C}$-representation of $G$ with $\pi^{\mathfrak B} \neq 0$. 
\begin{enumerate}
\item$\pi$ is isomorphic to a subrepresentation of $\Ind_B^G\sigma$, where $\sigma$ is a $\mathbb{C}$-character of $Z$
  trivial on $Z_0$.
\item \textbf{Casselman's criterion:} $\pi$ is square integrable modulo center \(as defined in \cite[\S 2.5]{Cas}\) if and only if its central
  character is unitary and 
$$|\chi (\mu(\varpi))|_{\mathbb{C}}  < 1$$ 
for all $\mu \in X_*(T)^+$ such that $-\mu \not\in X_*(T)^+$, and all characters $\chi$ of $T$ 
contained in $\delta_B^{-1/2}\pi _U$.
\end{enumerate}
\end{prop}

\begin{proof}
(i) Since $\pi$ is irreducible and smooth, it is admissible by \cite[II.2.8]{Viglivre}, and \cite[3.3.1]{Cas} implies $\pi_U$ is admissible as well.  By Lemma \ref{lm:pro-p-invts-jacquet} and the assumption $\pi^{\mathfrak{B}}\neq 0$, we see that $\pi_U\neq 0$.  The claim now follows by choosing an irreducible quotient $\pi_U \rightarrow \sigma$ for which $\sigma^{Z_0} \neq 0$ and applying Frobenius reciprocity.

(ii) This follows from \cite[Thm.\ 6.5.1]{Cas}.
\end{proof}

\begin{definition} 
We say a simple right $H_ {\mathbb C}(G,\mathfrak B)$-module is \emph{discrete} if it is isomorphic to $\pi^{\mathfrak B}$ for
  an irreducible admissible square-integrable modulo center $\mathbb C$-representation $\pi$ of $G$.  We say a semisimple right $H_ {\mathbb C}(G,\mathfrak B)$-module is discrete if its simple subquotients are discrete.  
\end{definition}

\begin{proposition}\label{prop dis} 
A simple right $H_ {\mathbb C}(G,\mathfrak B)$-module $M$ is discrete if and only if any $\mathbb{C}$-character $\chi$ of $\mathcal{A}_{\mathbb{C}}$ contained in $M$ satisfies the following condition: the restriction of $\chi$ to $\Lambda_{Z(G)}$ is a unitary character, and
\begin{equation}
  |\chi  (\theta_{\lambda_\mu^{-1}})|_{\mathbb{C}} < 1\label{eq:10} 
\end{equation}
for  any $\mu \in X_*(T)^+$ such that $-\mu\not\in X_*(T)^+$.
\end{proposition}

\begin{proof}
  Note that $M = \pi^\B$ for an irreducible (admissible) $\C$-representation $\pi$ of $G$. Then $\pi$ has unitary central character
  if and only if $\Lambda_{Z(G)}$ acts by a unitary character on $M$. As any irreducible $\mathcal A_{\C}$-module is a character,
  by Casselman's criterion (Proposition \ref{prop fact2}) and Lemma~\ref{lm:pro-p-invts-jacquet}, $M$ is discrete if and only 
  condition~(\ref{eq:10}) holds.
\end{proof}

\begin{remark}\label{rem danger} 
Some authors view $\pi^\mathfrak B$ as a left $H_ {\mathbb C}(G,\mathfrak{B})$-module.  One may pass between left and right modules by using the anti-automorphism $T_w \mapsto T_{w^{-1}}$; that is, we may define 
$$T_w \cdot v = v\cdot T_{w^{-1}}$$ 
for $w\in \widetilde{W},  v \in \pi^\mathfrak{B}$.  The space $\pi^\mathfrak B$, viewed as either a left or right $H_ {\mathbb C}(G,\mathfrak B)$-module, is then called discrete if $\pi$ is square integrable modulo center. For left modules, the proposition above holds with ``dominant'' replaced by ``anti-dominant,'' and ``$\theta$'' replaced by ``$\widetilde{\theta}^+$'' (for the definition of $\widetilde{\theta}^+$, see the paragraph preceding Proposition 8 in \cite{Vig-mathann}).
\end{remark}

\begin{lemma} \label{lem chi} 
For a character $\chi:\mathcal{A}_{\mathbb{C}}\to \mathbb C$ such that $\chi|_{\Lambda_{Z(G)}}$ is unitary, the following conditions are equivalent:
\begin{enumerate}
\item $|\chi (\theta_{\lambda_\mu^{-1}})|_{\mathbb{C}} < 1$ for any $\mu \in X_*(T)^+$ such that $-\mu \not\in X_*(T)^+$,
\item $|\chi (\theta_{\lambda^{-1}})|_{\mathbb{C}} < 1$ for any $\lambda \in \Lambda^{\textnormal{sc},+}$ such that $\lambda^{-1} \not\in \Lambda^{\textnormal{sc},+}$,
\item $|\chi (\theta_{\lambda^{-1}})|_{\mathbb{C}} < 1$ for any $\lambda \in \Lambda^{+}$ such that $\lambda^{-1} \not\in \Lambda^{+}$.
  \end{enumerate}
 \end{lemma}

\begin{proof}
We first recall that the invertible elements in $\Lambda^+$ consist of $\Lambda_{\ker\nu}$, so $|\chi(\theta_\lambda)|_{\mathbb{C}} = 1$ for all invertible elements of $\Lambda^+$.

As $\Lambda_T \cong X_*(T)$, we see that (iii) implies (i) and (ii). To prove that (ii) implies (iii), 
we need to show that $|\chi (\theta_{\lambda^{-1}})|_{\mathbb{C}} = 1$ for $\lambda \in \Lambda^+$
implies $\lambda^{-1} \in \Lambda^+$. By Lemma \ref{lem fin} pick $n \ge 1$ such that $\lambda^n \in \Lambda_{Z(G)} \times \Lambda^{\textrm{sc},+}$. 
Then $\lambda^n \lambda_0 \in \Lambda^{\textrm{sc},+}$ for some $\lambda_0 \in \Lambda_{Z(G)}$.
As $|\chi (\theta_{\lambda^{-n} \lambda_0^{-1}})|_{\mathbb{C}} = 1$ we deduce
from (ii) that $\lambda^n \lambda_0 \in \Lambda^{\textrm{sc},+} \cap (\Lambda^{\textrm{sc},+})^{-1}$, which is contained in $\Lambda^+ \cap (\Lambda^+)^{-1}$. 
Therefore $\lambda^n \in \Lambda^+ \cap (\Lambda^+)^{-1}$.
From the definition of dominance it follows that $\lambda \in \Lambda^+ \cap (\Lambda^+)^{-1}$.

The proof that (i) implies (iii) is similar but easier.
\end{proof}

\begin{proposition} \label{prop dis2} 
A simple right $H_ {\mathbb C} (G, \mathfrak B)$-module $M$ is discrete if and only if
  $\Lambda_{Z(G) }$ acts on $M$ by a unitary character and if its restriction to $H_ {\mathbb C}(G\sc, \mathfrak B\sc)$
  is discrete.
\end{proposition}

\begin{proof}
  This follows from Proposition \ref{prop dis} and Lemma \ref{lem chi}.
\end{proof}

\subsection{Characters}
\label{sec:characters}

In this section we continue to assume $C$ is a field of characteristic $c$, and suppose further that $\bf G$ is absolutely simple and isotropic. We determine the characters $H = H(G\sc, \mathfrak{B}\sc) \to C$ which extend to $\widetilde{H} = H(G,\mathfrak{B})$. This is an exercise, which is already in the literature when $C=\mathbb C$ (cf.\ \cite{Borel}).

For distinct reflections $s,t\in S$, the order $n_{s,t}$ of $st$ is finite, except if the type of $\Sigma$ is $\textA_1$. In the finite case, the braid relations \eqref{eq braid} imply
\begin{alignat}{2} 
\label{eq TsTt 1} (T_sT_t)^r & =(T_tT_s)^r & &\quad \text{if}~ n_{s,t} = 2r,\\ 
\label{eq TsTt 2} (T_sT_t)^rT_s & = (T_tT_s)^rT_t &  &\quad  \text{if}~ n_{s,t} = 2r+1.
\end{alignat} 
The $T_s$ for $s\in S$ and the relations \eqref{eq Tsq}, \eqref{eq TsTt 1} and \eqref{eq TsTt 2} give a presentation of $H$. A
presentation of $\widetilde{H}$ is given by the $T_u, T_s$ for $u\in \Omega, s\in S$ and the relations \eqref{eq Tsq}, \eqref{eq TsTt 1}, \eqref{eq TsTt 2} and
\begin{alignat}{2}
\label{eq TsTu 1} T_uT_{u'} & =T_{uu'} & & \quad \text{if}~ u,u'\in \Omega,\\
\label{eq TsTu 2} T_uT_s & = T_{u(s)}T_u & &  \quad \text{if} \ u\in \Omega, s\in S,
 \end{alignat} 
where $u(s)$ denotes the action of $\Omega$ on $S$.

 We have a disjoint decomposition 
 $$S = \bigsqcup_{i=1}^m S_i,$$ 
 where $S_i$ is the intersection of $S$ with a conjugacy class of $W$. The $S_i$ are precisely the connected components of $\mathsf{Dyn}$ when all multiple edges are removed (see \cite[VI.4.3 Th.\ 4]{BkiGAL6} and \cite[3.3]{Borel}).  Thus, we have
$$m = \begin{cases} 1 & \\ 2 & \\ 3 & \\ \end{cases} \quad  \text{when the type of}~ \Sigma = \begin{cases} \textA_\ell~ (\ell \geq 2), \textD_{\ell}~ (\ell \geq 4), \textE_6, \textE_7, \textrm{or}~ \textE_8; & \\ \textA_1, \textB_{\ell}~ (\ell \geq 3), \textF_4, \textrm{or}~ \textG_2; & \\ \textC_\ell~ (\ell \geq 2). & \end{cases}$$
 When $m>1$, we fix a labeling of the $S_i$ such that $|S_1|\geq |S_2|$, and when the type of $\Sigma$ is $\textC_\ell~ (\ell
\geq 2)$, we let $S_2=\{s_2\}$ and $S_3=\{s_3\}$ denote the endpoints of $\mathsf{Dyn}$. (Note that there are two possible labelings in types $\textA_1$ and $\textC_\ell~ (\ell \geq 2)$.)
The parameters $d_s$ are equal on each component $S_i$; we denote this common value by $d_i$.

\begin{definition}
The unique character $\chi:H \to C$ with $\chi (T_s) = q_s$ (resp.,\ $\chi (T_s) = -1$) for all $s\in S$ is called the \emph{trivial} (resp.,\ \emph{special}) $C$-character. 
 \end{definition}

\begin{lemma} \label{lem Bor} 
Suppose $\{T_s\}_{s\in S}\to C$ is an arbitrary map.
\begin{enumerate}
\item When $c \neq p$, the above map extends to a character of $H$ if and only if it is constant on each $S_i$, and takes the value $-1$ or $q^{d_i}$ on each $T_s, s\in S_i$. There are $2^m$ characters if $q^{d_i} + 1\neq 0$ in $C$ for each $i$. \label{lem Bor-1}
\item When $c = p$, the above map extends to a character of $H$ if and only if its values are $-1$ or $0$ on each $T_s, s\in S$. There are $2^{|S|}$ characters. \label{lem Bor-2}
Such a character is supersingular if and only if it is not special or trivial. 
\end{enumerate}
\end{lemma}

\begin{proof}
(i) This follows from the presentation of $H$ and the fact that the $T_w$ are invertible (so that the map must be constant on conjugacy classes).  

(ii) This follows from \cite[Prop.\ 2.2]{Vigss}. The claim about supersingularity follows from \cite[Thm.\ 6.15]{Vigss}.
\end{proof}

We wish to determine which characters of $H$ extend to $\widetilde{H}$.  Since the elements $T_u$ for $u\in \Omega$ are invertible in $\widetilde{H}$, the relations \eqref{eq TsTu 2} imply that a character $\chi:H \to C$ extends to a character of $\widetilde H$ if and only if $\chi (T_s)=\chi (T_{u(s)})$ for all $s\in S$ and $u\in \Omega$.  For example, if the image $\Psi$ of $\Omega$ in $\Aut(W,S,d_i)$ is trivial, then any character of $H$ extends to $\widetilde H$. The extensions are not unique in general.  By their very definition, the trivial and special characters always extend, and we also refer to their extensions as trivial and special characters.

Let $\chi:H \rightarrow C$ denote a character, and suppose $c \neq p$.  By Lemma \ref{lem Bor}\ref{lem Bor-1}, the value of $\chi$ on $T_s$ for $s\in S_i$ is constant for each $1 \leq i \leq m$.  We 
define $\chi_i := \chi (T_s) \in C$ for $s\in S_i$, and identify the character $\chi$ with the $m$-tuple $(\chi_i)_{1\leq i \leq m}$.
   
\begin{lemma} \label{lem ext}
Assume $c\neq p$.  Let $\chi: H \to C$ denote a character of $H$, associated to the $m$-tuple $(\chi_i)_{1\leq i \leq m}$.  Then $\chi$ extends to a
  character of $\widetilde{H}$ except in the following cases:
\begin{itemize}
 \item type $\textA_1$, equal parameters $d_1 = d_2$, $\Psi \neq 1$, and $\chi_1\neq \chi_2$;
 \item type $\textC_\ell~ (\ell \geq 2)$, equal parameters $d_2 = d_3$, $\Psi \neq 1$, and $\chi_2\neq \chi_3$.
\end{itemize}
\end{lemma}

\begin{proof}
When $m = 1$, then $\chi(T_s) = \chi(T_{u(s)})$ for all $u\in \Omega$ and $s\in S$, so that $\chi$ extends to $\widetilde{H}$.  We may therefore assume $m > 1$.  We proceed type-by-type:
\begin{itemize}
\item \textit{Type $\textA_1$ with equal parameters $d_1= d_2$.}  The group $\Aut(W, S, d_s) \cong \mathbb{Z}/2\mathbb{Z}$ permutes $s_1$ and $s_2$. If $\Psi = 1$ or $\chi_1 = \chi_2$, then $\chi$ extends to $\widetilde{H}$, while if $\Psi \neq 1$ and $\chi_1 \neq \chi_2$, the character $\chi$ cannot extend.
\item \textit{Type $\textB_\ell~ (\ell \geq 3)$.} In this case, $\Aut(W, S, d_s) \cong \mathbb{Z}/2\mathbb{Z}$ stabilizes the sets $S_1$ and $S_2$, so that $\chi(T_s) = \chi(T_{u(s)})$ for all $u\in \Omega$ and $s\in S$.  Thus $\chi$ extends to $\widetilde{H}$.
\item \textit{Type $\textC_\ell~ (\ell \geq 2)$ with equal parameters $d_2 = d_3$.} The group $\Aut(W, S, d_s)\cong \mathbb{Z}/ 2 \mathbb{Z}$ permutes $s_2$ and $s_3$.  If $\Psi = 1$ or if $\chi_2 = \chi_3$, then $\chi$ extends to $\widetilde{H}$, while if $\Psi \neq 1$ and $\chi_2 \neq \chi_3$, the character $\chi$ cannot extend.
\item \textit{Type $\textA_1$ with unequal parameters $d_1\neq d_2$; Type $\textF_4$; Type $\textG_2$; Type $\textC_\ell~ (\ell \geq 2)$ with unequal parameters $d_2 \neq d_3$.} In these cases, $\Aut(W, S, d_s)$ (and consequently $\Psi$) is trivial, and thus $\chi$ extends to $\widetilde{H}$.
\end{itemize}
\end{proof}

Before stating the next result, we require a definition.

\begin{definition} 
  Let $R \subset \mathbb{C}$ be a subring of $\mathbb{C}$.  We say a right $\widetilde{H}_{\mathbb{C}}$-module $M$ is \emph{$R$-integral} if there exists an $\widetilde{H}_R$-submodule $M^\circ \subset M$ such that the natural map 
 $$\mathbb{C} \otimes_R M^\circ \rightarrow M$$ 
is an isomorphism of $\widetilde{H}_{\mathbb{C}}$-modules.  We call $M^\circ$ an \emph{R-integral structure} of $M$.  If $\mathfrak{p}$ is a maximal ideal of $R$, the $\widetilde{H}_{R/\mathfrak{p}}$-module $R/\mathfrak{p} \otimes_R M^\circ$ is called \emph{reduction of $M^\circ$ modulo $\mathfrak{p}$}.  We make similar definitions for the algebra $H_{\mathbb{C}}$.  
\end{definition}

The following proposition combines the above results.

\begin{proposition} \label{prop cha} \hfill
\begin{enumerate}
\item $H_\mathbb C$ admits $2^m$ $\mathbb{C}$-characters.  They are all $\mathbb Z$-integral,
  and their reductions modulo $p$ are supersingular except for the special and trivial characters.
\item Suppose $\chi:H_{\mathbb{C}} \rightarrow \mathbb{C}$ is a character, associated to the $m$-tuple $(\chi_i)_{1 \leq i \leq m}$, and suppose we are in one of the following two cases:
\begin{itemize}
 \item type $\textA_1$, equal parameters $d_1 = d_2$, $\Psi \neq 1$, and $\chi_1\neq \chi_2$;
 \item type $\textC_\ell~ (\ell \geq 2)$, equal parameters $d_2 = d_3$, $\Psi \neq 1$, and $\chi_2\neq \chi_3$.
\end{itemize}
Then the $H_{\mathbb{C}}$-module $\chi\oplus \overline{\chi}$ extends to a two-dimensional, $\mathbb{Z}$-integral simple \(left or right\) $\widetilde{H}_{\mathbb{C}}$-module with supersingular
  reduction modulo $p$, where $\overline{\chi} = (\chi_2, \chi_1)$ in the $\textA_1$ case and $\overline{\chi} = (\chi_1, \chi_3, \chi_2)$ in the $\textC_\ell$ case.  
\item Suppose $\chi:H_{\mathbb{C}} \rightarrow \mathbb{C}$ is a character which does not fall into either of the two cases of the previous point.  Then $\chi$ extends to a $\mathbb Z$-integral complex character of $\widetilde{H}_{\mathbb{C}}$, and its reduction modulo $p$ is supersingular if $\chi$ is not special or trivial.
  \end{enumerate}
\end{proposition}

\begin{proof}
The claims regarding integrality in (i) and (iii) are immediate.

  (i) This follows from Lemma~\ref{lem Bor}.

  (ii) and (iii): Let $\chi_0 : H \to \Z$ denote the underlying $\Z$-integral structure of $\chi$.
If we are not in one of the two exceptional cases, the result follows from Lemmas \ref{lem Bor}, \ref{lem ext} and \ref{lm:fact}\ref{lm:fact 3}.  Otherwise, the character $\chi_0$ of $H$ extends to a character $\chi'_0$ of $H' := \mathbb \Z[\Lambda_{\ker\nu}] \otimes H$ that is trivial on $\Lambda_{\ker\nu}$.  The tensor product $\chi'_0 \otimes_{H'} \widetilde{H}$ is a right $\widetilde{H}$-module that is free of rank $2$ (since the subgroup $\Lambda_{\ker\nu}$ of $\Omega$ has
  index $|\Psi| = 2$, by \eqref{eq fix 1}).
If $\chi' : H' \to \C$ denotes the base change of $\chi'_0$ to $\C$, then $\chi' \otimes_{H'_\C} \widetilde{H}_\C$ is simple and its
restriction to $H_{\mathbb{C}}$ is equal to $\chi \oplus \overline{\chi}$.  Note that the characters $\chi$ and $\overline{\chi}$ in (ii) are neither special nor trivial, since the $\chi_i$ are unequal by assumption and therefore have supersingular reduction modulo $p$. We conclude by Lemma~\ref{lm:fact}\ref{lm:fact 3}.
\end{proof}

\subsection{Discrete simple modules with supersingular reduction}
\label{sec:discr-simple-modul}

We continue to assume $\bf G$ is absolutely simple and isotropic.   Let $\p$ denote the maximal ideal of $\mathbb{Z}[q^{1/2}] \subset \C$ with residue field $\mathbb{F}_p$. We now discuss discrete, $\mathbb{Z}[q^{1/2}]$-integral $\widetilde{H}_{\mathbb{C}}$-modules with supersingular reduction modulo $\mathfrak{p}$.

The following is the key proposition of this section.

\begin{proposition}\label{prop key} 
Suppose the type of $\Sigma$ is not equal to $\textA_\ell$ with equal parameters.  Then there exists a right $\widetilde{H}_{\mathbb C} $-module $M_\mathbb C$ such that:
\begin{itemize}
\item $M_{\mathbb{C}}$ is simple and discrete as an $\widetilde{H}_{\mathbb C} $-module;
\item $M_{\mathbb{C}}$ has a $\mathbb Z[q^{1/2}]$-integral structure $M$ which is furthermore free over $\mathbb{Z}[q^{1/2}]$; 
\item $M$ has supersingular reduction modulo $\p$.  
\end{itemize}
\end{proposition}

The proposition will follow from Propositions \ref{prop 12}, \ref{prop discsinc}, and \ref{spin8} below.  We sketch the main ideas of the proof.  

Consider first the special character $\chi: H_{\mathbb{C}} \rightarrow \mathbb{C}$.  It is $\mathbb{Z}[q^{1/2}]$-integral, its reduction modulo $\mathfrak{p}$ is non-supersingular, and $\mathfrak{T}(\chi)$ is equal to the Steinberg representation of $G\sc$ over $\mathbb{C}$, so that $\chi$ is discrete.  Any discrete, non-special character of $H_{\mathbb{C}}$ is $\mathbb{Z}[q^{1/2}]$-integral (in fact, $\mathbb{Z}$-integral) and Lemma~\ref{lem Bor} implies that its reduction modulo $\mathfrak{p}$ is supersingular (since the trivial character of $H_{\mathbb{C}}$ is not discrete).  Thus, we first attempt to find a discrete non-special character of $H_{\mathbb{C}}$; these have been classified by Borel in \cite[\S 5.8]{Borel}.  (Note that in \cite{Borel}, the Iwahori subgroup is the pointwise stabilizer $\widetilde{Z}_0\mathfrak B$ of an alcove; recall again that if $\bG$ is $F$-split or semisimple and simply connected we have $\widetilde{Z}_0= Z_0$.)  We describe these characters in Proposition \ref{prop 12}, and use Proposition \ref{prop cha} to determine which of these characters extend to $\widetilde{H}_{\mathbb{C}}$.

When $m = 1$, there do not exist any discrete non-special characters of $H_{\mathbb{C}}$, and we use instead a reflection module of $\widetilde{H}_{\mathbb{Z}[q^{1/2}]}$ (see Proposition \ref{prop discsinc}).  It is free of rank $|S|$ over $\mathbb Z[q^{1/2}]$ and has supersingular reduction modulo $\p$.  When the type is $\textA_\ell$, this module is non-discrete, which is why we must omit this type.  (We also use reflection modules in Proposition \ref{spin8} to handle certain groups of type $\textB_3$ for which Proposition \ref{prop 12} does not apply.)

 We now proceed with the required propositions.

\begin{proposition} \label{prop 12} 
Suppose the type of $\Sigma$ is $\textB_\ell~ (\ell \geq 4)$, $\textC_\ell~ (\ell \geq 2)$, $\textF_4$, $\textG_2$, $\textA_1$ with parameters $d_1\neq d_2$, or $\textB_3$ with parameters $(d_1,d_2)\neq (1,2)$. 
Then the algebra $\widetilde{H}_{\mathbb{C}}$ admits a discrete non-special simple right module $M_\mathbb{C}$, induced from or extending a character of $H_{\mathbb{C}}$, which is $\mathbb{Z}[q^{1/2}]$-integral.
Moreover, the dimension of $M_{\mathbb{C}}$ is $1$, unless $\Psi \neq 1$ and the type is 
\begin{itemize}
\item $\textC_2$ with parameters $(1,1,1), (2,1,1)$, or $(3,2,2)$;
\item $\textC_3$ with parameters $(1,1,1), (1,2,2)$, or $(2,3,3)$;
\item $\textC_4$ with parameters $(1,2,2)$, or $(2,3,3)$;
\item $\textC_5$ with parameters $(1,2,2)$.
\end{itemize}
In these cases, $M_\mathbb{C}$ extends the $H_{\mathbb{C}}$-module $(-1,-1,q^d)\oplus (-1,q^d,-1)$ where $d := d_2 = d_3$, and thus the dimension of $M_{\mathbb{C}}$ is $2$.  
\end{proposition}

\begin{proof}
  When $m=1$, the only discrete character of $H_{\mathbb{C}}$ is the special one (\cite[\S 5.7]{Borel}).

Suppose $m > 1$.  For each choice of irreducible root system $\Sigma$, we list in Tables \ref{table1} and \ref{table2} the possible parameters $(d_1,d_2)$ or $(d_1,d_2,d_3)$ for $G$ (from the tables in \cite[\S 4]{Tits}), and describe if $H_{\mathbb{C}}$ has a discrete non-special character (using \cite[\S 5.8]{Borel}).  

We start with $m = 2$ in Table \ref{table1}.  For every entry marked ``Y,'' the given discrete non-special character extends to a character of $\widetilde{H}_{\mathbb{C}}$ using the condition of Lemma \ref{lem ext}.

\begin{table}[h]
\caption{\textbf{$m = 2$}}
\label{table1}
\begin{center}
\begin{tabular}{|c|c|c|}\hline
$\Sigma$ & Parameters & \begin{tabular}{c} $\exists$ discrete non-special \\ character of $H_{\mathbb{C}}$?\end{tabular} \\ \hline\hline
$\textA_1$ & $(d,d)~(d \geq 1)$ & N \\ \cline{2-3}
  & $(1,3)$ & Y \\ \cline{2-3}
  & $(2,3)$ & Y \\ \cline{2-3}
  & $(1,2)$ & Y \\ \cline{2-3}
  & $(1,4)$ & Y \\ \cline{2-3}
  & $(3,4)$ & Y \\ \hline\hline
 $\textB_\ell~(\ell \geq 3)$ & $(1,1)$ & Y \\ \cline{2-3}
 & $(1,2)$ & Y (if $\ell \geq 4$),\quad N (if $\ell = 3$)\\ \cline{2-3}
 & $(2,1)$ & Y \\ \cline{2-3}
 & $(2,3)$ & Y \\ \hline\hline
$\textF_4$ & $(1,1)$ & Y \\ \cline{2-3}
 & $(1,2)$ & Y \\ \cline{2-3}
 & $(2,1)$ & Y \\ \hline\hline
$\textG_2$ & $(1,1)$ & Y \\ \cline{2-3}
& $(1,3)$ & Y \\ \cline{2-3}
& $(3,1)$  & Y \\ \hline
\end{tabular}
\end{center}
\end{table}

We now consider $m = 3$ (that is, type $\textC_\ell$) in Table \ref{table2}.  In this case, the tables in \cite[\S 5.8]{Borel} show that $H_{\mathbb{C}}$ always admits a discrete, non-special character.  Note also that Borel omitted the parameters $(3,2,2)$ for type $\textC_2$.  In order to obtain this missing case, we use the criterion of \cite[Eqn.\ 5.6(2)]{Borel} to see that the only discrete non-special characters of $H_{\mathbb{C}}$ are $(-1,-1,1)$ and $(-1,1,-1)$ (in the notation of \cite{Borel}).  Note that the characters corresponding to parameters with $d_2 \neq d_3$ automatically extend to $\widetilde{H}_{\mathbb{C}}$, by Lemma \ref{lem ext}.  

(We have one more remark about the tables in \cite[\S 5.8]{Borel}: the character $(-1,-1,1)$ for parameters $(2,1,4)$ only works for $\ell \ge 3$.)

\begin{table}[h]
\caption{\textbf{$m = 3$}}
\label{table2}
\begin{center}
\begin{tabular}{|c|c|c|}\hline
$\Sigma$ & Parameters &  \begin{tabular}{c} Condition that some discrete \\ non-special character of $H_{\mathbb{C}}$ extends to $\widetilde{H}_{\mathbb{C}}$\end{tabular} \\ \hline\hline
$\textC_\ell~(\ell \geq 2)$ & $(1,1,1)$ & $\ell \geq 4$, or $\Psi = 1$ \\ \cline{2-3}
 & $(2,1,1)$ & $\ell \geq 3$, or $\Psi = 1$ \\ \cline{2-3}
 & $(2,3,3)$ & $\ell = 2, \ell \geq 5$, or $\Psi = 1$\\ \cline{2-3}
 & $(2,1,3)$ & none  \\ \cline{2-3}
 & $(1,1,2)$ & none  \\ \cline{2-3}
 & $(2,2,3)$ & none  \\ \cline{2-3}
 & $(2,1,2)$ & none  \\ \cline{2-3}
 & $(1,2,2)$ & $\ell = 2, \ell \geq 6$, or $\Psi = 1$ \\ \cline{2-3}
 & $(2,1,4)$ & none  \\ \cline{2-3}
 & $(2,3,4)$ & none \\ \hline
$\textC_2$ & $(3,2,2)$ & $\Psi = 1$ \\ \hline
\end{tabular}
\end{center}
\end{table}

Finally, we remark that in all cases, Propositions \ref{prop dis2} and \ref{prop cha} imply that the $\widetilde{H}_{\mathbb{C}}$-module $M_\mathbb{C}$ constructed above (either as the extension of a character of $H_{\mathbb{C}}$, or as the induction of a character from $H_{\mathbb{C}}$ to $\widetilde{H}_{\mathbb{C}}$) is discrete and $\mathbb{Z}[q^{1/2}]$-integral.  
\end{proof}

We consider now the types $\textD_{\ell}~ (\ell \geq 4)$, $\textE_6$, $\textE_7$, and $\textE_8$.  The tables in \cite[\S 4]{Tits} imply that $\bf G$ is $F$-split, so that $d_s=1$ for all $s\in S$, and for distinct $s,t\in S$, the order $n_{s,t}$ of $st$ is $2$ or $3$.

\begin{proposition}\label{prop discsinc} 
Assume that the type of $\Sigma$ is $\textD_{\ell}~ (\ell \geq 4)$, $\textE_6$, $\textE_7$, or $\textE_8$.  Let $M$ denote the right $\widetilde{H}_{\mathbb{Z}[q^{1/2}]}$-module obtained as the twist of the \(left\) reflection $\widetilde{H}_{\mathbb{Z}[q^{1/2}]}$-module by the anti-automorphism $T_w \mapsto (-1)^{\ell(w)}T_{w^{-1}}^*$.  Then $M$ is free of rank $|S|$ over $\mathbb{Z}[q^{1/2}]$, has supersingular reduction modulo $\p$, and $M_\mathbb C$ is a discrete simple right $\widetilde{H}_{\mathbb C}$-module.
\end{proposition}

\begin{proof}
  The left reflection $\widetilde{H}_{\mathbb Z[q^{1/2}]}$-module is the free $\mathbb{Z}[q^{1/2}]$-module with basis $\{e_t\}_{t\in S}$, with $\widetilde{H}_{\mathbb{Z}[q^{1/2}]}$-module structure given by
  \begin{align*}
 T_s\cdot e_t & = \begin{cases}  - e_t & \text{if } s=t,\\ {q}e_t & \text{if } s\neq t, \ n_{s,t}= 2,\\ { q}e_t+ q^{1/2}e_s& \text{if } s\neq t, \ n_{s,t}= 3, \end{cases} \\
 T_u\cdot e_t & = e_{u(t)},
  \end{align*}
  where $s,t\in S, u\in \Omega$.  Twisting this module by the automorphism $T_w\mapsto (-1)^{\ell (w)}T_w^*$ gives a left $\widetilde{H}_{\mathbb{Z}[q^{1/2}]}$-module $M'$, satisfying
  \begin{align*}
  T_s \cdot e_t & = \begin{cases} q e_t & \text{if } s=t,\\ - e_t & \text{if } s\neq t, \ n_{s,t}= 2,\\ - e_t- q^{1/2}e_s& \text{if } s\neq t, \ n_{s,t}= 3, \end{cases} \\
  T_u \cdot e_t & = e_{u(t)}.
  \end{align*}
  Finally, we define $M$ to be the right $\widetilde{H}_{\mathbb{Z}[q^{1/2}]}$-module obtained from $M'$ by applying the anti-automorphism $T_w \mapsto T_{w^{-1}}$.  The $\widetilde{H}_{\C}$-module $M_{\C}$ is simple (even as an $H_{\C}$-module, cf.\ \cite[\S 3.13]{Lusztig}).
  
  By applying Lemma~\ref{lm:fact}\ref{lm:fact 3} and Proposition~\ref{prop dis2} twice, we may assume that $\bG$ is adjoint in order to prove the required properties of $M$.  
The reduction modulo $\p$ of $M$ is the $\mathbb F_p$-vector space with basis $\{e_t\}_{t\in S}$, with the structure of a right $\widetilde{H}_{\mathbb{F}_p}$-module given by
 \begin{align*}
e_t \cdot T_s& = \begin{cases} 0& \text{if } s=t,\\ -e_t& \text{if } s\neq t, \end{cases} \\
e_t \cdot T_u & = e_{u^{-1}(t)}.
  \end{align*}
  The restriction to $H_{\mathbb{F}_p}$ of this $\widetilde{H}_{\mathbb{F}_p}$-module is the direct sum of the characters $\{\chi_s\}_{s\in S}$, where 
  $$\chi_s (T_t)=\begin{cases} 0 & \text{if } s=t,\\ -1& \text{if } s\neq t. \end{cases}$$ 
By Lemmas \ref{lm:fact}\ref{lm:fact 3} and \ref{lem Bor}, we deduce that $M_{\mathbb{F}_p}$ is supersingular.    Further, one checks that the right action of $\theta_{\lambda_\mu^{-1}}$ ($\mu\in X_*(T)^+$) on $M_{\mathbb{C}}$ is equal to the left action of $(-1)^{\ell(\lambda_\mu)}\hat{T}_{\lambda_\mu^{-1}}^{-1}$ on (the base change to $\mathbb{C}$ of) the reflection module, where $\hat{T}_{\lambda_\mu^{-1}}$ is defined in \cite[\S 4.3]{Lusztig} (note that with respect to our normalizations, the elements $\omega_i$ of \emph{op.\ cit.} are anti-dominant).  The discreteness of $M_{\mathbb{C}}$ now follows from Proposition \ref{prop dis} and \cite[\S 3.2, Thm.\ 4.7]{Lusztig}.  (See also \cite[\S4.23]{Lusztig}.)
\end{proof}

Finally, we consider one of the omitted cases from Proposition \ref{prop 12}, namely type $\textB_3$ with parameters $(1,2)$.

\begin{prop}\label{spin8}
Assume that the type of $\Sigma$ is $\textB_3$ with parameters $(1,2)$.
Then $\widetilde{H}_{\mathbb{Z}[q^{1/2}]}$ admits a right module $M$, such that $M$ is free of rank $3$ over $\mathbb{Z}[q^{1/2}]$, has supersingular reduction modulo $\p$, and $M_\mathbb C$ is a discrete simple right $\widetilde{H}_{\mathbb C}$-module.
\end{prop}

\begin{proof} 
In this case, the group $\mathbf{G}\sc$ is an unramified non-split form of $\mathbf{Spin}_8$, by the tables in \cite{Tits}.  
We will use the reflection module as defined in \cite[\S 7]{GS}.

Denote by $\widetilde{\Delta}_{\textnormal{long}}$ the subset of simple affine roots $\widetilde{\Delta}$ which are long.  We define an action of $H_{\mathbb{Z}[q^{1/2}]}$ on the free $\mathbb{Z}[q^{1/2}]$-module of rank~3 with basis $\{e_\beta\}_{\beta\in \widetilde{\Delta}_{\textnormal{long}}}$ as follows.  If $\alpha\in \widetilde{\Delta}_{\textnormal{long}}$, we set
$$T_{s_\alpha} \cdot e_\beta = \begin{cases} - e_\beta & \textnormal{if}~\alpha = \beta, \\ qe_\beta & \textnormal{if}~\alpha\neq \beta,~ n_{s_\alpha,s_\beta} = 2, \\ qe_\beta + q^{1/2}e_\alpha & \textnormal{if}~\alpha\neq \beta,~ n_{s_\alpha,s_\beta} = 3,\end{cases}$$
and if $\alpha$ is the unique short root in $\widetilde{\Delta}$, we set 
$$T_{s_\alpha} \cdot e_\beta = q^2e_\beta.$$
Twisting this reflection module by the automorphism $T_w \mapsto (-1)^{\ell(w)}T_w^*$ gives a new left $H_{\mathbb{Z}[q^{1/2}]}$-module $M'$, with action given by
$$T_{s_\alpha} \cdot e_\beta =  \begin{cases}  qe_\beta & \textnormal{if}~\alpha = \beta, \\ -e_\beta & \textnormal{if}~\alpha\neq \beta,~ n_{s_\alpha,s_\beta} = 2, \\ -e_\beta - q^{1/2}e_\alpha & \textnormal{if}~\alpha\neq \beta,~ n_{s_\alpha,s_\beta} = 3,\end{cases}$$
if $\alpha \in  \widetilde{\Delta}_{\textnormal{long}}$, and 
$$T_{s_\alpha} \cdot e_\beta = -e_\beta$$
if $\alpha \in  \widetilde{\Delta}$ is short.  
We extend the action of $H_{\mathbb{Z}[q^{1/2}]}$ on $M'$ to $\widetilde{H}_{\mathbb{Z}[q^{1/2}]}$ by declaring that
$$T_u \cdot e_{\alpha} = e_{u(\alpha)}.$$
As the algebra $\widetilde{H}_{\mathbb{Z}[q^{1/2}]}$ is generated by $H_{\mathbb{Z}[q^{1/2}]}$ and the elements $T_u$, $u \in \Omega$, subject to the relations
$T_{uv} = T_u T_v$ and $T_u T_{s_\alpha} T_u^{-1} = T_{s_{u(\alpha)}}$ for $u,v \in \Omega$ and $\alpha \in \wt\Delta$, we see that $M'$ is a well-defined module
of $\widetilde{H}_{\mathbb{Z}[q^{1/2}]}$.  
Finally, we define $M$ to be the right $\widetilde{H}_{\mathbb{Z}[q^{1/2}]}$-module obtained from $M'$ by applying the anti-automorphism $T_w \mapsto T_{w^{-1}}$.  One checks directly that $M_{\C}$ is simple (even as an $H_{\C}$-module).

By Lemma~\ref{lm:fact}\ref{lm:fact 3} and Proposition~\ref{prop dis2} we are now reduced to the case where $\bG$ is simply connected.  The reduction modulo $\p$ of $M$ is the $\mathbb{F}_p$-vector space with basis $\{e_{\beta}\}_{\beta\in \widetilde{\Delta}_{\textnormal{long}}}$, with the structure of a right $H_{\mathbb{F}_p}$-module given by
$$ e_\beta \cdot T_{s_\alpha} =  \begin{cases} 0  & \textnormal{if}~\alpha = \beta, \\ -e_\beta & \textnormal{if}~\alpha\neq \beta,\end{cases}$$
for $\alpha \in  \widetilde{\Delta}$.  Therefore $M_{\mathbb{F}_p}$ is equal to the direct sum of the characters $\{\chi_\beta\}_{\beta \in \widetilde{\Delta}_{\textnormal{long}}}$, where
$$\chi_{\beta}(T_{s_\alpha}) = \begin{cases} 0 & \textnormal{if}~\alpha = \beta,\\ -1 & \textnormal{if}~\alpha \neq \beta.\end{cases}$$
for $\alpha \in \widetilde{\Delta}$.  Lemma \ref{lem Bor} therefore implies that $M_{\mathbb{F}_p}$ is supersingular.

Once again, we see that the right action of $\theta_{\lambda_\mu^{-1}}$ ($\mu\in X_*(T)^+$) on $M_{\mathbb{C}}$ is equal to the left action of $(-1)^{\ell(\lambda_\mu)} q_{\lambda_\mu}^{1/2} T_{\lambda_\mu^{-1}}^{-1}$ on (the base change to $\mathbb{C}$ of) the reflection module.  Section 8.5 of \cite{GS} gives an explicit description of Hecke operators associated to the fundamental anti-dominant coweights in terms of $T_u$ and the $T_{s_\alpha}$. 
 Using this description along with Proposition \ref{prop dis}, we see that the $\widetilde{H}_{\mathbb{C}}$-module $M_{\mathbb{C}}$ is discrete.  (See also \cite[Prop.\ 7.11]{GS}.)
\end{proof}

\subsection{Admissible integral structure via discrete cocompact subgroups}
\label{S discrete}

Let $E$ be a number field with ring of integers $\O_E$, $\p$ a maximal ideal of $\O_E$ with residue field $k:=\O_E/\p$, and $C/E$
a field extension.

For any extension of fields, the scalar extension functor commutes with the $\mathfrak{B}$-invariant functor and its left adjoint $\mathfrak{T}$ (\cite[Lem.\ III.1(ii)]{HenV}).
Therefore, if $\tau$ is an $E$-structure of a $C$-representation $\pi$, then $\tau^{\mathfrak{B}}$ is an $E$-structure of $\pi^{\mathfrak{B}}$.  Conversely, if $M$ is an $E$-structure of an $\widetilde{H}_{C}$-module $N$, then $\mathfrak{T}(M)$ is an $E$-structure  of $\mathfrak{T}(N)$.

\begin{definition} \label{def O int}
  We say that an admissible $C$-representation $\pi$ of $G$ is \emph{$\O_E$-integral} if $\pi$ contains a $G$-stable $\O_E$-submodule $\tau^\circ$
  such that, for any compact open subgroup $K$ of $G$, the $\O_E$-module $(\tau^\circ)^K$ is finitely generated, and the natural map
  $$\varphi: C\otimes_{\O_E}\tau^\circ \to \pi$$
   is an isomorphism.  We call $\varphi$ (and more often $\tau^\circ$) an \emph{$\O_E$-integral structure} of $\pi$. The $G$-equivariant map $\tau^\circ \to k\otimes_{\O_E}\tau^\circ$ (and more often the $k$-representation $k\otimes_{\O_E} \tau^\circ$ of $G$) is called the \emph{reduction of $\tau^\circ$ modulo $\p$}.  We say that
  $\tau^\circ$ is \emph{admissible} if $k\otimes_{\O_E}\tau^\circ$ is admissible for all $\p$.
\end{definition}

For any commutative ring $R$ and any discrete cocompact subgroup $\Gamma$ of $G$, we define
$$C^\infty(\Gamma \backslash G,R) := \left\{ f :G \to R~   \big | \ f(\gamma g k)= f(g) \ \begin{array}{l} \textnormal{for all}~ \gamma\in \Gamma, g\in G,\\ \textnormal{and}~ k\in K_f\end{array} \right\},$$
where $K_f$ is some compact open subgroup of $G$ depending on $f$. Letting $G$ act on this space by right translation, we obtain a smooth $R$-representation $\rho^\Gamma_R$.  The complex representation $\rho_{\mathbb C}^\Gamma$ of $G$ has an admissible $\O_E$-integral structure given by $\rho^\Gamma_{\O_E}$: the reduction of $\rho^\Gamma_{\O_E}$ modulo $\mathfrak{p}$ is the admissible representation $\rho_{k}^\Gamma$.

\begin{proposition} \label{prop emb} 
Assume $\chr F =0$ and $\bf G$ semisimple. If $\pi$ is a square-integrable $\mathbb{C}$-representation of $G$, then there exists a discrete cocompact subgroup $\Gamma$ of $G$ such that 
$$\Hom_G(\pi, \rho^\Gamma_{\mathbb C})\neq 0.$$
\end{proposition}

\begin{proof}
  Since $\chr F =0$, there exists a decreasing sequence $(\Gamma_n)_{n\in \mathbb N}$ of discrete cocompact subgroups of $G$ with trivial intersection, such that each is
  normal and of finite index in $\Gamma_0$. (See \cite[Thm.\ A]{BorelHarder}. The construction there is global, and we obtain the required decreasing sequence
  by passing to congruence subgroups.)
   For any discrete cocompact subgroup $\Gamma$, the normalized multiplicity of
  $\pi$ in $\rho^\Gamma_{\mathbb C}$ is
  $$m_{\Gamma, dg}(\pi):= \vol_\Gamma\cdot \dim_{\mathbb C}\big(\Hom_ G(\pi, \rho^\Gamma_{\mathbb C})\big),$$ 
  where $\vol_\Gamma$ is the volume of $\Gamma \backslash G$ for a $G$-invariant measure induced by a fixed Haar measure on $G$.  By the square-integrability assumption on $\pi$ and the limit multiplicity formula, the sequence $(m_{\Gamma_n, dg}(\pi))_{n \in \mathbb{N}}$ converges to a nonzero real number (see \cite[App.\ 3, Prop.]{DKV} and \cite[Thm.\ K]{Kazhdan}).
\end{proof}

\begin{proposition} \label{prop mod} 
Assume $\chr F=0$. Let $\pi$ be an irreducible $\mathbb C$-representation of $G$ and $\Gamma$ a discrete cocompact subgroup of $G$.
  \begin{enumerate}
  \item If 
  $\varphi: \mathbb C\otimes_{E}\tau \congto \pi$
  is an $E$-structure of $\pi$, then 
  the natural map 
  $$\mathbb C \otimes_E\Hom_{E[G]}(\tau,\rho^\Gamma_{E}) \to \Hom_{ \mathbb{C}[G]}( \pi,\rho^\Gamma_{ \mathbb C})$$ 
  is an isomorphism. \label{prop mod 1}
  \item Any irreducible subrepresentation $\tau$ of $\rho^\Gamma_E$ admits an admissible $\O_E$-integral structure
    $\tau \cap \rho_{\O_E}^\Gamma$, whose reduction modulo $\p$ is contained in $\rho^\Gamma_k$. \label{prop mod 2}
  \end{enumerate}
\end{proposition}

\begin{proof}
  We recall a general result in algebra from \cite[\S 12.2 Lem.\ 1]{BkiA8}: let $C'/C$ be a field extension and $A$ a
  $C$-algebra. For $A$-modules $M,N$, the natural map
  \begin{equation}\label{eq:K}
  C' \otimes_{C}\Hom_{A} (M,N) \to \Hom_{C'\otimes_{C}A} (C'\otimes_{C}M, C'\otimes_{C}N) \end{equation} 
  is injective, and bijective if $C'/C$ is finite or the $A$-module $M$ is finitely generated.

  (i) Take $C'/C = \mathbb C/E$, $A = E[G]$, $(M,N)= (\tau, \rho_E^\Gamma)$. Then \eqref{eq:K} is an
  isomorphism because $\tau$ is an irreducible $E$-representation of $G$.

  (ii) For any compact open subgroup $K$ of $G$, the $\O_E$-module $( \rho_{\O_E}^\Gamma)^K$ is finite free and $\rho_{\O_E}^\Gamma$ contains $\tau^\circ := \tau \cap \rho_{\O_E}^\Gamma$ as $\O_E$-representations of $G$. Since the ring $\O_E$ is noetherian, these facts imply the $\O_E$-submodule $(\tau^\circ)^K$ of $(\rho_{\O_E}^\Gamma)^K$ is finitely generated.  The natural linear $G$-equivariant isomorphism 
  $$E\otimes_{\O_E}\rho^\Gamma_{\O_E} \congto \rho^\Gamma_{E}$$ 
  restricts to a linear $G$-equivariant isomorphism 
  $$E\otimes_{\O_E} \tau^\circ \congto \tau,$$ 
  and therefore $\tau^\circ$ is an $\O_E$-integral structure of $\tau$. It remains to verify that the injection $\tau^\circ \into \rho_{\O_E}^\Gamma$ stays
  injective after reduction modulo $\p$. (As $\rho^\Gamma_{k}$ is admissible, this will also imply that $k\otimes_{\O_E} \tau^\circ$ is admissible.)
  More generally, suppose that $0 \to M' \to M \to M'' \to 0$ is any exact sequence of $\O_E$-modules with $M''$ torsion-free.
  Then $M''$ is the direct limit of its finitely generated submodules, and finitely generated torsion-free modules are projective, as $\O_E$ is Dedekind. Hence
  $\Tor_1^{\O_E}(M'',k) = 0$, as Tor functors commute with direct limits, so the sequence stays exact after reduction modulo~$\p$.
\end{proof}

The above result will be used in our construction of irreducible, admissible, supersingular $C$-representations.  It also has the following consequence, which may be of independent interest.  

\begin{corollary} \label{cor sca} Assume $\chr F=0$ and $\bf G$ semisimple. Then any irreducible supercuspidal $\mathbb{C}$-representation admits an admissible $\O_E$-integral structure whose reduction modulo $\p$ is contained in $\rho_k^\Gamma$, for some discrete cocompact subgroup $\Gamma$ of $G$.
\end{corollary}

\begin{proof}
When $\bf G$ is semisimple, any irreducible admissible supercuspidal $\mathbb C$-representation $\pi$ of $G$ descends to a number field (see \cite[II.4.9]{Viglivre}).  Since $\pi$ is in particular square-integrable, Proposition \ref{prop emb} implies that $\pi$ embeds into $\rho_{\mathbb{C}}^\Gamma$ for some discrete cocompact subgroup $\Gamma$ of $G$.  The claim the follows from points \ref{prop mod 1} and \ref{prop mod 2} of Proposition \ref{prop mod}.
\end{proof}

\subsection{Reduction to rank 1 and \texorpdfstring{$\PGL_n(D)$}{PGL\_n(D)}}
\label{sec:reduction-rank-1-pgln}

We now prove that most $p$-adic reductive groups admit irreducible admissible supersingular (equivalently, supercuspidal) representations.

\begin{theorem}\label{thm main} 
  Assume that $c = p$ and $\chr F = 0$.  Suppose $\bG$ is an isotropic, absolutely simple, connected adjoint $F$-group, not isomorphic to any of the following groups:
  \begin{enumerate}
  \item $\bPGL_n(D)$, where $n \ge 2$ and $D$ a central division algebra over $F$;
  \item $\bPU(h)$, where $h$ is a split hermitian form in $3$ variables over a ramified quadratic extension of $F$
    or a non-split hermitian form in $4$ variables over the unramified quadratic extension of $F$.
  \end{enumerate}
  Then $G$ admits an irreducible admissible supercuspidal $C$-representation.
\end{theorem}

\begin{proof}
  We first note by the tables in \cite{Tits} that the above exceptional groups are precisely the ones where $\Sigma$ is of type
  $\textA_\ell$ with equal parameters.  (In that reference our exceptional groups have names $\textA_{m-1}$,
  ${}^d\textA_{md - 1}$ for $m \ge 2$, $d \ge 2$ in case (i) and $\text{C-BC}_1$, ${}^2\textA_3''$ in case (ii).)

By Proposition \ref{prop key} there exists a right $\widetilde{H}_{\Z[q^{1/2}]}$-module $M$ which is free over $\mathbb{Z}[q^{1/2}]$, whose base change $M_{\C}$
  is a discrete simple $\widetilde{H}_{\mathbb{C}}$-module, and whose reduction $M_{\fp}$ is supersingular.  Set $E := \mathbb{Q}(q^{1/2})$, so that $\mathbb{Z}[q^{1/2}] \subset \mathcal{O}_E$.  Let $\pi := \mathfrak{T}(M_{\mathbb{C}})$ denote the irreducible square-integrable $\C$-representation of $G$ corresponding to $M_{\C}$; then $\tau := \mathfrak{T}(M_E)$ is a $E$-structure of $\pi$.  We know by Proposition~\ref{prop emb} that $\pi$ injects into $\rho^\Gamma_{\mathbb{C}}$ for some discrete cocompact subgroup $\Gamma$ of $G$, and therefore $\tau$ injects into $\rho^\Gamma_E$ by Proposition~\ref{prop mod}\ref{prop mod 1}.  We identify $\pi$ and $\tau$ with their images in $\rho^\Gamma_{\mathbb{C}}$ and $\rho^\Gamma_E$, respectively.  Proposition~\ref{prop mod}\ref{prop mod 2} then ensures that $\tau^\circ := \tau \cap \rho^\Gamma_{\mathcal{O}_E}$ is an admissible $\mathcal{O}_E$-integral structure of $\tau$.  In particular, we have a $G$-equivariant map $\mathbb{C} \otimes_{\mathcal{O}_E} \tau^\circ  \congto \pi$.

Define $M' := (\tau^\circ)^{\mathfrak{B}}$; since $E$ is a localization of $\mathcal{O}_E$, the isomorphism above implies $\mathbb{C}  \otimes_{\mathcal{O}_E} M' \congto \pi^{\mathfrak{B}} \cong M_{\mathbb{C}}$, so that $M'$ is an $\mathcal{O}_E$-integral structure of $M_\C$.  Let $\mathfrak{p} \subset \mathcal{O}_E$ denote the prime ideal lying over $p$, and let $\mathcal{O}_{\mathfrak{p}} \subset \mathbb{C}$ denote the localization of $\mathcal{O}_E$ at $\mathfrak{p}$.  Then $M'_{\mathcal{O}_{\mathfrak{p}}} = \mathcal{O}_{\mathfrak{p}} \otimes_{\mathcal{O}_E} M'$ is a finitely generated, torsion-free module over the discrete valuation ring $\mathcal{O}_{\mathfrak{p}}$, which implies it is free.

Both $M_{\mathcal{O}_{\mathfrak{p}}}$ and $M'_{\mathcal{O}_{\mathfrak{p}}}$ are $\widetilde{H}_{\mathcal{O}_{\mathfrak{p}}}$-modules which are free over $\mathcal{O}_{\mathfrak{p}}$, and they are isomorphic over $\mathbb{C}$.  Thus, we see that the reductions $M_{\mathcal{O}_{\mathfrak{p}}/\mathfrak{p}}$ and $M'_{\mathcal{O}_{\mathfrak{p}}/\mathfrak{p}}$ agree up to semisimplification by the Brauer--Nesbitt theorem.  In particular, $M'_{\mathcal{O}_{\mathfrak{p}}/\mathfrak{p}} = M'_{\mathcal{O}_{E}/\mathfrak{p}}$ is supersingular (since the same is true of $M_{\mathcal{O}_{\mathfrak{p}}/\mathfrak{p}} = M_{\mathbb{F}_p}$) and, by construction, $M'_{\mathcal{O}_{E}/\mathfrak{p}}$ is a submodule of $(\rho^{\Gamma}_{\mathcal{O}_{E}/\mathfrak{p}})^{\mathfrak{B}}$ (this uses the final claim of Proposition~\ref{prop mod} \ref{prop mod 2}). Therefore we can pick a non-zero supersingular element  $v$ of $(\rho^{\Gamma}_{\mathcal{O}_{E}/\mathfrak{p}})^{\mathfrak{B}}$. The $G$-representation $\rho^\Gamma_{\mathcal{O}_{E}/\mathfrak{p}}$ is admissible, as $\Gamma$ is cocompact, and hence so is its subrepresentation $\ang{G\cdot v}$ generated by $v$.  Any
  irreducible quotient of $\ang{G\cdot v}$ (which exists by Zorn's lemma) is admissible by \cite[\S4, Thm.\ 1]{henniart:adm}, as $F$ is of characteristic
  zero, and supersingular by Proposition~\ref{prop ss=sc}, as it contains (the nonzero image of) $v$.  The theorem now follows from Proposition~\ref{prop Cfinite}.    
  
 \end{proof}

The two exceptional cases will be dealt with in Sections~\ref{sec:supers-repr-rank}, \ref{sec:supers-repr-pgl_n} below.
Assuming this, we can now prove our main result.

\begin{proof}[Proof of Theorem~\ref{thm:main}]
  Suppose that $\bG$ is a connected reductive group over $F$. We want to show that $G = \bG(F)$ admits an irreducible
  admissible supercuspidal representation over any field $C$ of characteristic $p$. By Proposition~\ref{prop Cfinite} we may
  assume that $C$ is finite and as large as we like. Then by Proposition~\ref{prop redadjoint} we may assume that $\bG$ is
  isotropic, absolutely simple, and connected adjoint. The result then follows from Theorem~\ref{thm main}, Corollary \ref{cor rk1}, and
  Corollary~\ref{cor:main-PGLn}.
\end{proof}

\section{Supersingular representations of rank 1 groups}
\label{sec:supers-repr-rank}

In this section we verify Theorem~\ref{thm:main} when $\bG$ is a connected reductive $F$-group of relative semisimple rank 1.
In particular, this deals with the second exceptional case in Theorem~\ref{thm main}.

\subsection{Preliminaries}
\label{sec:prelim-rk1}
We suppose in this section that $C$ is a finite extension of $\mathbb{F}_p$ which contains the $|G|_{p'}$-th roots of unity, where $|G|_{p'}$ denotes the prime-to-$p$ part of the pro-order of $G$.

Suppose that $\chr F = 0$.  We will show that $G$ admits irreducible, admissible, supercuspidal $C$-representations.  By Proposition \ref{prop redadjoint}, it suffices to assume $\mathbf{G}$ is an absolutely simple and adjoint group of relative rank 1.  We make one further reduction.  Let $\mathbf{G}\sc$ denote the simply-connected cover of $\mathbf{G}$:
$$1 \rightarrow \mathbf{Z}(\mathbf{G}\sc) \rightarrow \mathbf{G}\sc \rightarrow \mathbf{G} \rightarrow 1.$$
By Proposition \ref{prop HtoG}, we see that $G\sc$ admits an irreducible, admissible, supercuspidal representation on which $Z(G\sc)$ acts trivially if and only if $G$ does.  Therefore, we may assume that our group $\mathbf{G}$ is absolutely simple, simply connected, and has relative rank equal to 1.  We will then construct irreducible, admissible, supercuspidal representations of $G$ on which its (finite) center acts trivially.

\subsection{Parahoric subgroups}
\label{sec:parahoric-subgroups}
Let $\mathscr{B}$ denote the adjoint Bruhat--Tits building of $G$.  By our assumptions on $\mathbf{G}$, $\mathscr{B}$ is a one-dimensional contractible simplicial complex, i.e., a tree.  Recall that $\mathcal{C}$ denotes the chamber of $\mathscr{B}$ corresponding to the Iwahori subgroup $\mathfrak{B}$, and let $x_0$ and $x_1$ denote the two vertices in the closure of $\mathcal{C}$.  We let $K_0$ and $K_1$ denote the pointwise stabilizers of $x_0$ and $x_1$, respectively.  We then have $\mathfrak{B} = K_0 \cap K_1$.  

The vertices $x_0$ and $x_1$ are representatives of the two orbits of $G$ on the set of vertices of $\mathscr{B}$, and the edge $\mathcal{C}$ is a representative of the unique orbit of $G$ on the edges of $\mathscr{B}$.  By \cite[\S 4, Thm.\ 6]{serre:trees}, we may therefore write the group $G$ as an amalgamated product:
$$G \cong K_0 *_{\mathfrak{B}} K_1.$$

Since the group $\mathbf{G}$ is semisimple and simply connected, the stabilizers of vertices and edges in $\mathscr{B}$ are parahoric subgroups (see, e.g., \cite[\S 3.7]{Vig1}).  For $i \in \{0,1\}$, we let $K_i^+$ denote the pro-$p$ radical of $K_i$, that is, the largest open, normal, pro-$p$ subgroup of $K_i$.  The quotient $\mathds{G}_i := K_i/K_i^+$ is isomorphic the group of $k_F$-points of a connected reductive group over $k_F$ (see \cite[\S 3.7]{HV1}).  Likewise, the pro-$p$-Sylow $\mathfrak{U}$ is the largest open, normal, pro-$p$ subgroup of $\mathfrak{B}$, and $\mathds{Z} := \mathfrak{B}/\mathfrak{U}$ is isomorphic to the group of $k_F$-points of a torus over $k_F$.  The image of $\mathfrak{B}$ in $\mathds{G}_i$ is equal to a minimal parabolic subgroup $\mathds{B}_i$, with Levi decomposition $\mathds{B}_i = \mathds{Z}_i\mathds{U}_i$.  Thus, we identify the quotient $\mathds{Z}$ with $\mathds{Z}_i$.

\subsection{Pro-\texorpdfstring{$p$}{p} Iwahori--Hecke algebras}
\label{sec:pro-p-iwahori}
We work in slightly greater generality than in \S\ref{sec:IHalgs}.  Let 
$$\mathcal{H}_C := H_C(G,\mathfrak{U}) = \End_{G}C[\mathfrak{U}\backslash G]$$ 
denote the pro-$p$ Iwahori--Hecke algebra of $G$ with respect to $\mathfrak{U}$.  We view $\mathcal{H}_C$ as the convolution algebra of $C$-valued, compactly supported, $\mathfrak{U}$-bi-invariant functions on $G$ (see \cite[\S 4]{Vig1} for more details).  For $g\in G$, we let $T_g$ denote the characteristic function of $\mathfrak{U}g\mathfrak{U}$.  The algebra $\mathcal{H}_C$ is generated by two operators $T_{\tilde{s}_0}, T_{\tilde{s}_1}$, where $\tilde{s}_0$ and $\tilde{s}_1$ are lifts to the pro-$p$ Iwahori--Weyl group $\mathcal{N}/(Z\cap \mathfrak{U})$ of affine reflections $s_0, s_1$ fixing $x_0, x_1$, respectively, along with operators $T_{z}$ for $z\in \mathds{Z}$.  (Note that this labeling is different than the labeling in \S \ref{sec:characters}.)  For $i \in \{0, 1\}$, we let $\mathcal{H}_{C,i}$ denote the subalgebra of $\mathcal{H}_C$ generated by $T_{\tilde{s}_i}$ and $T_z$ for $z\in \mathds{Z}$; this is exactly the subalgebra of functions in $\mathcal{H}_C$ with support in $K_i$, i.e., 
$$\mathcal{H}_{C,i} = H_C(K_i,\mathfrak{U}) = \End_{K_i}C[\mathfrak{U}\backslash K_i].$$  
The algebra $\mathcal{H}_{C,i}$ is canonically isomorphic to the finite Hecke algebra $H_C(\mathds{G}_i, \mathds{U}_i)$ (see \cite[\S 6.1]{CabEng}).

Since $K_i^+$ is an open normal pro-$p$ subgroup of $K_i$, the irreducible smooth $C$-representations of $K_i$ and $\mathds{G}_i$ are in bijection.  Further, the finite group $\mathds{G}_i$ possesses a strongly split BN pair of characteristic $p$ (\cite[Prop.\ 3.25]{Vig1}).  Therefore, by \cite[Thm.\ 6.12]{CabEng}, the functor $\rho \mapsto \rho^{\mathfrak{U}}$ induces a bijection between isomorphism classes of irreducible smooth $C$-representations of $K_i$ and isomorphism classes of simple right $\mathcal{H}_{C,i}$-modules, all of which are one-dimensional.

We briefly recall some facts about supersingular $\mathcal{H}_C$-modules (compare Lemma \ref{lem Bor}).  We refer to \cite[Def.\ 6.10]{Vigss} for the precise definition (which is analogous to Definition \ref{def ssing}) and give instead the classification of simple supersingular $\mathcal{H}_C$-modules.  Since $\mathbf{G}$ is simply connected, every supersingular $\mathcal{H}_C$-module is a character.  The characters $\Xi$ of $\mathcal{H}_C$ are parametrized by pairs $(\chi,J)$, where $\chi: \mathds{Z} \rightarrow C^\times$ is a character of the finite torus and $J$ is a subset of 
$$S_\chi := \big\{s\in \{s_0, s_1\}: \chi(c_{\tilde{s}}) \neq 0\big\}$$
(here $c_{\tilde{s}}$ is a certain element of $C[\mathds{Z}]$ which appears in the quadratic relation for $T_{\tilde{s}}$; note also that the definition of $S_\chi$ is independent of the choice of lift $\tilde{s}\in \mathcal{N}/(Z\cap \mathfrak{U})$ of $s$).  The correspondence is given as follows (cf.\ \cite[Thm.\ 1.6]{Vigss}): for $z\in \mathds{Z}$, we have $\Xi(T_z) = \chi(z)$, and for $s\in \{s_0,s_1\}$, we have
$$\Xi(T_{\tilde{s}}) = \begin{cases}0 & \textnormal{if}~s\in J,\\ \chi(c_{\tilde{s}}) & \textnormal{if}~s\not\in J.\end{cases}$$
Since $\mathbf{G}$ is simple, \cite[Thm.\ 1.6]{Vigss} implies that $\Xi$ is supersingular if and only if 
$$(S_\chi,J)\neq (\{s_0,s_1\},\emptyset),~~ (\{s_0,s_1\},~\{s_0,s_1\}).$$

\subsection{Diagrams}
\label{sec:diagrams}  

Since the group $G$ is an amalgamated product of two parahoric subgroups, the formalism of diagrams used in \cite{KoziolXu} applies to the group $G$.  We recall that a diagram $D$ is a quintuple $(\rho_0, \rho_1, \sigma, \iota_0, \iota_1)$ which consists of smooth $C$-representations $\rho_i$ of $K_i$ ($i \in \{0,1\}$), a smooth $C$-representation $\sigma$ of $\mathfrak{B}$, and $\mathfrak{B}$-equivariant morphisms $\iota_i:\sigma \rightarrow \rho_i|_{\mathfrak{B}}$. We depict diagrams as
\begin{center}
\begin{tikzcd}
 & \rho_0\\
 \sigma \ar[ur,"\iota_0"] \ar[dr,"\iota_1"'] & \\
 & \rho_1
\end{tikzcd}
\end{center}
Morphisms of diagrams are defined in the obvious way (i.e., so that the relevant squares commute).

Let $\Xi$ denote a supersingular character of $\mathcal{H}_C$, associated to a pair $(\chi,J)$.  We define a diagram $D_\Xi$ as follows:
\begin{itemize}
\item set $\sigma := \chi^{-1}$, which we view as a character of $\mathfrak{B}$ by inflation;
\item we let $\rho_{\Xi,i}$ denote an irreducible smooth $C$-representation of $K_i$ such that $\rho_{\Xi,i}^{\mathfrak{U}} \cong \Xi|_{\mathcal{H}_{C,i}}$ as $\mathcal{H}_{C,i}$-modules (by the discussion above, $\rho_{\Xi,i}$ is unique up to isomorphism);
\item let $\iota_i$ denote the $\mathfrak{B}$-equivariant map given by $\sigma = \chi^{-1} \congto \rho_{\Xi,i}^{\mathfrak{U}} \hookrightarrow \rho_{\Xi,i}|_{\mathfrak{B}}$.
\end{itemize}
Pictorially, we write
\begin{center}
$D_\Xi \ \  = \quad$ 
$\left( \begin{tikzcd}
 & \rho_{\Xi,0}\\
   \chi^{-1} \ar[ur,"\iota_0"] \ar[dr,"\iota_1"'] & \\
  & \rho_{\Xi,1}
\end{tikzcd} \right)$
\end{center}

We now wish to construct an auxiliary diagram $D'$ into which $D_\Xi$ injects.  This will be done with the use of injective envelopes.  Recall that if $\mathcal{G}$ is a profinite group and $\tau$ is a smooth $C$-representation of $\mathcal{G}$, an injective envelope consists of a smooth injective $C$-representation $\textnormal{inj}_{\mathcal{G}}\tau$ of $\mathcal{G}$ along with a $\mathcal{G}$-equivariant injection $j: \tau \hookrightarrow \textnormal{inj}_\mathcal{G}\tau$ which satisfies the following property: for any nonzero $C$-subrepresentation $\tau' \subset \textnormal{inj}_\mathcal{G}\tau$, we have $j(\tau) \cap \tau' \neq 0$.  This data exists and is unique up to (non-unique) isomorphism.  

\begin{lem}[{\cite[Lem.\ 6.13]{paskunas:diags}}]\label{pask1}
Let $\tau$ denote a smooth $C$-representation of $\mathcal{G}$, and let $j: \tau \hookrightarrow \textnormal{inj}_{\mathcal{G}}\tau$ denote an injective envelope.  Let $\mathfrak{I}$ denote an injective representation of $\mathcal{G}$, and suppose we have an injection $\phi: \tau \hookrightarrow \mathfrak{I}$.  Then $\phi$ extends to an injection $\widetilde{\phi}: \textnormal{inj}_{\mathcal{G}}\tau \hookrightarrow \mathfrak{I}$ such that $\phi = \widetilde{\phi}\circ j$.  
\end{lem}

\begin{lem}\label{pask2}
Suppose $\mathcal{G}$ has an open, normal subgroup $\mathcal{G}^+$.  Let $\tau$ denote a smooth $C$-representation of $\mathcal{G}$ such that $\mathcal{G}^+$ acts trivially, and let $j: \tau \hookrightarrow \textnormal{inj}_{\mathcal{G}}\tau$ denote an injective envelope of $\tau$ in the category of $C$-representations of $\mathcal{G}$.  Then $\tau \hookrightarrow (\textnormal{inj}_{\mathcal{G}}\tau)^{\mathcal{G}^+}$ is an injective envelope of $\tau$ in the category of $C$-representations of $\mathcal{G}/\mathcal{G}^+$.
\end{lem}

\begin{proof}
This is \cite[Lem.\ 6.14]{paskunas:diags}; its proof does not require that $\tau$ be irreducible or that $\mathcal{G}^+$ be pro-$p$, as we assume that $\mathcal{G}^+$ acts trivially.
\end{proof}

We now begin constructing $D'$.

\begin{lem}
\label{injres}
Let $i \in \{0,1\}$.  We then have
$$(\textnormal{inj}_{K_i}C[\mathds{G}_i])|_{\mathfrak{B}} \cong \bigoplus_\xi \textnormal{inj}_{\mathfrak{B}}\xi^{\oplus |\mathds{B}_i \backslash \mathds{G}_i|},$$
where $\xi$ runs over all $C$-characters of $\mathfrak{B}$ \(or, equivalently, of $\mathds{Z}_i$\), and we have fixed choices of injective envelopes.  
\end{lem}

\begin{proof}
Consider the $\mathfrak{B}$-representation $(\textnormal{inj}_{K_i}C[\mathds{G}_i])^{\mathfrak{U}}$.  The action of $\mathfrak{B}$ factors through the quotient $\mathfrak{B}/\mathfrak{U} \cong \mathds{Z}$, which is commutative of order coprime to $p$.  Therefore, we obtain a $\mathfrak{B}$-equivariant isomorphism
\begin{equation}\label{inj invts}
(\textnormal{inj}_{K_i}C[\mathds{G}_i])^{\mathfrak{U}} \cong \bigoplus_\xi \xi^{\oplus m_\xi}
\end{equation}
for non-negative integers $m_\xi$ satisfying
\begin{align*}
m_\xi & = \dim_C \Hom_{\mathfrak{B}}(\xi, \textnormal{inj}_{K_i}C[\mathds{G}_i])\\
 & = \dim_C \Hom_{\mathfrak{B}}\big(\xi, (\textnormal{inj}_{K_i}C[\mathds{G}_i])^{K_i^+}\big) \\
 & = \dim_C \Hom_{\mathds{B}_i}(\xi, \textnormal{inj}_{\mathds{G}_i}C[\mathds{G}_i])\\
 & = \dim_C \Hom_{\mathds{Z}_i}\big(\xi, (\textnormal{inj}_{\mathds{G}_i}C[\mathds{G}_i])^{\mathds{U}_i}\big).
\end{align*}
(The third equality follows from Lemma \ref{pask2}.) Since $C[\mathds{G}_i]$ is injective as a representation of $\mathds{G}_i$, we have isomorphisms of $\mathds{Z}_i$-representations
$$(\textnormal{inj}_{\mathds{G}_i}C[\mathds{G}_i])^{\mathds{U}_i} \cong C[\mathds{U}_i\backslash \mathds{G}_i] \cong \bigoplus_{\xi} \xi^{\oplus |\mathds{B}_i\backslash \mathds{G}_i|},$$
so that $m_\xi = |\mathds{B}_i\backslash \mathds{G}_i|$.

The isomorphism \eqref{inj invts} implies we have a $\mathfrak{B}$-equivariant injection
$$\bigoplus_{\xi} \xi^{\oplus |\mathds{B}_i\backslash \mathds{G}_i|} \hookrightarrow (\textnormal{inj}_{K_i}C[\mathds{G}_i])|_{\mathfrak{B}}.$$
As $\mathfrak{B}$ is open, \cite[\S I.5.9 d)]{Viglivre} implies that the representation on the right-hand side is injective.  Lemma \ref{pask1} then says that the above morphism extends to a split injection between injective $\mathfrak{B}$-representations
$$\bigoplus_{\xi} \textnormal{inj}_{\mathfrak{B}}\xi^{\oplus |\mathds{B}_i\backslash \mathds{G}_i|} \hookrightarrow (\textnormal{inj}_{K_i}C[\mathds{G}_i])|_{\mathfrak{B}}.$$
Since the $\mathfrak{U}$-invariants of both representations agree, the above injection must be an isomorphism.  
\end{proof}

\begin{lem}
Set $a := \textnormal{lcm}\(|\mathds{B}_0\backslash \mathds{G}_0|,~|\mathds{B}_1\backslash \mathds{G}_1|\)$.  There exists a diagram $D'$ of the form
\begin{center}
$D' \ \ = \quad$
$\left( \begin{tikzcd}
  & \textnormal{inj}_{K_0}C[\mathds{G}_0]^{\oplus a\cdot |\mathds{B}_0\backslash \mathds{G}_0|^{-1}}\\
  \bigoplus_\xi \textnormal{inj}_{\mathfrak{B}}\xi^{\oplus a} \ar[ur,"\kappa_0"] \ar[dr,"\kappa_1"'] & \\
  & \textnormal{inj}_{K_1}C[\mathds{G}_1]^{\oplus a\cdot |\mathds{B}_1\backslash \mathds{G}_1|^{-1}}
\end{tikzcd} \right)$
\end{center}
where $\kappa_0$ and $\kappa_1$ are isomorphisms, and a morphism of diagrams
\begin{center}
\begin{tikzcd}
 & & \rho_{\Xi,0} \ar[rr, hook, "\psi_{K_0}"]& & \textnormal{inj}_{K_0}C[\mathds{G}_0]^{\oplus a\cdot |\mathds{B}_0\backslash \mathds{G}_0|^{-1}}\\
 \psi: & \chi^{-1} \ar[rr, hook, "\psi_{\mathfrak{B}}"] \ar[ur, "\iota_0"] \ar[dr, "\iota_1"'] & & \bigoplus_\xi \textnormal{inj}_{\mathfrak{B}}\xi^{\oplus a} \ar[ur,"\kappa_0"] \ar[dr,"\kappa_1"'] & \\
& & \rho_{\Xi,1} \ar[rr, hook, "\psi_{K_1}"] &  & \textnormal{inj}_{K_1}C[\mathds{G}_1]^{\oplus a\cdot |\mathds{B}_1\backslash \mathds{G}_1|^{-1}}
\end{tikzcd}
\end{center}
in which all arrows are injections.  
\end{lem}

\begin{proof}
We fix the following injections, which are equivariant for the relevant groups: 
\begin{itemize}
\item injective envelopes $j_\xi: \xi \hookrightarrow \textnormal{inj}_{\mathfrak{B}}\xi$ for each $C$-character $\xi$ of $\mathfrak{B}$;
\item injective envelopes $j_i: C[\mathds{G}_i]^{\oplus a\cdot |\mathds{B}_i\backslash \mathds{G}_i|^{-1}} \hookrightarrow \textnormal{inj}_{K_i}C[\mathds{G}_i]^{\oplus a\cdot |\mathds{B}_i\backslash \mathds{G}_i|^{-1}}$ for $i \in \{0,1\}$;
\item an inclusion $c:\chi^{-1} \hookrightarrow \bigoplus_\xi \xi^{\oplus a}$;
\item an inclusion $c_i:\rho_{\Xi,i} \hookrightarrow C[\mathds{G}_i]^{\oplus a\cdot |\mathds{B}_i \backslash \mathds{G}_i|^{-1}}$ for $i\in \{0,1\}$.
\end{itemize}

Let $i\in \{0,1\}$.  We first construct the $\kappa_i$.  We have a $\mathfrak{B}$-equivariant sequence of maps
$$\chi^{-1} \xhookrightarrow{\iota_i} \rho_{\Xi,i} \xhookrightarrow{c_i} C[\mathds{G}_i]^{\oplus a\cdot |\mathds{B}_i\backslash \mathds{G}_i|^{-1}} \xhookrightarrow{j_i} \textnormal{inj}_{K_i}C[\mathds{G}_i]^{\oplus a\cdot |\mathds{B}_i \backslash \mathds{G}_i|^{-1}}$$
and thus we obtain
$$\chi^{-1} \xhookrightarrow{j_i\circ c_i\circ\iota_i} \big(\textnormal{inj}_{K_i}C[\mathds{G}_i]^{\oplus a\cdot |\mathds{B}_i \backslash \mathds{G}_i|^{-1}}\big)^{\mathfrak{U}}.$$
By Lemmas \ref{injres} and \ref{pask2}, we have $\bigoplus_\xi \xi^{\oplus a} \cong (\textnormal{inj}_{K_i}C[\mathds{G}_i]^{\oplus a\cdot |\mathds{B}_i\backslash \mathds{G}_i|^{-1}})^{\mathfrak{U}}$.  We \emph{fix} an isomorphism $\alpha_i: \bigoplus_\xi \xi^{\oplus a} \congto (\textnormal{inj}_{K_i}C[\mathds{G}_i]^{\oplus a\cdot |\mathds{B}_i\backslash \mathds{G}_i|^{-1}})^{\mathfrak{U}}$ such that 
\begin{equation}
\label{eqn1}
\alpha_i \circ c = j_i \circ c_i \circ \iota_i.
\end{equation}

Now consider the maps of $C$-representations of $\mathfrak{B}$:
$$\bigoplus_\xi \xi^{\oplus a} \xhookrightarrow{\alpha_i} \big(\textnormal{inj}_{K_i}C[\mathds{G}_i]^{\oplus a\cdot |\mathds{B}_i\backslash \mathds{G}_i|^{-1}}\big)^{\mathfrak{U}} \hookrightarrow \big(\textnormal{inj}_{K_i}C[\mathds{G}_i]^{\oplus a\cdot |\mathds{B}_i\backslash \mathds{G}_i|^{-1}}\big)|_{\mathfrak{B}}.$$
By Lemma \ref{pask1}, the above map extends to an $\mathfrak{B}$-equivariant split injection
$$\kappa_i: \bigoplus_\xi \textnormal{inj}_{\mathfrak{B}}\xi^{\oplus a} \hookrightarrow \big(\textnormal{inj}_{K_i}C[\mathds{G}_i]^{\oplus a\cdot |\mathds{B}_i\backslash \mathds{G}_i|^{-1}}\big)|_{\mathfrak{B}}$$
such that 
\begin{equation}
\label{eqn2}
\kappa_i \circ \big(\bigoplus_\xi j_{\xi}^{\oplus a}\big) = \alpha_i.
\end{equation}
Since both $\bigoplus_\xi \textnormal{inj}_{\mathfrak{B}}\xi^{\oplus a}$ and $(\textnormal{inj}_{K_i}C[\mathds{G}_i]^{\oplus a\cdot |\mathds{B}_i\backslash \mathds{G}_i|^{-1}})|_{\mathfrak{B}}$ are injective $C$-representations of $\mathfrak{B}$ and $\kappa_i$ induces an isomorphism between their $\mathfrak{U}$-invariants (cf.\ Lemma \ref{injres}), we see that $\kappa_i$ must in fact be an isomorphism.

We now construct the morphism of diagrams.  Set $\psi_{K_i} := j_i \circ c_i$ and $\psi_{\mathfrak{B}} := (\bigoplus_\xi j_{\xi}^{\oplus a}) \circ c$.  We have 
$$\psi_{K_i} \circ \iota_i  \stackrel{\eqref{eqn1}}{=} \alpha_i \circ c \stackrel{\eqref{eqn2}}{=}  \kappa_i \circ \psi_{\mathfrak{B}},$$
and therefore we obtain the desired morphism of diagrams.  
\end{proof}

\subsection{Supersingular representations via homology}  
\label{sec:superc-repr-via}

Recall that a $G$-equivariant coefficient system $\mathcal{D}$ consists of $C$-vector spaces $\mathcal{D}_{\mathcal{F}}$ for every facet $\mathcal{F} \subset \mathscr{B}$, along with restriction maps for every inclusion of facets.  This data is required to have a compatible $G$-action such that each $\mathcal{D}_{\mathcal{F}}$ is a smooth $C$-representation of the $G$-stabilizer of $\mathcal{F}$.  The functor sending $\mathcal{D}$ to the quintuple $(\mathcal{D}_{x_0},\mathcal{D}_{x_1},\mathcal{D}_{\mathcal{C}},\iota_0,\iota_1)$, where the $\iota_i$ are the natural restriction maps, is an equivalence of categories between $G$-equivariant coefficient systems and diagrams (cf.\ \cite[\S 6.3]{KoziolXu}).

We let $\mathcal{D}_\Xi$ and $\mathcal{D}'$ denote the $G$-equivariant coefficient systems on $\mathscr{B}$ associated to $D_\Xi$ and $D'$, respectively.  The homology of $G$-equivariant coefficient systems gives rise to smooth $C$-representations of $G$, and we define
$$\pi := \textnormal{im}\left(H_0(\mathscr{B},\mathcal{D}_\Xi) \xrightarrow{\psi_*} H_0(\mathscr{B},\mathcal{D}')\right),$$
where $\psi_*$ denotes the map on homology induced by $\psi$.

\begin{thm}
\label{thm main-rk1}
Suppose $\chr F = 0$.  Then the $C$-representation $\pi$ of $G$ admits an irreducible, admissible, supercuspidal quotient.  
\end{thm}

\begin{proof}
We use language and notation from \cite{paskunas:diags} and \cite{KoziolXu}.  

\setcounter{step}{0}
\step{} 
The representation $\pi$ is nonzero. 

Fix a basis $v$ for $\chi^{-1}$. Let $\omega_{0,\iota_0(v)}$ denote the 0-chain with support $x_0$ satisfying $\omega_{0,\iota_{0}(v)}(x_0) = \iota_0(v)$
and let $\bar{\omega}_{0,\iota_0(v)}$ denote its image in $H_0(\mathscr{B},\mathcal{D}_\Xi)$.  
Set $\bar{\omega} := \psi_*(\bar{\omega}_{0,\iota_0(v)}) = \bar{\omega}_{0,\psi_{K_0}\circ\iota_0(v)} \in \pi \subset H_0(\mathscr{B},\mathcal{D}')$.  This is the image in $H_0(\mathscr{B},\mathcal{D}')$ of a $\mathcal{D}'_{x_0}$-valued 0-chain supported on $x_0$, and since the maps $\kappa_0,\kappa_1$ are isomorphisms and $\psi$ is injective, we have $\bar{\omega}\neq 0$ \cite[Lem.\ 5.7]{paskunas:diags}. We also note that therefore $\bar{\omega}_{0,\iota_0(v)} \ne 0$.

\step{} 
The representation $\pi$ is admissible.  

Since $\kappa_0, \kappa_1$ are isomorphisms, \cite[Prop.\ 5.10]{paskunas:diags} gives
$$\pi|_{\mathfrak{B}} \subset H_0(\mathscr{B},\mathcal{D}')|_{\mathfrak{B}} \cong \mathcal{D}'_{\mathcal{C}} \cong \bigoplus_\xi \textnormal{inj}_{\mathfrak{B}}\xi^{\oplus a},$$
which by Lemma \ref{pask2} implies $\pi^{\mathfrak{U}} \hookrightarrow \bigoplus_\xi \xi^{\oplus a}$, so that $\pi$ is admissible.

\step{}
The $\mathcal{H}_C$-module $\pi^{\mathfrak{U}}$ contains $\Xi$.  

The element $\bar{\omega}_{0,\iota_0(v)}\in H_0(\mathscr{B},\mathcal{D}_\Xi)$ is $\mathfrak{U}$-invariant and stable by the action of $\mathcal{H}_C$, and the vector space it spans is isomorphic to $\Xi$ as an $\mathcal{H}_C$-module (for all of this, see the proof of \cite[Prop.\ 7.3]{KoziolXu}).  Since $\psi_*$ is $G$-equivariant, the same is true for $\bar{\omega}\in \pi$.

\step{}
The vector $\bar{\omega}$ generates $\pi$.  

Since $\bar{\omega}_{0,\iota_0(v)}$ generates $H_0(\mathscr{B},\mathcal{D}_\Xi)$ as a $G$-representation and $\psi_*$ is $G$-equivariant, $\bar{\omega}$ generates $\pi$ as a $G$-representation.

\step{} 
We construct the quotient $\pi'$ and list its properties.

By the previous step, the representation $\pi$ is generated by $\bar{\omega}$.  Proceeding as in the end of the proof of Theorem \ref{thm main}, we see that any irreducible quotient of $\pi = \langle G\cdot \bar{\omega}\rangle$ is admissible (since $\chr F = 0$, and such quotients exist by Zorn's lemma).  Let $\pi'$ be any such quotient.

\step{} 
We prove $\pi'$ is supercuspidal.  

Since $\bar{\omega}$ generates $\pi$, its image in $\pi'$ is nonzero.  Thus, we obtain an injection of $\mathcal{H}_C$-modules $\Xi \cong C\bar{\omega} \hookrightarrow (\pi')^{\mathfrak{U}}$, and supercuspidality follows from Proposition \ref{prop ss=sc}.
\end{proof}

\begin{corollary}
\label{cor rk1}
Suppose $\chr F = 0$ and $\mathbf{G}$ is a connected reductive $F$-group of relative semisimple rank 1.  Then $G$ admits an irreducible admissible supercuspidal $C$-representation.  
\end{corollary}

\begin{proof}
By the reductions in \S \ref{sec:prelim-rk1}, it suffices to assume $\mathbf{G}$ is absolutely simple and simply connected, and to construct a supercuspidal $C$-representation on which $Z(G)$ acts trivially.  Since the center of $G$ is finite, it is contained in $\mathfrak{B} \cap Z = Z_0$.  Hence, taking $\Xi$ to be associated to $(\mathbf{1}_{\mathds{Z}},J)$, where $\mathbf{1}_{\mathds{Z}}$ is the trivial character of ${\mathds{Z}}$ and $J\neq \emptyset, \{s_0,s_1\}$ (noting that $S_{\mathbf{1}_{\mathds{Z}}} = \{s_0,s_1\}$), Theorem \ref{thm main-rk1} produces an irreducible admissible supercuspidal $C$-representation $\pi'$ with trivial action of the center.  This gives the claim.
\end{proof}

\begin{rk}
The construction of $\pi'$ above shares some similarities with the construction in \S \ref{sec:reduction-rank-1-pgln}.  Therein, supercuspidal representations are constructed as subquotients of $C^\infty(\Gamma\backslash G,C) \cong \Ind_{\Gamma}^G\mathbf{1}_\Gamma$, where $\Gamma$ is a discrete, cocompact subgroup of $G$ and $\mathbf{1}_\Gamma$ denotes the trivial character of $\Gamma$.  Taking $\Gamma$ to be torsion-free, we use the Mackey formula to obtain
$$(\Ind_{\Gamma}^G\mathbf{1}_\Gamma)|_{K_i}  \cong \bigoplus_{\Gamma\backslash G /K_i}\Ind_{\{1\}}^{K_i}\mathbf{1} \cong \textnormal{inj}_{K_i}C[\mathds{G}_i]^{\oplus a_i'}$$
where $a_i' = |\Gamma \backslash G/ K_i|$.  (The last isomorphism follows from the fact that $\Ind_{\{1\}}^{K_i}\mathbf{1}$ is injective, by \cite[\S I.5.9 b)]{Viglivre}, and $(\Ind_{\{1\}}^{K_i}\mathbf{1})^{K_i^+} \cong C[\mathds{G}_i]$; we may then proceed as in the proof of Lemma \ref{injres}.)  The construction of Theorem \ref{thm main-rk1} produces supercuspidal representations as subquotients of $H_0(\mathscr{B},\mathcal{D}')$, for which we have
$$H_0(\mathscr{B},\mathcal{D}')|_{K_i} \cong \textnormal{inj}_{K_i}C[\mathds{G}_i]^{\oplus a_i},$$
where $a_i = a\cdot |\mathds{B}_i\backslash \mathds{G}_i|^{-1}$ (cf.\ \cite[Prop.\ 5.10]{paskunas:diags}).  
\end{rk}

\section{Supersingular representations of \texorpdfstring{$\PGL_n(D)$}{PGL\_n(D)}}
\label{sec:supers-repr-pgl_n}

In this section we verify Theorem~\ref{thm:main} when $\bG = \bPGL_n(D)$, where $n \ge 2$ and $D$ a central division algebra over $F$.
In particular, this deals with the first exceptional case in Theorem~\ref{thm main}.

\subsection{Notation and conventions}
\label{sec:notation-conventions}

Throughout Section~\ref{sec:supers-repr-pgl_n}, we let $\qpb$ denote a fixed algebraic closure of $\qp$, with ring of integers $\zpb$ and residue
field $\fpb$. We normalize the valuation $\val$ of $\qpb$ such that $\val(p) = 1$.

Let $D$ denote a central division algebra over $F$ of dimension $d^2$. 
Let $B = ZU$ denote the upper-triangular Borel subgroup of $\GL_n(D)$ with diagonal minimal Levi
subgroup $Z \cong (D\s)^n$ and unipotent radical $U$. Let $T \cong (F\s)^n$ denote the diagonal
maximal split torus, $\cN$ its normalizer in $\GL_n(D)$, and $U\op$ the lower-triangular unipotent matrices. 

Let $\O_D$ denote the ring of integers of $D$, $\m_D$ the maximal ideal of $\O_D$, and $k_D$ the
residue field, so $[k_D : k_F] = d$. Let $D(1) := 1+\m_D$, so $D(1) \lhd D\s$.  
Let $\val_D : D\s \onto \Z$ denote the normalized valuation of $D$.
Let $I(1)$ denote the pro-$p$ Iwahori subgroup
\begin{equation*}
  I(1) := \{ g \in \GL_n(\O_D) : \text{$\o g \in \GL_n(k_D)$ is upper-triangular unipotent} \}.
\end{equation*}

For any field $K$ let $\Gamma_K$ denote the absolute Galois group for a choice of separable closure.
If $K'/K$ is a finite separable extension, then $\Gamma_{K'}$ is a subgroup of $\Gamma_K$, up to conjugacy,
hence the restriction of a $\Gamma_K$-representation to $\Gamma_{K'}$ is well defined up to isomorphism.

If $K/\qp$ is finite we let $I_K$ denote the inertia subgroup of $\Gamma_K$ and $k_K$
the residue field of $K$.
If $\rho : \Gamma_K \to \GL_n(\qpb)$ is de Rham and $\tau : K \to \qpb$ is continuous, we let
$\HT_\tau(\rho)$ denote multi-set of $\tau$-Hodge--Tate weights. We normalize Hodge--Tate weights so
that the cyclotomic character $\ve$ has $\tau$-Hodge--Tate weight $-1$ for any $\tau$.
We let $\WD(\rho)$ denote the associated Weil--Deligne representation of $W_K$ over $\qpb$
(defined by Fontaine, cf.\ Appendix B.1 of \cite{MR1639612}).

We normalize local class field theory so that uniformizers correspond to geometric Frobenius
elements under the local Artin map. Let $\rec_F$ denote the local Langlands correspondence from isomorphism
classes of irreducible smooth representations of $\GL_n(F)$ over $\C$ to isomorphism classes of $n$-dimensional Frobenius
semisimple Weil--Deligne representations of $W_F$ over $\C$. (See \cite{ht}.) 

If $L$ is a global field, we let $|\cdot|_L$ denote the normalized absolute value of $\A_L$.

\subsection{On the Jacquet--Langlands correspondence}
\label{sec:jacq-langl-corr}

We recall some basic facts about the representation theory of $\GL_n(D)$ and the local Jacquet--Langlands correspondence.
All representations in this section will be smooth and over $\C$. 

For a finite-dimensional central simple algebra $A$ let $\Nrd : A\s \to Z(A)\s$ (or $\Nrd_A$ for
clarity) denote the reduced norm.
Let $\nu$ denote the smooth character $|\Nrd|_F$ of $\GL_m(D)$ for any $m$. 
If $\pi_i$ are smooth representations of $\GL_{n_i}(D)$, let $\pi_1 \times \cdots \times \pi_r$
denote the normalized parabolic induction of $\pi_1 \otimes \cdots \otimes \pi_r$ to $\GL_{\sum n_i}(D)$.
In particular these notions
also apply to general linear groups over $F$ (by setting $D = F$).

We will say that a representation is \emph{essentially unitarizable} if some twist of it is unitarizable.

The Jacquet--Langlands correspondence \cite{DKV} is a canonical bijection $\JL$ between 
irreducible essentially square-integrable representations of $\GL_n(D)$ and 
irreducible essentially square-integrable representations of $\GL_{nd}(F)$ that
is compatible with character twists and preserves central characters. (For short, we say ``square-integrable'' instead
of ``square-integrable modulo center''.)

On the other hand, Badulescu \cite{badu-jlc} defined a map $|\LJ_{\GL_n(D)}|$ in the other direction, 
from irreducible essentially unitarizable representations of $\GL_{nd}(F)$ to
irreducible essentially unitarizable representations of $\GL_{n}(D)$ or zero, which in general
is neither injective nor surjective. (More precisely,
\cite{badu-jlc} only considers unitarizable representations, but we can extend it by twisting.)
In the split case $|\LJ_{\GL_{n}(F)}|$ is the identity.
It follows from Thm.\ 2.2 and Thm.\ 2.7(a) in \cite{badu-jlc} that $|\LJ_{\GL_n(D)}|(\JL(\pi)) \cong \pi$
for any essentially square-integrable representation $\pi$ of $\GL_n(D)$.

If $\rho$ is a supercuspidal representation of $\GL_m(F)$ and $\ell \ge 1$, then
$Z^u(\rho,\ell)$ is by definition the unique irreducible quotient of
$\rho \nu^{(1-\ell)/2} \times \rho \nu^{(3-\ell)/2} \times \cdots \times \rho \nu^{(\ell-1)/2}$.
It is an essentially square-integrable representation of $\GL_{m\ell}(F)$. All essentially square-integrable representations of $\GL_n(F)$ arise in this way,
for some decomposition $n = m\ell$.

If $\rho'$ is a supercuspidal representation of $\GL_m(D)$, we can write $\JL(\rho') \cong Z^u(\rho,s)$
for some supercuspidal representation $\rho$ and integer $s \ge 1$. 
Then $Z^u(\rho',\ell)$ is by definition the unique irreducible quotient of
$\rho' \nu^{s(1-\ell)/2} \times \rho' \nu^{s(3-\ell)/2} \times \cdots \times \rho' \nu^{s(\ell-1)/2}$.
It is an essentially square-integrable representation of $\GL_{m\ell}(D)$. All essentially square-integrable representations of $\GL_n(D)$ arise in this way,
for some decomposition $n = m\ell$ (a result of Tadi\'c, cf.\ \cite[\S2.4]{badu-jlc}). 
Moreover, $\JL(Z^u(\rho',\ell)) \cong Z^u(\rho,\ell s)$ \cite[\S3.1]{badu-jlc}.

If $\pi$ is a smooth representation of $\GL_n(D)$ let $\pi_U$ denote its (unnormalized) Jacquet
module. The following lemma was proved earlier, cf.\ Remark~\ref{rk:pro-p-invts-jacquet}.

\begin{lm}\label{lm:pro-p-invts-jacquet2}
  Suppose that $\pi$ is an admissible representation of $\GL_n(D)$ over $\C$.
  Then the natural map $p_U : \pi \to \pi_U$ induces an isomorphism
  $\pi^{I(1)} \to (\pi_U)^{Z \cap I(1)}$.
\end{lm}

The following results will be needed in Section~\ref{sec:lift-non-supers}.

\begin{lm}\label{lm:glnd-iwahori}
  Suppose that $\Pi$ is an irreducible generic smooth representation of $\GL_{nd}(F)$ over $\C$ that
  is essentially unitarizable and such that the representation $\pi := |\LJ_{\GL_n(D)}|(\Pi)$ of $\GL_n(D)$ is non-zero.
  If $\pi^{I(1)} \ne 0$, then there exist irreducible representations $\rho_1'$, \dots, $\rho_n'$ of
  $D\s/D(1)$ such that $\pi$ is a subquotient of $\rho_1' \times \cdots \times \rho_n'$ and
  $\rec_F(\Pi)|_{W_F} \cong \bigoplus_{i = 1}^n \rec_F(\JL(\rho_i'))|_{W_F}$.
\end{lm}

\begin{proof}
  After a twist we may assume that $\Pi$ is unitarizable.  As $\Pi$ is moreover generic, we know
  that $\Pi \cong \sigma_1 \nu^{\alpha_1} \times \cdots \times \sigma_r \nu^{\alpha_r}$ for some
  square-integrable $\sigma_i$ of $\GL_{n_i}(F)$ and real numbers $\alpha_i \in (-\frac 12,\frac
  12)$ satisfying $\alpha_i + \alpha_{r+1-i} = 0$ and $\sigma_i = \sigma_{r+1-i}$ if $\alpha_i \ne
  0$ (see e.g.\ \cite[Lemma I.3.8]{ht}). Since $|\LJ_{\GL_n(D)}|(\Pi) \ne 0$ by assumption, it follows that $d \mid
  n_i$ for all $i$ and $\pi = |\LJ_{\GL_n(D)}|(\Pi) \cong \sigma_1' \nu^{\alpha_1} \times \cdots
  \times \sigma_r' \nu^{\alpha_r}$, where $\sigma_i'$ is the square-integrable representation of
  $\GL_{n_i/d}(D)$ such that $\JL(\sigma_i') \cong \sigma_i$. (See \cite[\S3.5]{badu-jlc}.)
  Let $n_i' := n_i/d$.
  
  From $\pi^{I(1)} \ne 0$ and Lemma~\ref{lm:pro-p-invts-jacquet2} it follows that the supercuspidal
  support of $\pi$ is a tame representation of $Z$ (up to conjugacy), so each $\sigma_i'$ is of the
  form $Z^u(\rho_i'',n_i')$, where $\rho_i''$ is an irreducible representation of $D\s/D(1)$. We
  write $\JL(\rho_i'') \cong Z^u(\rho_i,e_i)$ with $\rho_i$ irreducible supercuspidal, so $\sigma_i
  \cong Z^u(\rho_i, e_i n_i')$. In particular, $\pi$ is a subquotient of the normalized induction of
  $\bigotimes_{1 \le i \le r, 0 \le j \le n_i'-1} \rho_i'' \nu^{\alpha_i + e_i((n_i'-1)/2-j)}$.  On
  the other hand, $\Pi$ is a subquotient of the normalized induction of $\bigotimes_{1 \le i \le r,
    0 \le j \le e_i n_i'-1} \rho_i \nu^{\alpha_i + (e_i n_i'-1)/2-j}$.  As $\rec_F(\Pi)|_{W_F}$ only
  depends on the supercuspidal support of $\Pi$ (see the paragraph before Thm.\ VII.2.20 in
  \cite{ht}), we obtain
  \begin{equation*}
    \rec_F(\Pi)|_{W_F} \cong \bigoplus_{1 \le i \le r, 0 \le j \le e_i n_i'-1}
    |\cdot|_F^{\alpha_i + (e_i n_i'-1)/2-j}\rec_F(\rho_i)|_{W_F}.
  \end{equation*}
  Similarly, $\rec_F(\JL(\rho_i''))|_{W_F} \cong \bigoplus_{k=0}^{e_i-1} |\cdot|_F^{(e_i-1)/2-k}\rec_F(\rho_i)|_{W_F}$.
  Denoting by $\rho_1'$, \dots, $\rho_n'$ the representations $\rho_i'' \nu^{\alpha_i +
    e_i((n_i'-1)/2-j)}$ in any order, a straightforward computation confirms that $\bigoplus_{i = 1}^n
  \rec_F(\JL(\rho_i'))|_{W_F} \cong \rec_F(\Pi)|_{W_F}$.
\end{proof}

We now recall a result of Bushnell--Henniart concerning explicit functorial transfers of irreducible representations of
$D\s/D(1)$.  An \emph{admissible tame pair} $(E/F,\zeta)$ consists of an unramified extension of degree $f$ dividing $d$, and
a tamely ramified smooth character $\zeta : E\s \to \C\s$ such that all $\Gal(E/F)$-conjugates of $\zeta$ are distinct.  In
that case, after choosing an $F$-embedding of $E$ into $D$ (which is unique up to conjugation by $D^\times$), $B := Z_D(E)$ is a central simple $E$-algebra of dimension $e^2$, where $e := d/f$. Define a smooth character
$\Lambda : B\s (1+\m_D) \to \C\s$ by declaring it to be $\zeta \circ \Nrd_{B}$ on $B\s$ and trivial on $1+\m_D$. Then
define $\pi_D(\zeta) := \Ind_{B\s (1+\m_D)}^{D\s} \Lambda$ is an irreducible representation of $D\s/D(1)$ (of dimension $f$).

\begin{prop}\label{prop:bushnell-henniart}\

  \begin{enumerate}
  \item Any irreducible representation of $D\s/D(1)$ is isomorphic to $\pi_D(\zeta)$ for some admissible tame pair $(E/F,\zeta)$.
  \item The element $\varpi \in F$ acts as the scalar $\zeta(\varpi)^e$ on $\pi_D(\zeta)$.
  \item If $(E/F,\zeta)$ is an admissible tame pair, then
    \begin{equation*}
      \rec_F(\JL(\pi_D(\zeta))) \cong \Sp_{e}(\Ind_{W_{E}}^{W_F} (\eta_{E}^{e(f-1)} \zeta)),
    \end{equation*}
    where $\eta_{E}$ is the unramified quadratic character of $E\s$.
  \end{enumerate}
\end{prop}

We recall that the special Weil--Deligne representation $\Sp_e(\sigma)$, for $\sigma$ an irreducible representation of $W_F$,
is indecomposable and satisfies $\Sp_e(\sigma)|_{W_F} \cong \bigoplus_{k=0}^{e-1} \sigma |\cdot|_F^{\frac{e-1}2-k}$.

\begin{proof}
  For (i), see \cite[\S1.5]{bushnell-henniart-level0}. Part (ii) follows from the definition. Part (iii) is 
  the main result of \cite{bushnell-henniart-level0}.
\end{proof}

\subsection{On lifting non-supersingular Hecke modules}
\label{sec:lift-non-supers}

Let $\HH := \HH(\GL_n(D),I(1))$ the corresponding pro-$p$ Iwahori--Hecke algebra over $\Z$
\cite{Vig1} and for a commutative ring $R$ let $\HH_R := \HH \otimes R$. Similarly we define
$\HH_Z := \HH(Z,Z \cap I(1))$ and $\HH_{Z,R} := \HH_Z \otimes R$. Note that the pro-$p$ Iwahori
subgroup $Z \cap I(1)$ is normal in $Z$. All Hecke modules we will consider are right modules.
A finite-dimensional $\HH_{\qpb}$-module is said to be \emph{integral} if it arises by base change
from a $\HH_{\zpb}$-module that is finite free over $\zpb$.

Let $W(1) := \mathcal N/Z \cap I(1)$, $\Lambda(1) := Z/Z \cap I(1)$, and define monoids
\begin{equation*}
  Z^+ := \{\diag(\delta_1,\dots,\delta_n) \in Z : \val_D(\delta_1) \ge \cdots \ge \val_D(\delta_n)\}
\end{equation*}
and $\Lambda(1)^+ := Z^+/Z \cap I(1)$. 

We recall that $\HH$ has an Iwahori--Matsumoto basis $T_w$ for $w \in W(1)$ and 
a Bernstein basis $E_w$ for $w \in W(1)$, which in fact depends on a choice of spherical orientation.
We choose our spherical orientation such that $E_w = T_w$ for $w \in \Lambda(1)^+$. (This is possible
by \cite[Ex.\ 5.30]{Vig1}. It is the opposite of our convention in \S\ref{sec:bernstein-elements}.)
Similarly, $\HH_Z$ has basis $T_w^Z$ for $w \in \Lambda(1)$.

For $w \in W(1)$ we have integers $q_w \in q^{\Z_{>0}}$, as recalled in \S\ref{sec:bernstein-elements}.
(Note that our base alcove $\cC$ is the one fixed pointwise by $I(1)$.)

\begin{lm}\label{lm:length-formula}
  For $z = \diag(\delta_1,\dots,\delta_n) \in Z$ with $\delta_i \in D\s$ we have $$q_z = q^{d\sum_{i < j} |\val_D(\delta_i)-\val_D(\delta_j)|}.$$
\end{lm}

\begin{proof}
  As the Iwahori--Hecke algebra has equal parameters $q^d$ we deduce that $q_z = q^{d \ell(z)}$, where $\ell$ is the length
  function relative to the alcove $\cC$. By using the action of the finite Weyl group $\cN/Z$ and the first length formula in
  \cite[Cor.\ 5.11]{Vig1}, we may assume that $z \in Z^+$. By \cite[\S3.9]{Vig1} we then have $q_z = (I(1)zI(1):I(1)) =
  (I(1):I(1) \cap z I(1)z^{-1}) = (U_0 : z U_0 z^{-1})$, where $U_0 := U \cap I(1)$. Hence $q_z = q^{d\sum_{i < j}
    (\val_D(\delta_i)-\val_D(\delta_j))}$, as required.
\end{proof}

Let $W_0 \cong S_n$ denote the Weyl group of $T$.  Recall from \cite[\S5, \S1.3]{Vigss} that $\cA_0(\Lambda_T)$ is the free
module with basis $E_{\mu(\varpi)}$ for $\mu \in \Lambda_T := X_*(T)$ and that the central subalgebra
$\cZ_T := \cA_0(\Lambda_T)^{W_0}$ of $\HH$ has a basis consisting of the sums $\sum_\mu E_{\mu(\varpi)}$ with $\mu$ running
over the $W_0$-orbits in $X_*(T)$. For $I \subset \{1,\dots,n\}$ let $E_I := E_{\mu_I(\varpi)}$, where
$\mu_I \in X_*(T) \cong \Z^n$ is defined by $\mu_{I,i} = 1$ if $i \in I$ and $\mu_{I,i} = 0$ otherwise. For $1 \le i \le n$
let $\mathscr Z_i := \sum_{I, |I| = i} E_I$.  By induction and \cite[Cor.\ 5.28]{Vig1} we see that the algebra $\cZ_T$ is
generated by $\mathscr Z_1$, \dots, $\mathscr Z_{n-1}$, $\mathscr Z_n^{\pm 1}$.
  
The following lemma follows from \cite[Prop.\ 6.9]{Vigss}.

\begin{lm}\label{lm:supersing-crit}
  A finite-dimensional $\HH_{\fpb}$-module $M$ is supersingular if and only if the action of $\mathscr Z_i$ on $M$ is
  nilpotent for all $1 \le i \le n-1$.
\end{lm}

\begin{lm}\label{lm:ollivier-vigneras}
  There exists a unique injective algebra homomorphism $\wt\theta : \HH_{Z,\qpb} \to \HH_{\qpb}$ such that
  $\wt\theta(T_w^Z) = T_w$ for all $w \in \Lambda(1)^+$. We have
  \begin{equation}\label{eq:8}
    E_I = q^{d^2(\sum_{i \in I} i - \binom{|I|+1}2)} \wt\theta(T^Z_{\mu_I(\varpi)}).
  \end{equation}
\end{lm}

\begin{proof}
  The first assertion follows from \cite[\S2.5.2, Rk.\ 2.20]{OV}.
  We claim that for any $\mu \in X_*(T)$,
  \begin{equation}\label{eq:7}
    E_{\mu(\varpi)} = q^{d^2 \sum_{r<s : \mu_r < \mu_s}(\mu_s-\mu_r)} \wt\theta(T^Z_{\mu(\varpi)}),
  \end{equation}
  which implies~\eqref{eq:8} by taking $\mu = \mu_I$.

  Note that $X_*(T)^+ = \{\mu \in X_*(T) : \mu_1 \ge \cdots \ge \mu_n\}$. If $\mu \in X_*(T)^+$, 
  then $\mu(\varpi) \in Z^+$ and hence $E_{\mu(\varpi)} = T_{\mu(\varpi)}$ and formula~(\ref{eq:7})
  holds.  In general, choose $\mu' \in X_*(T)^+$ such that $\mu+\mu' \in X_*(T)^+$. 
  Then formula~(\ref{eq:7}) follows easily from the following three assertions:
  (1) $T^Z_{\mu(\varpi)} T^Z_{\mu'(\varpi)} = T^Z_{\mu(\varpi)\mu'(\varpi)}$;
  (2) $E_{\mu(\varpi)} E_{\mu'(\varpi)} = (q_{\mu(\varpi)} q_{\mu'(\varpi)}
  q_{\mu(\varpi)\mu'(\varpi)}^{-1})^{1/2} E_{\mu(\varpi)\mu'(\varpi)}$ in the notation of \cite[\S4.4]{Vig1}, 
  where we take the positive square root; 
  and (3) Lemma~\ref{lm:length-formula}.
  Assertion (1) is clear and assertion (2) is \cite[Cor.\ 5.28]{Vig1}. 
\end{proof}

The following simple and presumably well-known lemma will be used below.

\begin{lm}\label{lm:valuations}
  Suppose that $\rho : W_F \to \GL_n(\qpb)$ is a smooth representation. Then for any $\gamma \in W_F$ the valuations of the
  eigenvalues of $\rho(\gamma)$ depend only on the image of $\gamma$ in $W_F/I_F \cong \Z$.
\end{lm}

\begin{proof}
  Fix a geometric Frobenius element $\Frob_F \in W_F$, and let $v_1 \le \cdots \le v_n$ denote the valuations of the
  eigenvalues of $\rho(\Frob_F)$. We need to show that the eigenvalues of $\rho(\Frob_F^r g)$ have valuations $r v_1 \le
  \cdots \le r v_n$ for any $g \in I_F$.  As $\rho(I_F)$ is finite and normalized by $\rho(\Frob_F)$, we see that $\rho(\Frob_F)^{\ell r}$ and $\rho(I_F)$
  commute for some $\ell \ge 1$, so $\rho(I_F)$ preserves the generalized eigenspaces of $\rho(\Frob_F^{\ell r})$. Hence the
  valuations of the eigenvalues of $\rho(\Frob_F^{\ell r} g)$ are independent of $g \in I_F$, and the claim follows by
  passing to $\ell$-th powers.
\end{proof}

We now fix an isomorphism $\imath : \qpb \congto \C$.

\begin{prop}\label{prop:non-supersingular-reduction-and-frob-evals}
  Suppose that $\Pi$ is an irreducible generic smooth representation of $\GL_{nd}(F)$ over $\C$ that
  is essentially unitarizable and such that the representation $\pi := |\LJ_{\GL_n(D)}|(\Pi)$ of $\GL_n(D)$ is non-zero.
  Suppose that $\imath^{-1}(\pi^{I(1)})$ is a non-zero integral $\HH_{\qpb}$-module with non-supersingular
  reduction, and let $v_1 \le \cdots \le v_{nd}$ denote the valuations of the eigenvalues of a
  geometric Frobenius on $\imath^{-1}(\rec_F(\Pi))$. Then there exists $1 \le j \le n-1$ such that
  \begin{equation*}
    \sum_{i=1}^{jd} v_i = -\frac{d^2 j(n-j)}2 \val(q).
  \end{equation*}
\end{prop}

\begin{proof}

\setcounter{step}{0}
  \step{}\label{step:evals-Zi}
  We compute the action of $\mathscr Z_1$, \dots, $\mathscr Z_n$ on 
  the Hecke module $\imath^{-1}(\pi^{I(1)})$ and show in particular that it is scalar.

  Note by~Lemma \ref{lm:glnd-iwahori} that $\pi^{I(1)}$ is a subquotient of $(\rho_1' \times \cdots \times
  \rho_n')^{I(1)}$ for some irreducible representations $\rho_i'$ of $D\s/D(1)$, and $\rho_1' \times \cdots \times
  \rho_n' \cong \Ind_B^{\GL_n(D)}(\rho_1' \nu^{d(n-1)/2} \otimes \cdots \otimes \rho_n' \nu^{-d(n-1)/2})$
  (unnormalized induction).
  By \cite[Prop.\ 4.4]{OV} we have
  \begin{align}
    \imath^{-1}(\rho_1' \times \cdots \times \rho_n')^{I(1)} &\cong \imath^{-1}(\rho_1' \nu^{d(n-1)/2} \otimes \cdots 
    \otimes \rho_n' \nu^{-d(n-1)/2})^{Z \cap I(1)} \otimes_{\HH_{Z,\qpb}, \wt\theta} \HH_{\qpb}, \label{eq:9}   
  \end{align}   
  where we used the homomorphism $\wt\theta$ of Lemma~\ref{lm:ollivier-vigneras}.  

  By Proposition~\ref{prop:bushnell-henniart}(i) we can write $\rho_i' \cong \pi_D(\zeta_i)$ for some
  admissible tame pair $(F_i/F, \zeta_i)$. We let $f_i := [F_i : F]$ and $e_i := d/f_i$.
  Let $\zeta'_i := \imath^{-1}(\zeta_i)$.
  From equations (\ref{eq:8}), (\ref{eq:9}) and Proposition~\ref{prop:bushnell-henniart}(ii) we deduce that $\mathscr Z_j$ acts on $\imath^{-1}(\pi^{I(1)})$ as the scalar
  \begin{align}
    \lambda_j &:= \sum_{|I| = j} \bigg(q^{d^2(\sum_{i \in I} i - \binom{j+1}2)} q^{d^2((n+1)j/2 - \sum_{i \in I} i)} \prod_{i \in I} \zeta'_i(\varpi)^{-e_i}\bigg)\label{eq:6} \\
    & = q^{-d^2\binom{j}2} \sum_{|I| = j} \bigg(q^{d^2(n-1)j/2} \prod_{i \in I} \zeta'_i(\varpi)^{-e_i}\bigg).\notag
  \end{align}

  \step{}\label{step:end-of-pf}
  We complete the proof.  By
  assumption, the Hecke module $\imath^{-1}(\pi^{I(1)})$ is integral, so $\lambda_i \in \zpb$ for all $i$ and
  $\lambda_n \in \zpb\s$. Moreover, as the reduction of $\imath^{-1}(\pi^{I(1)})$ is non-supersingular we deduce
  by Lemma~\ref{lm:supersing-crit} that $\lambda_{n-j} \in \zpb\s$ for some $1 \le j \le n-1$.

  From now on assume for convenience that the $\zeta'_i$ are ordered such that the sequence
  $\val(\zeta'_i(\varpi)^{-e_i})$ is non-increasing. Consider the polynomial $\prod_{i=1}^n
  (1-q^{d^2(n-1)/2} \zeta'_i(\varpi)^{-e_i} X)$.  By~(\ref{eq:6}) its Newton polygon is defined by the points
  $(i,\val(\lambda_i) + d^2 \binom i2 \val(q))$ for $0 \le i \le n$. 
  From $\lambda_{n-j} \in \zpb\s$, $\lambda_i \in \zpb$, and the convexity of the quadratic function $x(x-1)/2$
  we deduce that $(n-j,d^2 \binom {n-j}2 \val(q))$ is a vertex of the Newton polygon.
  It follows for the sum of the largest $j$ root valuations that
  \begin{equation}\label{eq:1}
    \sum_{i=1}^j \val(q^{d^2(n-1)/2} \zeta'_i(\varpi)^{-e_i}) = d^2\bigg(\binom n2-\binom {n-j}2\bigg)\val(q).
  \end{equation}
  Again by convexity we obtain the root valuation bounds
  \begin{align}
    \val(q^{d^2(n-1)/2} \zeta'_i(\varpi)^{-e_i}) &\ge d^2(n-j) \val(q) \quad\forall i \le j,\label{eq:3} \\
    \val(q^{d^2(n-1)/2} \zeta'_i(\varpi)^{-e_i}) &\le d^2(n-j-1) \val(q) \quad\forall i > j.\label{eq:4}
  \end{align}

  From Lemma~\ref{lm:glnd-iwahori} and Proposition~\ref{prop:bushnell-henniart}(iii) we see that
  \begin{equation*}
    \rec_F(\Pi)|_{W_F} \cong \bigoplus_{i=1}^n \bigoplus_{k=0}^{e_i-1}
    \Ind_{W_{F_i}}^{W_F} (\eta_{F_i}^{e_i(f_i-1)} \zeta_i) |\cdot|_F^{(e_i-1)/2-k}.
  \end{equation*}
  If $\Frob_F$ denotes a geometric Frobenius of $W_F$, then $\Frob_F^{f_i}$ is a geometric Frobenius
  of $W_{F_i}$. We see that all eigenvalues of $\Frob_F$ on $\Ind_{W_{F_i}}^{W_F} (\eta_{F_i}^{e_i(f_i-1)}
  \zeta'_i)$ have valuation $\frac 1{f_i}\val(\zeta'_i(\Frob_F^{f_i})) = \frac
  1{f_i}\val(\zeta'_i(\varpi))$.  Hence, for $i \le j$ and $0 \le k \le e_i-1$ all eigenvalues of
  $\Frob_F$ on $\Ind_{W_{F_i}}^{W_F} (\eta_{F_i}^{e_i(f_i-1)} \zeta'_i) |\cdot|_F^{(e_i-1)/2-k}$ have
  valuation
  \begin{multline*}
    \frac 1{f_i}\val(\zeta'_i(\varpi)) - \Big(\frac{e_i-1}2-k\Big) \val(q) \le 
    \frac 1d\val(\zeta'_i(\varpi)^{e_i})+\Big(\frac{e_i-1}2\Big) \val(q) \\
    < d\Big(\frac{n-1}2-(n-j)\Big)\val(q) + \frac d2 \val(q) = \frac{d(2j-n)}2\val(q),
  \end{multline*}
  where we used~\eqref{eq:3} and that $e_i-1 < d$. Similarly, for $i > j$ and $0 \le k \le e_i-1$
  we find that the eigenvalues of $\Frob_F$ on $\Ind_{W_{F_i}}^{W_F} (\eta_{F_i}^{e_i(f_i-1)} \zeta'_i) |\cdot|_F^{(e_i-1)/2-k}$
  have valuation greater than $\frac{d(2j-n)}2\val(q)$.
  Therefore, from~\eqref{eq:1} we deduce that
  \begin{multline*}
    \sum_{i=1}^{jd} v_i = \sum_{i=1}^j \sum_{k=0}^{e_i-1} f_i \bigg(\frac 1{f_i}\val(\zeta'_i(\varpi)) - \Big(\frac{e_i-1}2-k\Big) \val(q)\bigg) = \sum_{i=1}^j \val(\zeta'_i(\varpi)^{e_i}) \\
    = d^2\bigg(\binom{n-j}2 - \binom n2+\frac{j(n-1)}2\bigg) \val(q) = -\frac{d^2 j(n-j)}2 \val(q).
  \end{multline*}
\end{proof}

\subsection{A reducibility lemma}
\label{sec:reducibility-lemma}

Let $F_0$ denote the maximal absolutely unramified intermediate field of $F/\qp$. The following lemma
generalizes \cite[Prop.\ 4.5.2]{egh}, which dealt with regular crystalline Galois representations.

\begin{lm}\label{lm:hodge-newton-reducibility}
  Suppose that $\rho : \Gamma_F \to \GL_n(\qpb)$ is a de Rham Galois representation.  Let $v_1 \le \cdots \le v_n$ denote the
  valuations of the eigenvalues of a geometric Frobenius element acting on $\WD(\rho)$, and for each embedding
  $\tau : F \to \qpb$ let $h_{\tau,1} \le \cdots \le h_{\tau,n}$ denote the $\tau$-Hodge--Tate weights of $\rho$.  Then
  $\sum_{i = 1}^j v_i \ge [F:F_0]^{-1} \sum_{i=1}^j \sum_{\tau : F \to \qpb} h_{\tau,i}$ for any $0 \le j \le n$.
  
  Suppose that $h_{\tau,1} < h_{\tau,n}$ for some $\tau$ and that for some $1 \le j \le n-1$ we have
  $\sum_{i = 1}^j v_i = [F:F_0]^{-1} \sum_{i=1}^j \sum_{\tau : F \to \qpb} h_{\tau,i}$. Then $\rho$ is reducible.
\end{lm}

\begin{proof}
  We first choose $E/\qp$ a sufficiently large finite subextension of $\qpb/\qp$, so that in particular $\rho$ can be defined
  over $E$ and all embeddings $\tau$ have image contained in $E$.  Choose $F'/F$ a finite Galois extension over which $\rho$
  becomes semistable. Let $D := \Dst(\rho|_{\Gamma_{F'}})$ be the covariantly associated free $F'_0 \otimes_{\qp} E$-module,
  equipped with actions of $\vp$, $N$, $\Gal(F'/F)$, where $F'_0$ denotes the maximal absolutely unramified intermediate
  field of $F'/\qp$. As usual, we write $D \cong \bigoplus_{\sigma : F'_0 \to E} D_\sigma$.  Fix any embedding
  $\sigma_0 : F'_0 \to E$ and let $f' := [F'_0:\qp]$. Note that $\vp^{f'}$ acts linearly on $D$ and stabilizes
  each $D_\sigma$.

  By construction of $\WD(\rho)$ and Lemma~\ref{lm:valuations}, the eigenvalues of $\vp^{f'}$ on $D_{\sigma_0}$ have valuations
  $r v_1 \le \cdots \le r v_n$, where $r := [F'_0:F_0]$.  For any $0 \le j \le n$, choose a $\vp^{f'}$-stable $E$-subspace
  $D'_{\sigma_0} \subset D_{\sigma_0}$ of dimension $j$ such that the eigenvalues of $\vp^{f'}$ on $D'_{\sigma_0}$ have
  valuations $r v_1 \le \cdots \le r v_j$. Then $D'_{\sigma_0}$ is also $N$-stable, since $N\vp = p\vp N$. Now for each
  $\sigma : F_0' \to E$ choose the unique $E$-subspace $D'_\sigma \subset D_\sigma$ that agrees with our choice of
  $D'_{\sigma_0}$ when $\sigma = \sigma_0$ and such that $D' := \bigoplus_{\sigma : F'_0 \to E} D'_\sigma$ is
  $\vp$-stable. Then $D'$ is stable under the actions of $F'_0 \otimes_{\qp} E$, $\vp$, $N$.  As in the proof of \cite[Prop.\
  4.5.2]{egh} (see also the proof of \cite[Prop.\ 5.1]{MR2359853}) we now compute that
  $t_N(D') = \frac{[E:\qp]}{[F_0:\qp]} \sum_{i = 1}^j v_i$ and that
  $t_H(D') \ge \frac{[E:\qp]}{[F:\qp]} \sum_{i = 1}^j \sum_{\tau : F \to E} h_{\tau,i}$.  By weak admissibility of $D$ we
  have
  \begin{equation}\label{eq:12}\tag{$*_j$}
    \sum_{i = 1}^j v_i \ge [F:F_0]^{-1} \sum_{i=1}^j \sum_{\tau : F \to E} h_{\tau,i},
  \end{equation}
  proving the first claim (with equality when $j = 0$ or $j = n$).
  
  Now suppose that equality holds in \eqref{eq:12} for some $1 \le j \le n-1$. If $v_j = v_{j+1}$ then monotonicity of the
  $h_{\tau,i}$ and $(*_{j+1})$ imply that equality holds in both $(*_{j-1})$, $(*_{j+1})$ and that $h_{\tau,j} = h_{\tau,j+1}$ for all $\tau$.
  Thus, by modifying $j$, we may assume without loss of
  generality that equality holds in \eqref{eq:12} and that $v_j < v_{j+1}$ (as $h_{\tau,1} < h_{\tau,n}$ for some $\tau$, by assumption).

  Let $D'$ be the sum of all generalized $\vp^{f'}$-eigenspaces in the $E$-vector space $D$ whose corresponding eigenvalues
  have valuation at most $r v_j$. As $v_j < v_{j+1}$ we see that $D'$ is a free $F'_0 \otimes_{\qp} E$-module of rank $j$,
  stable under the actions of $\vp$, $N$, and
  $\Gal(F'/F)$. Equality in~\eqref{eq:12} gives that $t_H(D') = t_N(D')$, so $\rho$ admits a $j$-dimensional
  subrepresentation.
\end{proof}

\begin{rk}
  The lemma can fail when $h_{\tau,1} = h_{\tau,n}$ for all $\tau$. For example, let $F/\qp$ be a quadratic extension and 
  $\chi : \Gamma_F \to \qpb\s$ a potentially unramified character that does not extend to $\Gamma_{\qp}$.  Then
  $\Ind_{\Gamma_F}^{\Gamma_{\qp}} \chi$ is irreducible and de Rham with all Hodge--Tate weights equal to 0 (since it is
  potentially unramified). In particular, $v_1 = v_2 = 0$. Concretely, via local class field theory, we can take
  $F = \Q_{p^2}$ and $\chi : \wh{\Q_{p^2}\s} \to \qpb\s$ tame and non-trivial on $\mu_{p+1}(\Q_{p^2})$.
\end{rk}

\subsection{On base change and descent for compact unitary groups}
\label{sec:base-change-descent}

The purpose of this section is to discuss base change and descent results for compact
unitary groups that go slightly beyond those in \cite{labesse}, namely allowing that
the unitary group is non-quasisplit at some finite places. The proofs will be provided
by Sug Woo Shin in Appendix~\ref{sec:appendix}.

Suppose that $\F/\F^+$ is a CM extension of number fields with $\F^+ \ne \Q$ and
$G$ a unitary group over $\F^+$ such that
\begin{enumerate}
\item $G_{/\F}$ is an inner form of $\GL_{nd}$;
\item $G(\F^+_u)$ is compact for any place $u \mid \infty$ of $\F^+$;
\item $G$ is quasi-split at all finite places that are inert in $\F/\F^+$.
\end{enumerate}
Let $c$ denote the complex conjugation of $\F/\F^+$.
Let $\Delta^+(G)$ denote the set of finite places of $\F^+$ where $G$ is not
quasi-split. This is a finite set of places that split or ramify in $\F$.
Let $\Delta(G)$ denote the set of places of $\F$ lying over a place of $\Delta^+(G)$.

\begin{prop}\label{prop:base-change}
  Suppose that $\pi$ is a \(cuspidal\) automorphic representation of $G(\A_{\F^+})$.  Then there
  exists a partition $n = n_1 + \cdots + n_r$ and discrete automorphic representations $\Pi_i$ of
  $\GL_{n_i d}(\A_{\F})$ satisfying $\Pi_i\dual \cong \Pi_i^c$ such that $\Pi := \Pi_1 \boxplus \cdots
  \boxplus \Pi_r$ is a weak base change of $\pi$. More precisely, at every finite split place $v = w
  w^c$ of $\F^+$ we have $|\LJ_{G(\F_w)}|(\Pi_w) \cong \pi_v$ as representations of $G(\F_w) \cong
  G(\F^+_v)$, and at infinity the compatibility is as in \cite[Cor.\ 5.3]{labesse}.
\end{prop}

\begin{prop}\label{prop:descent}
  Suppose that $\F/\F^+$ is unramified at all finite places and that $\Pi$ is a cuspidal automorphic
  representation of $\GL_{nd}(\A_{\F})$ such that $\Pi\dual \cong \Pi^c$, $\Pi_\infty$ is
  cohomological, and $\Pi_w$ is supercuspidal for all $w\in \Delta(G)$ \(in particular $|\LJ_{G(\F_w)}|(\Pi_w) \ne 0$\).
  Then there exists a \(cuspidal\)
  automorphic representation $\pi$ of $G(\A_{\F^+})$ such that at every finite
  split place $v = w w^c$ of $\F^+$ we have $|\LJ_{G(\F_w)}|(\Pi_w) \cong \pi_v$ as representations
  of $G(\F_w) \cong G(\F^+_v)$.
\end{prop}

\subsection{Supersingular representations of \texorpdfstring{$\GL_n(D)$}{GL\_n(D)}}
\label{sec:supers-repr-gl_nd}

We now prove the existence of supersingular (equivalently, supercuspidal) representations of $\GL_n(D)$ and $\PGL_n(D)$.

\begin{thm}\label{thm:main-PGLn}
  Suppose that $C$ is algebraically closed of characteristic $p$.
  For any smooth character $\zeta : F\s \to C\s$ there exists an irreducible admissible
  supercuspidal $C$-representation of $\GL_n(D)$ with central character $\zeta$.  In particular,
  there exists an irreducible admissible supercuspidal $C$-representation of $\PGL_n(D)$.
\end{thm}

\begin{coroll}\label{cor:main-PGLn}
  If $C$ is any field of characteristic $p$, then $\PGL_n(D)$ admits an irreducible admissible supercuspidal representation over
  $C$.
\end{coroll}

The proof uses Galois representations associated to automorphic representations on certain unitary groups.
We now make a few relevant definitions in preparation for the proof.

As in \S\ref{sec:lift-non-supers} we fix an isomorphism $\imath : \qpb \congto \C$.
Recall that if $\F/\F^+$ is a CM extension of number fields and $\Pi$ is a regular algebraic cuspidal
polarizable automorphic representation of $\GL_n(\A_{\F})$ (in the sense of \cite[\S2.1]{BLGGT}) we
have an associated semisimple potentially semistable $p$-adic Galois representation
$r_{p,\imath}(\Pi) : \Gamma_{\F} \to \GL_n(\qpb)$ that satisfies and is determined by local-global
compatibility with $\Pi$ at all finite places \cite[Thm.\ 2.1.1]{BLGGT}, \cite{blggt2}.

Suppose that $\F^+ \ne \Q$ and that $G$ is a unitary group over $\F^+$ as in
\S\ref{sec:base-change-descent}.  If $\pi$ is an automorphic representation of $G(\A_{\F^+})$, then
its weak base change $\Pi = \Pi_1 \boxplus \cdots \boxplus \Pi_r$ of
Proposition~\ref{prop:base-change} is regular algebraic and each $\Pi_i$ is polarizable.  By the
Moeglin--Waldspurger classification of the discrete spectrum and the previous paragraph it follows
that $\Pi$ has an associated semisimple potentially semistable $p$-adic Galois representation
$r_{p,\imath}(\pi) = r_{p,\imath}(\Pi) : \Gamma_{\F} \to \GL_{nd}(\qpb)$ that satisfies and is determined by
local-global compatibility with $\pi$ at all finite places of $\F$ that split over $\F^+$ and are not contained
in $\Delta(G)$. (We note that the Chebotarev density theorem shows that the set of Frobenius elements
at places $w$ of $\F$ that split over $\F^+$ is dense in $\Gamma_{\F}$.)
In particular, if $\Pi$ is not cuspidal, then $r_{p,\imath}(\pi)$ is reducible.

\begin{proof}[Proof of Theorem~\ref{thm:main-PGLn}]\

  \setcounter{step}{-1}
  \step{}\label{step:reduce-fpb}
  We show that it suffices to prove the theorem when $C = \fpb$.

  Given a smooth character $\zeta : F\s \to C\s$ we can define $\zeta' : F\s \to \fpb\s$ by
  extending $\zeta|_{\cO_F\s}$ (which is of finite order and hence takes values in $\fpb\s$)
  arbitrarily.  If Theorem~\ref{thm:main-PGLn} holds over $\fpb$, there exists an
  irreducible admissible supercuspidal $\fpb$-representation $\pi$ of $\GL_n(D)$ with central
  character $\zeta'$. Then by Step~\ref{step:sc-basechange} of the proof of Proposition~\ref{prop Cfinite}
  there exists an irreducible admissible supercuspidal $C$-representation $\pi'$ of $\GL_n(D)$
  with central character $\zeta'$. As $C$ is algebraically closed, a suitable unramified twist of
  $\pi'$ has central character $\zeta$.
  
  We will assume from now on that $C = \fpb$.

  \step{}\label{step:CMfield}
  We find a CM field $\F$ with maximal totally real subfield $\F^+ \ne \Q$ and a place $v
  \mid p$ of $\F^+$ such that
  \begin{enumerate}
  \item $\F/\F^+$ is unramified at all finite places;
  \item any place of $\F^+$ that divides $p$ splits in $\F$;
  \item $\F^+_v \cong F$;
  \end{enumerate}
  and a cyclic totally real extension $L^+/\F^+$ of degree $nd$ in which $v$ is inert.

  By Krasner's lemma we can find a totally real number field $H$, a place $u$ of $H$, and an
  isomorphism $H_u \congto F$. Now we apply \cite[Lemma 3.6]{henniart-gl3} and its proof to find
  finite totally real extensions $L^+/\F^+/H$ and a place $v$ of $\F^+$ above $u$ such that
  $L^+/\F^+$ is cyclic of degree $nd$, $\F^+_v = H_u$, and $v$ is inert in $L^+$. (We briefly recall
  the proof: pick a monic polynomial $Q$ of degree $nd$ over $F$ whose splitting field is the
  unramified extension of degree $nd$. Then let $L^+$ be the splitting field of a monic polynomial
  $P$ over $H$ that is $u$-adically very close to $Q$ and let $\F^+$ be the decomposition field of
  some place above $u$. We can use sign changes of $P$ at real places to ensure that $L^+$ is
  totally real.)

  Now pick any totally imaginary quadratic extension $\F/\F^+$ in which any place dividing $p$ splits.  By \cite[Lemma
  4.1.2]{cht} we can find a finite solvable Galois totally real extension $K^+/\F^+$ that is
  linearly disjoint from $L^+/\F^+$, such that $v$ splits in $K^+$, and such that for any prime $v'$
  of $\F^+$ that ramifies in $\F$ and any prime $w'$ of $K^+$ above $v'$ the extension
  $K^+_{w'}/\F^+_{v'}$ is isomorphic to the extension $\F_{v'}/\F^+_{v'}$. Then we can replace
  $\F/\F^+$ by $K^+\F/K^+$, $L^+$ by $K^+L^+$, and $v$ by any place of $K^+$ lying above $v$ to
  ensure that, without loss of generality, $\F/\F^+$ is unramified at all finite places.
  (In particular, we can always achieve $\F^+ \ne \Q$ in this way.)
  
  We let $w$ denote a place of $\F$ lying over $v$ and fix an isomorphism of topological fields $\F_w \congto F$.
  We let $L := L^+ \F$ and let $c$ denote the unique complex conjugation of $L$.

  \step{}\label{step:unitary-gp}
  Letting $v_1 \nmid p$ denote any place of $\F^+$ that is inert in $L^+$ and splits in $\F$,
  we now find a unitary group $G$ over $\F^+$ such that
  \begin{enumerate}
  \item $G_{/\F}$ is an inner form of $\GL_{nd}$;
  \item $G(\F^+_u)$ is compact for any place $u \mid \infty$ of $\F^+$;
  \item $G(\F_w) \cong \GL_n(D)$;
  \item $G$ is quasi-split at all finite places not contained in $\{v, v_1\}$.
  \end{enumerate}

  Let $G^*$ denote the unique quasi-split outer form of $\GL_{nd}$ over $\F^+$ that splits over $\F$.
  By \cite[\S2]{clozel-unitary} we can find an inner form $G$ of $G^*$ that satisfies all the above
  conditions. (If $nd$ is odd we do not need the auxiliary place $v_1$. If $nd$ is even we use $v_1$
  to ensure our local conditions can be globally realized.)

  The set $\Delta^+(G)$ (defined in \S\ref{sec:base-change-descent}) contains $v$ if $d > 1$ and is contained in $\{v, v_1\}$. 
  Any place of $\Delta(G)$ is inert in $L$ and splits over $\F^+$, and the set $\Delta_L(G)$ of places of $L$ lying above $\Delta(G)$
  is in canonical bijection with $\Delta(G)$.

  For any finite place $v' \not\in \Delta^+(G)$ of $\F^+$ that splits as $v' = w' w'^c$ in $\F$ 
  we obtain an isomorphism $\iota_{w'} : G(\F_{v'}^+) = G(\F_{w'}) \congto \GL_{nd}(\F_{w'})$ that is unique
  up to conjugacy. Moreover, $c\circ \iota_{w'}$ and $\iota_{w'^c}$ differ by an outer automorphism of
  $\GL_{nd}(\F_{w'^c})$. We also fix an isomorphism $\iota_w : G(\F_{v}^+) = G(\F_{w}) \congto \GL_{n}(D)$.
  (It is canonical, up to conjugacy, by condition (i).)

  \step{}\label{step:hecke-char}
  We find an algebraic Hecke character $\chi : \A_L\s/L\s \to \C\s$ with associated
  potentially crystalline $p$-adic Galois representation $\psi = r_{p,\imath}(\chi) : \Gamma_L \to
  \qpb\s$ (cf.\ \cite[Lemma 4.1.3]{cht}) such that
  \begin{enumerate}
  \item\label{item:conj-selfdual} $\psi \psi^c = \varepsilon^{-(nd-1)}$;
  \item\label{item:weil-irred} for any place $w' \in \Delta_L(G)$ the induced representation
    $\Ind_{W_{L_{w'}}}^{W_{\F_{w'}}} \chi_{w'}$ is irreducible;
  \item\label{item:regular} the representation $r := \Ind_{\Gamma_L}^{\Gamma_{\F}} \psi$ has regular Hodge--Tate
    weights, i.e., for each $\kappa' : \F \to \qpb$ the $nd$ integers
    $\HT_{\kappa'}(\Ind_{\Gamma_L}^{\Gamma_{\F}} \psi)$ are pairwise distinct;
  \item\label{item:irred} the restriction $\o r|_{\Gamma_{\F_w}}$ to $\Gamma_{\F_w}$ of the reduction $\o r \cong
    \Ind_{\Gamma_L}^{\Gamma_{\F}} \o\psi$ is irreducible.
  \end{enumerate}

  We first introduce some notation. Let $\Delta_p$ denote the places $w'$ of $L$ that divide $p$.  
  Note that, by construction, any place $w'\in \Delta_L(G) \cup \Delta_p$ splits over $L^+$, i.e.\ $w' \ne w'^c$.
  Let $S_K := \Hom\cts(K,\qpb)$ for any topological field $K$ of characteristic zero and $S_k := \Hom(k,\fpb)$ for any field $k$
  of characteristic $p$.

  Our strategy is to carefully choose continuous characters $\theta_{w'} : \Gamma_{L_{w'}} \to \qpb\s$ for any
  $w' \in \Delta_L(G) \cup \Delta_p$ that satisfy $(\theta_{w'} \theta^c_{w'^c})|_{I_{L_{w'}}} = \ve^{-(nd-1)}|_{I_{L_{w'}}}$
  and are potentially crystalline when $w' \in \Delta_p$. We then deduce by \cite[Lemma A.2.5(1)]{BLGGT} that there exists a
  character $\psi : \Gamma_L \to \qpb\s$ such that $\psi \psi^c = \varepsilon^{-(nd-1)}$ and
  $\psi|_{I_{L_{w'}}} = \theta_{w'}|_{I_{L_{w'}}}$ for all $w' \in \Delta_L(G) \cup\Delta_p$. In particular, $\psi$ is
  potentially crystalline, and we let $\chi$ be the associated algebraic Hecke character. It follows that
  condition~\ref{item:conj-selfdual} holds.

  For any $w' \in \Delta_L(G)$ we can choose a smooth character $\zeta_{w'} : \Gamma_{L_{w'}}\ab \cong \wh{L_{w'}\s}
  \to \qpb\s$ such that the $\Gal(L_{w'}/\F_{w'})$-conjugates of $\zeta_{w'}|_{\O_{L_{w'}}\s}$ are
  pairwise distinct. (For example, we can take a faithful character of $k_{L_{w'}}\s$
  and inflate it to $\O_{L_{w'}}\s$.) We may assume without loss of generality that
  $\zeta_{w'}\zeta_{w'^c}^c = 1$.

  Now suppose that $w' \in \Delta_p$. Suppose that we are given any integers $\lambda_\kappa$ ($\kappa \in S_{L}$) satisfying
  $\lambda_\kappa + \lambda_{\kappa c} = nd-1$ for all $\kappa \in S_L$.  Let $\theta\crs_{w'} : \Gamma_{L_{w'}} \to \qpb\s$
  be any crystalline character with $\HT_\kappa(\theta\crs_{w'}) = \lambda_\kappa$ for all
  $\kappa \in S_{L_{w'}} \subset S_L$.  Without loss of generality, by our constraint on the $\lambda_\kappa$, we may assume
  that $\theta\crs_{w'}(\theta\crs_{w'^c})^c = \ve^{-(nd-1)}$.

  For $w' \in \Delta_L(G) \cup \Delta_p$ define
  \begin{equation*}
    \theta_{w'} :=
    \begin{cases}
      \zeta_{w'} & \text{if $w' \in \Delta_L(G) - \Delta_p$;} \\
      \theta\crs_{w'}\zeta_{w'} & \text{if $w' \in \Delta_L(G) \cap \Delta_p$;} \\
      \theta\crs_{w'} & \text{if $w' \in \Delta_p - \Delta_L(G)$.} \\
    \end{cases}
  \end{equation*}

  This completes the construction of a potentially crystalline character $\psi$ and its associated algebraic Hecke character $\chi$.  By
  construction, for any $w' \in \Delta_L(G)$ the character $\imath \zeta_{w'}|_{I_{L_{w'}}}$ corresponds to
  $\chi_{w'}|_{\O_{L_{w'}}\s}$ under the local Artin map.  Therefore, since the $\Gal(L_{w'}/\F_{w'})$-conjugates of
  $\zeta_{w'}|_{\O_{L_{w'}}\s}$ are pairwise distinct, we deduce that condition~\ref{item:weil-irred} holds.

  Finally, we will choose the integers $\lambda_\kappa$ ($\kappa \in S_{L}$) so that conditions~\ref{item:regular} and
  \ref{item:irred} hold.  Note that condition~\ref{item:regular} is equivalent to the condition
  \begin{enumerate}
  \item[\mylabel{item:regp}{(iii$'$)}] for any $\kappa' \in S_{\F}$ the $nd$ integers $\{ \lambda_\kappa
    : \kappa \in S_L, \kappa|_{\F} = \kappa' \}$ are pairwise distinct.
  \end{enumerate}

  First choose the $\lambda_\kappa$ for those $\kappa \in S_L$ that do not induce either of the
  places $w$, $w^c$ on $L$ so that condition \ref{item:regp} holds for any $\kappa' \in
  S_{\F}$ not inducing either of the places $w$, $w^c$ on $\F$. It remains to choose the
  $\lambda_\kappa$ for those $\kappa$ that induce the place $w$ on $L$ (since the remaining
  $\lambda_\kappa$ are determined by the condition $\lambda_\kappa + \lambda_{\kappa c} = nd-1$ for
  all $\kappa$), i.e.\ for $\kappa \in S_{L_{w}}$.

  To choose the $\lambda_{\kappa}$ for $\kappa \in S_{L_{w}}$, we note that $\o r|_{\Gamma_{\F_w}} \cong
  \Ind_{\Gamma_{L_w}}^{\Gamma_{\F_w}} (\o\psi|_{\Gamma_{L_w}})$ is irreducible if and only if the
  $\Gal(L_{w}/\F_{w})$-conjugates of $\o\psi|_{\Gamma_{L_w}}$ are pairwise distinct, or equivalently if
  the characters $\o\psi|_{I_{L_w}}^{q^i}$ ($0 \le i \le nd-1$) are pairwise distinct.  (Recall that
  $q = \#k_F$.)  We have $\o\psi|_{I_{L_{w}}} \cong \o{\theta\crs_w \zeta_w}|_{I_{L_{w}}}$ if $d > 1$ or $\o\psi|_{I_{L_{w}}} \cong \o{\theta\crs_w}|_{I_{L_{w}}}$ otherwise.  By
  \cite[Cor.\ 7.1.2]{GHS}, noting our opposite conventions concerning Hodge--Tate weights, we have $\o{\theta\crs_w}|_{I_{L_{w}}} = \prod_{\sigma \in S_{k_{L_w}}}
  \omega_{\sigma}^{-b_{\sigma}}$, where $\omega_\sigma$ corresponds to the character 
  $\O_{L_w}\s \onto k_{L_w}\s \xrightarrow{\sigma} \fpb\s$ under local class field theory and
  $b_{\sigma} := \sum_{\kappa \in S_{L_{w}} : \o\kappa =
    \sigma} \lambda_{\kappa}$.  Fix any $\sigma \in S_{k_{L_w}}$ and $s \in \Z$. Then we can choose the $\lambda_\kappa$ for
    $\kappa \in S_{L_{w}}$ so that $\o\psi|_{I_{L_{w}}} = \omega_\sigma^s$. By taking $s$ so that the $\omega_\sigma^{sq^i}$
    ($i = 0, \dots, nd-1$) are pairwise distinct (taking, for example, $s=1$), condition~\ref{item:irred} holds.  Finally, we can ensure that condition~\ref{item:regp}
  holds for all $\kappa' \in S_{\F_w}$ while keeping $\o r|_{I_{\F_w}}$ unchanged by varying the
  $\lambda_{\kappa}$ (for $\kappa \in S_{L_{w}}$) modulo $q^{nd}-1$. This completes Step~\ref{step:hecke-char}.

  \step{}\label{step:autom-rep}
  Using automorphic induction and descent we define an automorphic representation
  $\pi'$ of $G(\A_{\F^+})$ with associated Galois representation $r = \Ind_{\Gamma_L}^{\Gamma_{\F}} \psi$.

  Let $\Pi''$ denote the automorphic induction of $\chi$ with respect to the cyclic extension $L/\F$.
  It is an automorphic representation of $\GL_{nd}(\A_{\F})$ that is parabolically induced from
  a cuspidal representation. (For the functoriality of automorphic induction in cyclic extensions we
  refer to \cite{henniart-induction}, which shows in particular that it is compatible with local
  automorphic induction at all places. Note the results of \cite{henniart-induction} apply to
  unitary representations, but by twisting they continue to hold for twists of unitary
  representations.)

  We claim that $\Pi''$ is cuspidal. This follows from \cite{henniart-induction}, Theorems 2, 3, and
  Proposition 2.5, provided that the Hecke characters $\{\chi^\sigma : \sigma \in \Gal(L/\F)\}$
  are pairwise distinct. Equivalently, the Galois characters $\{\psi^\sigma : \sigma \in \Gal(L/\F)\}$ 
  are pairwise distinct, which in turn is equivalent to the condition that
  $\Ind_{\Gamma_L}^{\Gamma_{\F}} \psi$ is irreducible. This is a consequence of condition \ref{item:irred} in Step~\ref{step:hecke-char},
  so $\Pi''$ is cuspidal.
  
  Let $\Pi' := \Pi'' \otimes |\det|_{\F}^{(nd-1)/2}$. By condition~\ref{item:conj-selfdual} in Step~\ref{step:hecke-char} we have $\chi \chi^c = |\cdot|_L^{-(nd-1)}$, hence
  $(\Pi')\dual \cong \Pi'^c$.  On the other hand, $\Pi'_\infty$ is cohomological by \cite[Lemma
  3.14]{MR1044819}, as it is regular by condition~\ref{item:regp} in Step~\ref{step:hecke-char}. It follows that $\Pi'$ is regular algebraic and polarizable in the sense of
  \cite[\S2.1]{BLGGT}, so we have an associated Galois representation $r_{p,\imath}(\Pi')$. By local-global
  compatibility at unramified places and Chebotarev we deduce that $r_{p,\imath}(\Pi') \cong \Ind_{\Gamma_L}^{\Gamma_{\F}}
  \psi$.

  For $w' \in \Delta(G)$ the local factor $\Pi'_{w'}$ is supercuspidal, as $\rec_{\F_{w'}}(\Pi'_{w'})$ is
  irreducible by condition~\ref{item:weil-irred} in Step~\ref{step:hecke-char}. It follows from what we recalled in
  \S\ref{sec:jacq-langl-corr} that $|\LJ_{G(\F_{w'})}|(\Pi'_{w'}) \ne 0$.

  By Proposition~\ref{prop:descent} we deduce that $\Pi'$ descends to a (cuspidal) automorphic representation $\pi'$ of
  $G(\A_{\F^+})$, such that for all finite places $v' \not\in \Delta^+(G)$ of $\F^+$ that split
  as $v' = w' w'^c$ in $\F$ we have $\pi'_{v'} \cong \Pi'_{w'}$ as representations of $G(\F^+_{v'})
  \cong \GL_{nd}(\F_{w'})$. We deduce that $r_{p,\imath}(\pi') \cong \Ind_{\Gamma_L}^{\Gamma_{\F}} \psi$.

  \step{}\label{step:local-candidate}
  We use the automorphic representation $\pi'$ to define an irreducible admissible
  $\fpb$-representation $\sigma$ of $G(\F^+_v) \cong \GL_n(D)$.

  Fix a maximal compact open subgroup $K_p$ of $\prod_{v' \mid p} G(\F^+_{v'})$.  If $U$ is any
  compact open subgroup of $K_p G(\A_{\F^+}^{\infty,p})$ and $\cW$ is any $\zpb[K_p]$-module, we let $S(U,\cW)$
  be the $\zpb$-module of functions $f : G(\F^+)\backslash G(\A_{\F^+}^\infty) \to \cW$ such that $f(gu) =
  u_p^{-1} f(g)$ for all $g \in G(\A_{\F^+}^\infty)$ and $u \in U$ (where $u_p$ denotes the projection of $u$ to $K_p$).

  Using the compactness of $G$ at infinity, we see as in \cite[Lemma 7.1.6]{egh} that there exists a $\qp$-algebraic representation
  $\cW\subalg$ of $\prod_{v' \mid p} G(\F^+_{v'})$ over $\qpb$ such that $\ilim_{U} S(U,\cW\subalg)$ contains
  $\imath^{-1} \pi'^\infty$ as $G(\A_{\F^+}^\infty)$-representation.
  Choose a $K_p$-invariant $\zpb$-lattice $\cW\subalg^\circ$ in $\cW\subalg$ and let $\o{\cW\subalg} := \cW\subalg^\circ \otimes_{\zpb} \fpb$. 

  Pick a compact open subgroup $U = \prod_{v'\nmid \infty} U_{v'}$ of $G(\A_{\F^+}^{\infty})$ such that
  \begin{enumerate}
  \item $(\pi'^\infty)^U \ne 0$;
  \item there exists a place $v'\nmid p\infty$ of $\F^+$ such that $U_{v'}$ contains no element of finite order other than the
    identity;
  \item the group $\prod_{v'\mid p} U_{v'}$ is contained in $K_p$ and acts trivially on $\o{\cW\subalg}$.  
  \end{enumerate}
  Note that condition (ii) implies that for any compact open subgroup $U' = U'_p \prod_{v'\nmid p\infty} U_{v'}$
  with $U'_p \le K_p$ we have $S(U',\cW) \cong \cW^{\oplus s}$ as $\zpb$-modules for some $s \ge 1$ depending only on $U'$.
  In particular,
  \begin{equation}
    S(U',\cW) \otimes_{\zpb} R \congto S(U',\cW \otimes_{\zpb} R)\label{eq:13}
  \end{equation}
  for any $\zpb$-algebra $R$ (see e.g.\ \cite[\S7.1.2]{egh}).
  We will apply this with $R = \qpb$ and $R = \fpb$.

  Let $\cP$ denote the set of places $w' \nmid p$ of $\F$ that split over a place $v'$ of $\F^+$ not contained in $\Delta^+(G)$, 
  and are such that $U_{v'}$ is a maximal compact subgroup of $G(\F^+_{v'})$. For each such $w'$ we conjugate the
  isomorphism $\iota_{w'}$ of Step~\ref{step:unitary-gp} so that $\iota_{w'}(U_{v'}) = \GL_{nd}(\O_{\F_{w'}})$. 
  Note that the set $\cP$ has finite complement in the set of places of $\F$ that split over $\F^+$. Let
  $\T^\cP$ denote the commutative polynomial $\zpb$-algebra in the variables $T_{w'}^{(i)}$
  for $w' \in \cP$ and $0 \le i \le nd$, acting on any $S(U,\cW)$ as double coset operators as in \cite[\S7.1.2]{egh}.
  Note that the ring $\T^\cP$ acts by scalars on $(\imath^{-1} \pi'^\infty)^U$ 
  inside $S(U,\cW\subalg)$ and stabilizes the $\zpb$-lattice $S(U,\cW\subalg^\circ)$. Therefore there exists a unique maximal ideal $\m$ of $\T^{\cP}$ with residue field $\fpb$ 
  such that $(\imath^{-1} \pi'^\infty)^U \subset S(U,\cW\subalg)_\m$.

  Applying~\eqref{eq:13} and localizing at $\m$ we obtain that $S(U,\o{\cW\subalg})_{\m} \ne 0$.
  Then 
  $$S(U,\fpb)_{\m} \otimes_{\fpb} \o{\cW\subalg} \cong S(U,\o{\cW\subalg})_{\m} \ne 0,$$
  where the isomorphism uses condition (iii) on $U$.
  Writing $U^v := \prod_{v'\ne v} U_{v'}$ and $S(U^v, \fpb) := \ilim_{U_v} S(U^v U_v,
  \fpb)$, we have $S(U^v, \fpb)_{\m} \ne 0$. This is a non-zero admissible smooth representation of $G(\F_v^+) \cong \GL_n(D)$,
  using the isomorphism $\iota_w$ of Step~\ref{step:unitary-gp}.
  Let $\sigma$ be an irreducible (admissible) $\GL_n(D)$-subrepresentation of 
  $S(U^v, \fpb)_{\m}$, which exists by the proof of Lemma 9.9 in \cite{bib:herzig-classification}
  or \cite[Lemma 7.10]{HV2}. 

  \step{}\label{step:ss}
  We show that $\sigma$ is supersingular, or equivalently, supercuspidal.

  By \cite[Thm.\ 3]{OV} it suffices to show that the $\HH_{\fpb}$-module $\sigma^{I(1)}$ is
  supersingular, where $I(1)$ denotes the pro-$p$ Iwahori subgroup of $\GL_n(D) \cong G(\F_v^+)$ defined in
  \S\ref{sec:lift-non-supers}.
  In fact, we will even show that $(S(U^v, \fpb)_{\m})^{I(1)} \cong S(U^v\cdot
  I(1),\fpb)_{\m}$ is supersingular. Assume by contradiction that this is false, so one of the operators $\mathscr Z_j$ for $1
  \le j \le n-1$ has a nonzero eigenvalue $\lambda_j$ on $S(U^v\cdot I(1),\fpb)_{\m}$. 

  Again from~\eqref{eq:13} we know that $S(U^v\cdot I(1),\zpb)\otimes_{\zpb} R \cong S(U^v\cdot I(1),R)$ for $R = \qpb$ and
  $R = \fpb$.  By applying \cite[Lemma 4.5.1]{egh} (a version of the Deligne--Serre lemma) with $A = \T^{\cP}[\mathscr Z_j]$,
  $M = S(U^v\cdot I(1),\zpb)$, $\n$ the maximal ideal of $A$ generated by $\m$ and $\mathscr Z_j-\lambda_j$, we deduce that
  there exists a 
  homomorphism $\theta : A \to \zpb$ such that the $\theta$-eigenspace of $S(U^v\cdot I(1),\qpb)$ is
  non-zero, $\ker(\o\theta|_{\T^\cP}) = \m$, and $\theta(\mathscr Z_j) \in \zpb^\times$.
  By \cite[Lemma 7.1.6]{egh}, there exists an automorphic representation $\pi$ of $G(\A_{\F^+})$ satisfying
  \begin{enumerate}
  \item $(\imath^{-1}\pi^\infty)^{U^v\cdot I(1)}$ has a non-trivial $\theta$-eigenspace;
  \item $\pi_\infty$ is trivial.
  \end{enumerate}
  It follows from (i) that $\imath^{-1}\pi_v^{I(1)} \ne 0$ is an integral $\HH_{\qpb}$-module whose reduction is non-supersingular.
  (A priori, we get that $(\imath^{-1}\pi_v^{I(1)})^{\oplus s}$ is integral for some $s \ge 1$, but then we can project to any factor.
  Note that any finitely generated submodule of a finite free $\zpb$-module is free.)

  By local-global compatibility and \cite[Cor.\ 3.1.2]{cht},
  for any $w' \in \cP$ the characteristic polynomial of $\o r(\Frob_{w'})$ equals $\sum_{i = 0}^{nd} (-1)^i
  (\mathbf{N}w')^{i(i-1)/2}T_{w'}^{(i)}X^{nd-i}$ modulo $\m$, where
  $\Frob_{w'}$ denotes a geometric Frobenius element at $w'$. The same is true for $\o r_{p,\imath}(\pi)$, 
  as $\ker(\o\theta|_{\T^\cP}) = \m$, and hence we deduce by the Chebotarev density theorem that $\o{r_{p,\imath}(\pi)} \cong \o r$.

  By Proposition~\ref{prop:base-change} we obtain an automorphic representation $\Pi$ of $\GL_{nd}(\A_{\F})$ with
  associated Galois representation $r_{p,\imath}(\Pi)$ lifting $\o r$ such that $|\LJ_{G(\F_{w'})}|(\Pi_{w'})
  \cong \pi_{v'}$ for all finite places $v'$ of $\F^+$ that split as $v' = w' w'^c$ in $\F$. As $\o
  r$ is irreducible by construction we know that $\Pi$ is cuspidal. In particular, $\Pi_{w'}$ is
  essentially unitarizable and generic for each finite place $w'$ of $\F$. Let $v_1 \le \cdots \le
  v_{nd}$ denote the valuations of the eigenvalues of a geometric Frobenius on
  $\imath^{-1}(\rec_F(\Pi_w))$. From Proposition~\ref{prop:non-supersingular-reduction-and-frob-evals}
  (applied to $\Pi_w$) we deduce that there exists $1 \le j \le
  n-1$ such that  
  \begin{equation}\label{eq:5}
    \sum_{i=1}^{jd} v_i = -\frac{d^2 j(n-j)}2 \val(q).
  \end{equation}
  Note that the infinitesimal character of $\Pi$ is the same as that of the
  trivial representation. By \cite[Thm.\ 2.1.1]{BLGGT} we deduce that
  $\HT_\tau(r_{p,\imath}(\Pi)|_{\Gamma_{\F_w}}) = \{0,1,\dots,nd-1\}$ for all $\tau : \F_w \to \qpb$
  and that $\imath \WD(r_{p,\imath}(\Pi)|_{\Gamma_{\F_w}})\Fss \cong \rec_F(\Pi_w \otimes |\det|_F^{(1-nd)/2})$.
  Together with~\eqref{eq:5} it follows that 
  \begin{equation*}
    \sum_{i=1}^{jd} v'_i = -\frac{d^2 j(n-j)}2 \val(q) + jd\val(q^{(nd-1)/2}) = \binom{jd}2 \val(q),
  \end{equation*}
  where $v'_1 \le \cdots \le v'_{nd}$ denote the valuations of the
  eigenvalues of a geometric Frobenius on $\WD(r_{p,\imath}(\Pi)|_{\Gamma_{\F_w}})$. By
  Lemma~\ref{lm:hodge-newton-reducibility}, noting that $\val(q) = [F_0:\qp]$, it follows that
  $r_{p,\imath}(\Pi)|_{\Gamma_{\F_w}}$ is reducible, which contradicts that its reduction $\o
  r|_{\Gamma_{\F_w}}$ is irreducible by Step~\ref{step:hecke-char}.

  \step{}\label{step:central-char}
  We fix the central character. 

  Suppose we are given a smooth character $\zeta : F\s \to \fpb\s$.  As in Step~\ref{step:reduce-fpb} it is enough to
  construct an irreducible admissible supercuspidal representation such that $\cO_F\s$ acts via
  $\zeta|_{\cO_F\s}$. 

  Note that $\sigma$ has a central character $\chi_\sigma$, as it is irreducible and admissible.  We
  claim that $\chi_\sigma|_{\O_F\s} = \det (\o r|_{I_{\F_w}}) \cdot \o\ve^{nd(nd-1)/2}$ under the
  local Artin map.  The central character of the $\GL_n(D)$-representation $\imath^{-1}\pi_v$ in
  Step~\ref{step:ss} lifts $\chi_\sigma$ and is equal to the central character of $\imath^{-1}\Pi_w$.  (This
  equality follows from the definition of $\LJ$ in \cite[\S2.7]{badu-jlc}, noting that $\Pi_w$ is
  generic and hence fully induced from an essentially square-integrable representation.) By
  local-global compatibility at $p$ (cf.\ Step~\ref{step:ss}) the latter character equals $\WD(\det
  r_{p,\imath}(\Pi)|_{\Gamma_{\F_w}})|_{I_{\F_w}}$ on $\O_F\s$, under the local Artin map.  As $
  r_{p,\imath}(\Pi)|_{\Gamma_{\F_w}}$ has parallel Hodge--Tate weights $0,1,\dots,nd-1$, we have $\det
  r_{p,\imath}(\Pi)|_{I_{\F_w}} = \ve^{-nd(nd-1)/2} \cdot \WD(\det
  r_{p,\imath}(\Pi)|_{\Gamma_{\F_w}})|_{I_{\F_w}}$ and hence deduce the claim.

  It thus suffices to show that in Step~\ref{step:hecke-char} above we can choose $r$ such that $\det (\o
  r|_{I_{\F_w}})$ is any prescribed character that is extendable to $\Gamma_{\F_w}$.  Let us fix any
  $\o\kappa \in S_{k_{L_w}}$ and write $\o\psi|_{I_{L_w}}=\omega_{\o\kappa}^s$ for some $s \in \Z$. Then
  the condition that the $\o\psi|_{I_{L_w}}^{q^i}$ ($i = 0, 1,\dots, nd-1$) are pairwise distinct
  means:
  \begin{equation}
    s \not\equiv 0 \pmod {\textstyle\frac{q^{nd}-1}{q^\ell-1}} \quad \forall \ell \mid nd, \ \ell < nd.\label{eq:2}
  \end{equation}
  On the other hand, $\det (\o r|_{I_{\F_w}}) = \prod_{i=0}^{nd-1} \o\psi|_{I_{L_w}}^{q^i} =
  \omega_{\o\kappa'}^s$, where $\o\kappa' \in S_{k_{\F_w}}$ is the restriction of $\o\kappa$ to
  $k_{\F_w}$. As any character $\Gamma_{\F_w} \to \fpb\s$ restricts to a power of $\omega_{\o\kappa'}$ on
  inertia, we can prescribe $\det (\o r|_{I_{\F_w}})$ if and only if we can choose $s$ in any
  residue class modulo $q-1$.  Since $\frac{q^{nd}-1}{q^\ell-1} \ge q+1$ for any $\ell \mid nd$, $\ell <
  nd$, it follows that we can pick any $s$ in the interval $[1,q-1]$, completing the proof.
\end{proof}

\begin{proof}[Proof of Corollary~\ref{cor:main-PGLn}]
  Going back to Step~\ref{step:local-candidate} of the proof of Theorem~\ref{thm:main-PGLn}, it is clear that
  the representation $S(U^v, \fpb)_{\m} \ne 0$ is defined over a finite field (as $\o r$ is),
  and hence so is its irreducible subrepresentation $\sigma$. This proves the corollary when $C$ is a sufficiently
  large finite field of characteristic $p$. 
  We conclude by Proposition~\ref{prop Cfinite}.
\end{proof}

\appendix
\address[S.\ W.\ Shin]{Department of Mathematics\\
University of California, Berkeley\\
901 Evans Hall, Berkeley, CA 94720, USA
/ Korea Institute for Advanced Study\\
 Dongdaemun-gu, Seoul 130-722, Republic of Korea}
\email{sug.woo.shin@berkeley.edu}
\makeatletter
\tracingmacros=1
\def\@tocwrite#1#2{}%
\section{Base change}\label{sec:appendix}
\let\@tocwrite=\original@@tocwrite 

\@tocwrite{section}{Base change (by Sug Woo Shin)}
\makeatother

\begin{center}
Sug Woo Shin\footnote{Supported by NSF grant DMS-1802039.}
\end{center}
\markboth{Sug Woo Shin}{Base change}
\medskip

In this appendix we will prove Propositions~\ref{prop:base-change} and \ref{prop:descent}.

We need a character identity for the Jacquet--Langlands correspondence. We fix compatible Haar measures on $\GL_{nd}(F)$ and $\GL_n(D)$ in the sense of \cite[p.631]{Kottwitz1988}.
We say that $f\in C^\infty_c(\GL_n(D))$ and $f^*\in C^\infty_c(\GL_{nd}(F))$ are associated, or that $f^*$ is a transfer of $f$, if the orbital integral identity $O_\delta(f)=O_{\delta^*}(f^*)$ holds for every regular semisimple elements $\delta\in \GL_n(D)$ and $\delta^*\in \GL_{nd}(F)$ with the same characteristic polynomial. (We use the same Haar measures on the centralizers of $\delta$ and $\delta^*$ in $\GL_n(D)$ and $\GL_{nd}(F)$, respectively, to compute the orbital integrals.) A well-known fact, proven in \cite{DKV}, is that every $f\in C^\infty_c(\GL_n(D))$ admits a transfer in $C^\infty_c(\GL_{nd}(F))$. (This is a special case of the Langlands--Shelstad transfer.) Let $e(\GL_n(D))\in \{\pm 1\}$ denote the Kottwitz sign \cite{Kottwitz1983}. Explicitly $e(\GL_n(D))=(-1)^{nd-n}$.

\begin{prop}\label{prop:JL-char-id}
  Let $\pi^*$ be an irreducible unitarizable representation of $\GL_{nd}(F)$. For every associated pair $f\in C^\infty_c(\GL_n(D))$ and $f^*\in C^\infty_c(\GL_{nd}(F))$, we have
  $$\tr \pi^*(f^*) =e(\GL_n(D))\cdot \tr \left(|\LJ_{\GL_n(D)}|(\pi^*)\right)(f).$$
\end{prop}

\begin{proof}
  This follows from \cite[Prop 3.3]{badulescu07} and the Weyl integration formula \cite[A.3.f]{DKV} for $\GL_n(D)$ and $\GL_{nd}(F)$.
\end{proof}

We assume that the CM extension $\F/\F^+$ and the unitary group $G$ over $\F^+$ are as in Section~\ref{sec:base-change-descent}.

Write $G^*$ for a quasi-split inner twist of $G$ over $\F^+$ (with an isomorphism between $G^*$ and $G$ over an algebraic closure of $\F^+$). By convention,
every trace considered on $p$-adic or adelic points of $G^*$ over $\F$ (as opposed to over $\F^+$) will mean the \emph{twisted} trace relative to the action of $\Gal(\F/\F^+)$ on $\Res_{\F/\F^+} G^*$ (with the Whittaker normalization), unless specified otherwise.

\begin{proof}[Proof of Proposition~\ref{prop:base-change}]
  This proposition is implied by \cite[Cor.\ 5.3]{labesse} except possibly the relation
  $|\LJ_{G(\F_w)}|(\Pi_w) \cong \pi_v$.\footnote{In fact this assertion is implicit in \cite[Cor.\ 5.3]{labesse} where it reads
  ``Aux places non ramifi\'ees ou d\'ecompos\'ees la correspondance $\sigma_v\mapsto \pi_v$ est donn\'ee par le changement de base local.''
  However when $v=w w^c$ the author introduced the notion of local base change (\S4.10 of \emph{op.\ cit.}) only when $\mathbf{U}$ is a general linear group at $v$ (in his notation).
  We need the case when $\mathbf{U}$ is a nontrivial inner form of a general linear group at $v$.}
  We elaborate on this point. Thus we assume $v=ww^c$ as in the proposition. We will omit the subscript for $|\LJ|$ when there is little danger of confusion.
  
  Let $S$ be a finite set of places of $\F^+$ containing all infinite places as well as all finite place where either $\pi$ or $G$ is ramified. Denote by $S_{\mathrm{f}}$ the subset of finite places in $S$. In particular $S_{\mathrm{f}}\supset \Delta^+(G)$.
  For an irreducible admissible representation $\sigma$ of $G(\A_{\F^+})$ unramified outside $S$, we write $\BC(\sigma^S)=\Pi^S$ to mean that the local unramified base change of $\sigma_u$ is $\Pi_u$ at all places $u\notin S$. (The unramified base change is defined via the Satake transform.) Using the Langlands parametrization at archimedean places, we write $\BC(\sigma_\infty)=\Pi_\infty$ to mean that the local base change of $\sigma_\infty$ is $\Pi_\infty$.
  
   For each finite place $u$ and $f_u\in C^\infty_c(G(\F^+_u))$ let $f^*_u\in C^\infty_c(G^*(\F^+_u))$ denote a transfer. There exists $\phi_u\in C^\infty_c(G^*(\F \otimes_{\F^+} \F^+_u))$ whose base change transfer is $f^*_u$ by \cite[Lem.\ 4.1]{labesse}. 
   Write $f_{S_{\mathrm{f}}} :=\prod_{u\in S_{\mathrm{f}}} f_u$  and $\phi_{S_{\mathrm{f}}} :=\prod_{u\in S_{\mathrm{f}}} \phi_u$.
   
  Let $\Pi_v:=\Pi_w\otimes \Pi_{w^c}$ be the $v$-component of $\Pi$, which is a representation of $G^*(\F \otimes_{\F^+} \F^+_v)$.
  Let $\pi^*_v:=\Pi_w$ via the isomorphism $G^*(\F^+_v)\cong G^*(\F_w)$. Then we have the following character identities, where $\tr \Pi_v(\phi_v)$ means the twisted trace by abuse of notation:
  \begin{equation}\label{BC-LJ-identity}
 \tr \Pi_v(\phi_v)= \tr \pi^*_v(f^*_v)= e(G(\F^+_v))\cdot \tr\left( |\LJ|(\pi^*_v)\right)(f_v).
      \end{equation}
      The first equality holds by the same computation as for \cite[Prop.\ 4.13.2 (a)]{Roga}. The second equality is Proposition \ref{prop:JL-char-id}.
        On the other hand, the trace formula argument of \cite[Thm.\ 5.1]{labesse} shows
  \begin{equation}\label{trace-formula-1}
 \sum_{\sigma} m(\sigma) \tr \sigma_{S_{\mathrm{f}}}(f_{S_{\mathrm{f}}}) =  c\cdot \tr \Pi_{S_{\mathrm{f}}}(\phi_{S_{\mathrm{f}}}),
    \end{equation}
  with a constant $c$ and the automorphic multiplicity $m(\sigma)\in \Z_{\ge 0}$, where
  the sum runs over $\sigma$ such that $\BC(\sigma^S)=\Pi^S$ and $\BC(\sigma_\infty)=\Pi_\infty$. Again the trace on the right-hand side is the twisted trace.
  Since \eqref{trace-formula-1} holds for each $f^\infty=\prod_{u\nmid \infty} f_u$ (and $f^*_u$ and $\phi_u$ constructed from $f_u$ at each $u$ as above), we choose 
  $f_u$ to be the characteristic function on a sufficiently small compact open subgroup of $G(\F^+_u)$ at $u\in S_{\mathrm{f}}\backslash \{v\}$. Then $\tr \sigma_u(f_u)\ge 0$, so we obtain
  $$\sum_{\sigma} n(\Pi_v,\sigma)\tr \sigma_v(f_v)=c'\cdot \tr \Pi_v(\phi_v), \quad\mbox{with}~~n(\Pi_v,\sigma)\ge 0,$$
  where $c'$ is a new constant and the sum runs over $\sigma$ such that $\BC(\sigma^S)=\Pi^S$, $\BC(\sigma_\infty)=\Pi_\infty$, and $\tr \sigma_u(f_u)\neq 0$ at every $u\in S_{\mathrm{f}}\backslash \{v\}$. Notice that $\sigma=\pi$ contributes to the sum with $n(\Pi_v,\pi)>0$, by choice of $f_u$ at $u\in S_{\mathrm{f}}\backslash \{v\}$. By choosing a suitable $f_v$ we deduce that $c' \ne 0$.
   Substituting \eqref{BC-LJ-identity} we obtain
  $$\sum_{\sigma} n(\Pi_v,\sigma)\tr \sigma_v(f_v)=c'\cdot e(G(\F^+_v))\cdot \tr \left( |\LJ|(\pi^*_v)\right) (f_v),$$
  with the sum running over the same set of $\sigma$. Since the sum is not identically zero, $|\LJ|(\pi^*_v)$ is irreducible (rather than $0$). By linear independence of characters of $G(\F^+_v)$, we deduce that the coefficients on the left-hand side are zero unless $\sigma_v\cong |\LJ|(\pi^*_v)$. Since $n(\Pi_v,\pi)>0$, we must have $\pi_v\cong |\LJ|(\pi^*_v)$, noting that no cancellation takes place in the sum as the coefficients are non-negative.
\end{proof}

\begin{proof}[Proof of Proposition~\ref{prop:descent}]
  The proposition would follow from \cite[Thm.\ 5.4]{labesse} but we need some care since our $G$ is not quasi-split\footnote{Contrary to the assumption on $\mathbf U$ above \cite[Thm.\ 5.4]{labesse} that $U$ is quasi-split at all \emph{inert} places, it seems the assumption ought to be that $\mathbf U$ is quasi-split at \emph{all} finite places. We believe that ``Le second membre \'etant non identiquement nul'' (in the proof of \cite[Thm.\ 5.4]{labesse}, between the second and third displays) is not always true, e.g.\ if $\Pi_w$ is a principal series representation at a non-quasi-split place that splits in $\F$. (See the third paragraph of the current proof.) If it were true, we could deduce Proposition \ref{prop:descent} directly from \cite[Thm.\ 5.4]{labesse}.}; we also need some more information at split places. Thus we sketch the trace formula argument. Again we drop the subscript from $|\LJ|$. 

  The argument of  \cite[Thm.\ 5.4]{labesse} shows the identity (adapted to our notation)
  \begin{equation}\label{character-identity}
  \sum_{\sigma} m(\sigma) \tr \sigma(f) = \tr \Pi (\phi)
  \end{equation}
  with the functions $\phi=\prod_u \phi_u$ on $G^*(\A_{\F})$ and $f=\prod_u f_u$ on $G(\A_{\F^+})$ as in the proof there,
  where the sum runs over automorphic representations $\sigma$ of $G(\A_{\F^+})$ with multiplicity $m(\sigma)$ whose weak base change is $\Pi$.
  The right-hand side is interpreted as the twisted trace by the convention mentioned earlier.

  The key point to show is that the right-hand side does not always vanish. There is a subtlety when $G$ is not quasi-split, because not every test function $\phi$ may be allowed in \eqref{character-identity}. The potential problem is that a base change transfer of $\phi_u$ at $u$ from $G^*(\F_u)$ to $G^*(\F^+_u)$ is not in the image of endoscopic transfer from $G(\F^+_u)$ to $G^*(\F^+_u)$. We make a choice of test functions avoiding this problem.

  At $\infty$ one does the same as in Labesse's proof so that $\tr \Pi_\infty(\phi_\infty)\neq 0$. At finite places $u$, we recall that $f_u$ and $\phi_u$ are related as follows: writing $f^*_u$ for a transfer of $f_u$ from $G(\F^+_u)$ to $G^*(\F^+_u)$, the functions $f^*_u$ and $\phi_u$ are associated in the sense of \cite[4.5]{labesse}.
  There is no problem when $u\notin \Delta^+(G)$ as $G$ and $G^*$ are isomorphic outside $\Delta^+(G)$; more precisely we choose $\phi_u$ on $G(\F\otimes_{\F^+} \F^+_v)$ such that
  $$\tr \Pi_u(\phi_u)\neq 0$$
   and choose $f_u$ to be a base change transfer to $G(\F^+_v)$ (which exists by \cite[Lem.\ 4.1]{labesse}, where it is called an ``associated'' function).
  At each $v=w w^c \in \Delta^+(G)$, choose $f_v$ and let $f^*_v$ be a transfer. 
  Write $\pi^*_v:=\Pi_w$ via the chosen isomorphism $G^*(\F_w) \cong
  G^*(\F^+_v)$. Then by Proposition \ref{prop:JL-char-id},
   $$\tr \pi^*_v(f^*_v)=e(G(\F^+_v))\cdot \tr \left(|\LJ|(\pi^*_v)\right)(f_v).$$
  Note that $|\LJ|(\pi^*_v)$ is irreducible (i.e.\ nonzero) since $\pi^*_v$ is supercuspidal by assumption.
   If we choose $f_v$ such that $\tr \left( |\LJ|(\pi^*_v)\right) (f_v)\neq 0$ then the above identity tells us that
  $\tr \pi^*_v(f^*_v)\neq 0$. Choosing $\phi_v$ to be a function associated with $f^*_v$ (such a $\phi_v$ exists by either \cite[Lem.\ 4.1]{labesse}), we have as in \eqref{BC-LJ-identity},
    $$\tr \Pi_v(\phi_v)= \tr \pi^*_v(f^*_v)\neq 0.$$

  We have exhibited a choice of $f$ and $\phi$ above such that \eqref{character-identity} is valid with the right-hand side non-vanishing. Therefore there exists some $\pi$ on the left-hand side contributing with positive multiplicity. Let $S$ be the set of places of $\F^+$  containing all infinite places and the finite places where $G$ and $\Pi$ are ramified. Write $S_{\mathrm{f}}$ for the subset of finite places in $S$. As we are free to choose $\phi_u$ in the unramified Hecke algebra at each $u\in S_{\mathrm{f}}$, we may assume that $\pi^S$ is unramified with $\BC(\pi^S)=\Pi^S$. The nonvanishing of $\tr \pi_\infty(f_\infty)$ tells us that $\BC(\pi_\infty)=\Pi_\infty$. 
  Thus \eqref{character-identity} is reduced to a formula of the form \eqref{trace-formula-1}, with $\pi$ contributing nontrivially to the sum. Arguing as in the proof of the preceding proposition, we deduce that $|\LJ|(\pi^*_v)\cong \pi_v$.
  \end{proof}

\bibliography{eos}
\bibliographystyle{amsalpha}

\end{document}